\documentclass[11pt]{amsart}

\usepackage{graphicx,epsfig,amssymb}

\usepackage{fullpage}

\chardef\bslash=`\\ % p. 424, TeXbook % Normalized (nonbold, nonitalic) \tt font, to avoid font substitution
\makeatother \hfuzz1pc % Don't bother to report overfull boxes if overage is < 1pc

\newtheorem{theorem}{Theorem}[section]

\newtheorem{corollary}[theorem]{Corollary}

\newtheorem{lemma}[theorem]{Lemma}

\newtheorem{definition}{Definition}[section]
 \newtheorem{remark}[theorem]{Remark}

\theoremstyle{remark} % not used here: makes heading italic

\numberwithin{equation}{section}

\newcounter{list_count}
\newenvironment{rlist}{\begin{list}
{(\roman{list_count})}{\usecounter{list_count}
\setlength{\rightmargin}{\leftmargin}}}{\end{list}}

\def\qed{\hfill \QED\medskip}
\def\epsilon{\varepsilon}
\def\eps{\varepsilon}
\def\a{\alpha}

\def\inter{\bigcap}
\def\union{\bigcup}

\def\comp{\circ}

\def\d{\delta}
\def\wig{\sim}
\newcommand{\overequiv}[1]{\stackrel{#1}{\wig}}

\def\C1+{C^{1+}}

\def\sol{{\mathbf S}}
\def\TT{{\mathbf B}}

\newcommand{\cA}{{\mathcal A}}
\newcommand{\cB}{{\mathcal B}}
\newcommand{\cC}{{\mathcal C}}
\newcommand{\cD}{{\mathcal D}}

\newcommand{\cG}{{\mathcal G}}
\newcommand{\cH}{{\mathcal H}}

\newcommand{\cL}{{\mathcal L}}

\newcommand{\cO}{{\mathcal O}}

\newcommand{\cP}{{\mathcal P}}
\newcommand{\cQ}{{\mathcal Q}}
\newcommand{\cR}{{\mathcal R}}
\newcommand{\cS}{{\mathcal S}}
\newcommand{\cT}{{\mathcal T}}

\def\reals{{\mathbb R}}
\def\integers{{\mathbb Z}}

\newcommand{\be}{\begin{equation}}
\newcommand{\ee}{\end{equation}}

\def\C1{$C^1$}

\def\enddemo{\medskip}

\def\endproof{\medskip
\global\parfillskip0pt plus 1fil\relax
\gdef\enddemo{\medskip}}

\def\qed{\vbox{\hrule\hbox{\vrule height6pt\hskip6pt\vrule}\hrule}}
\def\iint{{\mathrm int}\,}

\def\i{\iota} 
\def\ip{{\iota^\prime}}
\def\su{\{ s,u \}}
\def\s{\sigma}

\def\L{\Lambda}

\def\Str{\cS}

      \def\l{\ell}

\newcommand{\rint}{\operatorname{int}}
 
  \def\Msol{{\rm Msol}} 
 
\def\S{\Sigma} \def\rat{\rho}
\def\II{\hat{I}} \def\Co{{\rm Co}} \def\Gap{{\rm Gap}} \def\msc{{\rm msc}}

\begin{document} \def\Thfin{\Theta_{fin}}
\def\S{Section} \def\s{\sigma} \def\msol{{\rm
msol}} \def\Sol{{\rm scl}}

%May be not necessary: \def\sol{{\rm sol}}
\def\cD{{\mathcal D}} \def\cDU{{\mathcal DU}}
\def\cB{{\mathcal B}}

\title[geometric measures]{Geometric measures
for hyperbolic sets on surfaces} \author{A.
A. Pinto} \author{D. A. Rand}\address[A. A. Pinto]{Faculdade de
Ciencias, Universidade do Porto\\ 4000 Porto,
Portugal.}
%% Note the doubled @@:
\email[A. A. Pinto]{aapinto@@fc.up.pt}

%% Second author \author{D. A. Rand}
\address[D. A. Rand]{Mathematics Institute,
University of Warwick\\ Coventry CV4 7AL,
UK.} %%\curraddr[]{}% when away from home %%Note the doubled @@:
\email[D. A. Rand]{dar@@maths.warwick.ac.uk}
%%\thanks{Research of the second author was supported in part by }

\begin{abstract} We present a moduli space for
all hyperbolic basic sets of diffeomorphisms on surfaces that
have an invariant measure that is absolutely
continuous with respect to Hausdorff measure.
To do this we introduce two new invariants:
the measure solenoid function and the cocycle-gap pair.
We extend the
eigenvalue formula of A. N.  Liv\v sic  and Ja. G. Sinai
for  Anosov diffeomorphisms which preserve an
absolutely continuous measure  to
hyperbolic basic sets on surfaces which possess an invariant measure absolutely
continuous with respect to Hausdorff measure.
We characterise the Lipschitz conjugacy classes
of such hyperbolic systems
 in a number of ways, for example,
in terms of eigenvalues of periodic points and Gibbs measures.
 \end{abstract}

\maketitle

%VT
\thispagestyle{empty}
\def\IMSmarkvadjust{0 pt}
\def\IMSmarkhadjust{0 pt}
\def\IMSmarkhpadding{0 pt}
\def\IMSpubltext{Published in modified form:}
\def\SBIMSMark#1#2#3{
 \font\SBF=cmss10 at 10 true pt
 \font\SBI=cmssi10 at 10 true pt
 \setbox0=\hbox{\SBF \hbox to \IMSmarkhpadding{\relax}
                Stony Brook IMS Preprint \##1}
 \setbox2=\hbox to \wd0{\hfil \SBI #2}
 \setbox4=\hbox to \wd0{\hfil \SBI #3}
 \setbox6=\hbox to \wd0{\hss
             \vbox{\hsize=\wd0 \parskip=0pt \baselineskip=10 true pt
                   \copy0 \break%
                   \copy2 \break% 
                   \copy4 \break}}
 \dimen0=\ht6   \advance\dimen0 by \vsize \advance\dimen0 by 8 true pt
                \advance\dimen0 by -\pagetotal
	        \advance\dimen0 by \IMSmarkvadjust
 \dimen2=\hsize \advance\dimen2 by .25 true in
	        \advance\dimen2 by \IMSmarkhadjust

%
%   Check for publication info
%
%  \newread\jref
  \openin2=publishd.tex
  \ifeof2\setbox0=\hbox to 0pt{}
  \else 
     \setbox0=\hbox to 3.1 true in{
                \vbox to \ht6{\hsize=3 true in \parskip=0pt  \noindent  
                {\SBI \IMSpubltext}\hfil\break
                \input publishd.tex 
                \vfill}}
  \fi
  \closein2
  \ht0=0pt \dp0=0pt
 \ht6=0pt \dp6=0pt
 \setbox8=\vbox to \dimen0{\vfill \hbox to \dimen2{\copy0 \hss \copy6}}
 \ht8=0pt \dp8=0pt \wd8=0pt
 \copy8
 \message{*** Stony Brook IMS Preprint #1, #2. #3 ***}
}

\SBIMSMark{2006/03}{May 2006}{}

\tableofcontents

\newpage

\section{Introduction} \label{dfegrggggd}

We say  that $(f,\Lambda)$ is a $C^{1+}$
\emph{hyperbolic diffeomorphism} if it has
the  following properties:
\begin{rlist}
\item $f:M \to M$ is a  $C^{1+\alpha}$
diffeomorphism  of a compact surface $M$
with respect to a $C^{1+\alpha}$ structure $\cC_f$
on $M$, for some $\alpha >0$;
\item $\Lambda$ is a hyperbolic
invariant subset of $M$; and
\item
$f|\Lambda$ is topologically transitive and
that $\Lambda$ has a local product structure.
\end{rlist}
\def\fL{f, \L}
We denote by $\cT(\fL )$ the set
of all $C^{1+}$ hyperbolic diffeomorphisms
$(g,\Lambda_g)$ such that $(g,\Lambda_g)$ and
$(f,\Lambda)$ are topologically conjugated
by a homeomorphism $h_{f,g}$ (see \S ~\ref{dfgdbbccees}).
We allow both the case
where $\Lambda = M$  and the case where
$\Lambda$ is a proper subset of $M$. If
$\Lambda = M$ then  $f$ is Anosov and $M$ is
a torus \cite{FranksF4,Newhouse}. The best
known examples where $\Lambda$ is a proper
subset of $M$ are the Smale horseshoes, and
the  codimension one attractors  such as  the
Plykin attractor and the derived-Anosov
diffeomorphisms.

For every $g \in \cT(\fL )$, we denote by
$\delta_{g,s}$ (resp.\ $\delta_{g,u}$) the
Hausdorff dimension of the local stable
(resp.\ local unstable) leaves of $g$
intersected with $\L$. Let $\lambda_{g,s}(x)$
and $\lambda_{g,u}(x)$ denote the stable and
unstable eigenvalues of the periodic orbit of $g$
containing a point $x $. A. N.  Liv\v sic  and Ja. G. Sinai
\cite{Livsic} proved that an Anosov
diffeomorphism $g$ has an invariant measure
that is absolutely continuous with respect
to Lesbegue measure if, and only if,
$\lambda_{g,s}(x)\lambda_{g,u}(x)=1$ for
every  periodic point  $x$. We extend
 the theorem of A. N.  Liv\v sic  and Ja. G. Sinai to $C^{1+}$ hyperbolic
diffeomorphisms with hyperbolic sets on
surfaces such as Smale horseshoes and
codimension one attractors.

\begin{theorem} \label{dfsdfaaa1}
A $C^{1+}$
hyperbolic diffeomorphism $g \in \cT(\fL )$
has a $g$-invariant probability measure which
is absolutely continuous to the
Hausdorff measure on $\L_g$
if and only if for every  periodic
point  $x$ of $g|\L_g$,
$$\lambda_{g,s}(x)^{\delta_{g,s}}\lambda_{g,u
}(x)^{\delta_{g,u}}=1 \ . $$
\end{theorem}

The proof of all  theorems stated in the
introduction are given in \S ~\ref{fghtghrr}.

Since  $(f,\Lambda)$ is a $C^{1+}$ hyperbolic
diffeomorphism it admits a Markov partition
$\cR = \{ R_1,\ldots , R_k \}$. This
implies the existence of a two-sided
subshift  $\tau:\Theta \to \Theta$ of finite type $\Theta$ in the
symbol space $\{1 ,\ldots , k\}^\integers$,
and an inclusion $i:\Theta \to \Lambda$ such that
(a) $f \circ i = i \circ \tau$
and (b) $i(\Theta_{j})=R_j$ for
every $j=1,\ldots,k$, where $\Theta_j$ is the cylinder containing all
words $\ldots \epsilon_{-1} \epsilon_0 \epsilon_1 \ldots \in \Theta$ with $\epsilon_0=j$.
For every $g \in \cT(\fL )$,
the inclusion $i_g = h_{f,g}\circ i : \Theta
\to \Lambda_g$ is such that
$g \circ i_g = i_g
\circ \tau$.
We call such a map
$i_g:\Theta \to \Lambda_g$ a \emph{marking}  of
$(g,\Lambda_g)$.

\begin{definition} If $g \in \cT(\fL )$ is a $C^{1+}$
hyperbolic diffeomorphism as above and $\nu$
is a Gibbs measure on $\Theta$ then we say
that $(g,\Lambda_g,\nu)$ is a \emph{Hausdorff
realisation of} $\nu$ if $(i_g)_*\nu$ is absolutely continuous with
respect to the Hausdorff measure on $\Lambda_g$.
If this is the case then we will often just say that $\nu$
is a Hausdorff realisation for $(g,\Lambda_g)$.
\end{definition}

We note that if $g \in \cT(\fL )$ the
Hausdorff measure on $\Lambda_g$ exists and is unique.
However, a  Hausdorff realisation  need
not exist for $(g,\Lambda_g)$.

\def\cTf{\cT_{f,\L}}

 Let   $\cTf ( \delta_s,\delta_u)$  be the set
of all $C^{1+}$ hyperbolic diffeomorphisms
$(g,\L_g)$ in $\cT(\fL )$ such that (i)
$\delta_{g,s}=\delta_s$  and
$\delta_{g,u}=\delta_u$; (ii) there is a
$g$-invariant measure $\mu_g$ on $\Lambda_g$
which is absolutely continuous with respect to the Hausdorff
measure on $\Lambda_g$. We  denote by
$[\nu ]\subset \cTf ( \delta_s,\delta_u)$
the  subset of all
$C^{1+}$-realisations of a Gibbs measure
$\nu$ in $\cTf ( \delta_s,\delta_u)$.

De la
Llave, Marco and Moriyon
\cite{llave,L3,M1,M3} have shown that the set
of stable and unstable eigenvalues of all
periodic points is a complete invariant of
the $C^{1+}$ conjugacy classes of Anosov
diffeomorphisms.
 We extend their result to
the sets  $[\nu]\subset  \cTf ( \delta_s,\delta_u)$.

\begin{theorem} \label{dfsdfaaa3}
(i)  Any two elements of $[\nu ]\subset \cTf ( \delta_s,\delta_u)$
have the same set
of stable and unstable eigenvalues and these
sets are a complete invariant of $[\nu ]$ in the sense that
if $g_1, g_2 \in \cTf ( \delta_s,\delta_u)$ have the same eigenvalues
if, and only if, they are in the same subset $[\nu ]$.

(ii) The map $\nu
\rightarrow [\nu ]\subset \cTf ( \delta_s,\delta_u)$
gives a $1-1$ correspondence between
$C^{1+}$-Hausdorff realisable Gibbs measures
$\nu$ and Lipschitz  conjugacy classes in
$\cTf ( \delta_s,\delta_u)$.
\end{theorem}

In Theorem \ref{fgdgfgdds3}, we also prove
that the set of  stable and unstable
eigenvalues of all periodic orbits of a
$C^{1+}$ hyperbolic diffeomorphism $g \in
\cT(\fL )$ is a complete invariant of each
Lipschitz  conjugacy class. We  note that
for Anosov
diffeomorphisms
every Lipschitz  conjugacy class is a $C^{1+}$
conjugacy class. This can be proved by combining
Remark 8.1 with Lemmas 4.2 and 8.1.

\begin {remark}
We have restricted our discussion to Gibbs measures
because it follows from Theorem \ref{dfsdfaaa3} that,
if $g \in  \cTf( \delta_s,\delta_u)$ has a $g$-invariant
measure $\mu$ which is absolutely continuous with
respect to the Hausdorff measure then $\mu$ is a
$C^{1+}$-Hausdorff realisation of a Gibbs
measure $\nu$ so that $\mu =( i_g)_*
\nu$.
\end{remark}

E. Cawley \cite{cawley} characterised all
$C^{1+}$-Hausdorff realisable Gibbs measures
as Anosov diffeomorphisms using cohomology
classes on the torus. While it is possible that
her cocycles could give enough information to
characterise other hyperbolic systems on
surfaces up to lippeomorphism, it is clear that they
cannot encode enough for $C^{1+}$ conjugacy
because, for example, they do not encode
enough information about gaps and so do not
determine the smooth structure of stable leaves
in the case where they are Cantor sets.
To deal with all these cases in an integrated way
we use measure solenoid functions
and gap-cocycle pairs to
classify
$C^{1+}$-Hausdorff realisable Gibbs measures
of  all $C^{1+}$
hyperbolic diffeomorphisms on surfaces.

The \emph{stable
and unstable  measure solenoid functions} are
built from the Gibbs measure as we
show in  \S ~\ref{fghbjhnnjw}. For
Anosov diffeomorphisms, the domains $\msol^s$
and $\msol^u$ of the stable and unstable
measure solenoid functions are dense subsets
of finite disjoint unions of closed
intervals. Define a stable leaf segment of a
Markov rectangle to be a segment of a
stable leaf crossing the
Markov rectangle (see \S ~\ref{leafsegments},
\S ~\ref{spanning}, \S ~\ref{markgorprp} and
Figure ~\ref{Rectangle}). Every point in
$\msol^s$ consists of a pair $(I,J)$ of
adjacent stable spanning leaf segments of
Markov rectangles which are not contained in
a stable global leaf containing a stable
boundary of a Markov rectangle. The stable
measure solenoid function $\sigma_{\nu,s}$
associates to each pair $(I,J)$ the ratio
between the measure of $J$ and the measure of
$I$ computed with respect to the conditional
measure, determined by the Gibbs measure
$\nu$,  of  a stable leaf containing $I$ and
$J$. The construction of the unstable
solenoid function $\sigma_{\nu,u}$ is
similar. In \S ~\ref{uyggrersx}, we
define a \emph{boundary condition} which
consists of a finite set of simple algebraic
equalities that the continuous extensions of
the stable and unstable measure solenoid
functions have to satisfy, for the
corresponding Gibbs measures to be
$C^{1+}$-Hausdorff realisable in the Anosov
case. This is necessary because in this case the Markov
rectangles have common boundaries along the
stable and unstable  foliations.
We note that in the Anosov case, the Lebesgue measure is the Hausdorff measure.

\begin{theorem}
\label{anosov1}  (Anosov diffeomorphisms)
Suppose that $f$ is a $C^{1+}$ Anosov diffeomorphism
of the torus $\L$. Fix
a Gibbs measure $\nu$ on $\Theta$.
Then the following statements are equivalent:
\begin{rlist}
\item
The set  $\nu$,
$[\nu ]\subset \cTf (1,1)$
is non-empty and is precisely the set
of $g\in \cTf (1,1)$ such that
$(g,\Lambda_g,\nu)$ is a $C^{1+}$ Hausdorff realisation.
In this case $\mu =( i_g)_* \nu$ is absolutely continuous with
respect to Lesbegue measure.
\item
The stable measure solenoid function
$\sigma_{\nu,s}:\msol^s \to \mathbb{R}^+$
has a  non-vanishing H\"older continuous extension to
the closure of $\msol^s$
satisfying the boundary condition.
\item
The unstable measure solenoid function
$\sigma_{\nu,u}:\msol^u \to \Bbb{R}^+$
has a   non-vanishing H\"older continuous extension to
the closure of $\msol^s$
satisfying the boundary condition.
\end{rlist}
\end{theorem}

The treatment of codimension one
attractors has a number of
extra-dificulties due to the fact that the
invariant set $\L$ is locally a Cartesian
product of a Cantor set with an interval but
the stable and unstable measure solenoid
functions are built in a similar way to the
construction for Anosov diffeomorphisms. In the case
of codimension one attractors, the continuous
extension of the  stable measure solenoid
functions have to satisfy the
\emph{cylinder-cylinder condition} for the
corresponding Gibbs measures to be
$C^{1+}$-Hausdorff realisable (see \S
~\ref{fghbjhnnjw}) . The   cylinder-cylinder
condition, like the boundary condition,
consists of a finite set of simple algebraic
equalities and  is needed because
the Markov rectangles have common boundaries
along the stable laminations.
Hence, the cylinder-cylinder condition just applies to the stable measure function.

\begin{theorem}
\label{uyevvgd}
\label{codimension1} (Codimension one attractors)
Suppose that $f$ is a $C^{1+}$ surface diffeomorphism
and $\L$ is a codimension one hyperbolic attractor. Fix
a Gibbs measure $\nu$ on $\Theta$. Then
the following statements are equivalent:
\begin{rlist}
\item
For all $0<\delta_s<1$,
$[\nu ]\subset \cTf(\delta_s ,1)$
is non-empty and is precisely the set
of $g\in \cTf(\delta_s ,1)$ such that
$(g,\Lambda_g,\nu)$ is a $C^{1+}$ Hausdorff realisation.
In this case $\mu =( i_g)_* \nu$ is absolutely continuous with
respect to the Hausdorff measure on $\L_g$.
\item
The stable measure solenoid function
$\sigma_{\nu,s}:\msol^s \to \Bbb{R}^+$
has a non-vanishing H\"older continuous extension to
the closure of $\msol^s$   satisfying the  cylinder-cylinder condition.
\item
The unstable measure solenoid function
$\sigma_{\nu,u}:\msol^u \to \Bbb{R}^+$
has  a non-vanishing H\"older continuous extension
to the closure of $\msol^u$.
\end{rlist}
\end{theorem}

In the case of Smale horseshoes, there are no
extra conditions that the measure solenoid
functions have to satisfy for the
corresponding Gibbs measures to be
$C^{1+}$-Hausdorff realisable.

\begin{theorem}
\label{htrytgvg} (Smale horseshoes)
Suppose that $(f,\L)$ is a Smale horseshoe
and $\nu$ is a Gibbs measure on $\Theta$.
Then
for all $0<\delta_s , \delta_u<1$,
$[\nu ]\subset \cTf(\delta_s ,\delta_u)$
is non-empty and is precisely the set
of $g\in \cTf(\delta_s ,\delta_u)$ such that
$(g,\Lambda_g,\nu)$ is a $C^{1+}$ Hausdorff realisation.
In this case $\mu =( i_g)_* \nu$ is absolutely continuous with
respect to the Hausdorff measure on $\L_g$.
\end{theorem}

Using Theorems \ref{dfsdfaaa3}-\ref{uyevvgd},
for $\i \in \{ s,u\}$, we prove that
the map $\nu \to \sigma_{\nu,\i}$   gives a
one-to-one correspondence between the subsets
$[\nu]\subset \cTf(\delta_s ,\delta_u)$ and the
measure solenoid functions $\sigma_{g,\i}$
satisfying the conditions indicated in the
above Theorems \ref{anosov1} and
\ref{uyevvgd}. By Theorem \ref{dfsdfaaa3},
we have that the sets $[\nu]\subset \cTf(\delta_s,\delta_u)$
are precisely the Lipschitz conjugacy classes
contained in $\cTf(\delta_s ,\delta_u)$.

\begin{theorem}
\label{333dfsf}
The measure solenoid functions
determine a pair of   infinite-dimensional metric
spaces, that we denote by ${\mathcal SOL}^s$ and ${\mathcal SOL}^u$,
which parametrize all Lipschitz conjugacy
classes $[\nu]\subset \cTf(\delta_s ,\delta_u)$.
\end{theorem}

The  scaling functions, presented by M.
Feigenbaum in \cite{Feigenbaum1,Feigenbaum3}
and D. Sullivan in \cite{sullivan:scalingfn},
and the  solenoid functions, presented  in
\cite{PR1,HR,PS} (see also \S
~\ref{fggghhreew} and \S ~\ref{fgvdghbvc}), are
used to classify all $C^{1+}$ conjugacy
classes of expanding maps on train-tracks and
of hyperbolic diffeomorphisms for  a given a
topological conjugacy class. Both scaling function and
solenoid function are complete invariants of the smooth structure
but the solenoid functions have the great advantage that, unlike
the scaling functions, one knows which solenoid functions occur
as the solenoid functions of $C^{1+}$ expanding maps.

In  \S
~\ref{dfgggvvd}, we introduce the \emph{stable
and unstable cocycle-gap pairs}
$(\gamma_s,J_s)$ and $(\gamma_u,J_u)$ which
allow us to parametrize the $C^{1+}$
conjugacy classes inside of each
Lipschitz conjugacy
class $[\nu]\subset \cTf(\delta_s ,\delta_u)$.
The cocycle-gap pair  $(\gamma,J)$ consists of
 a   gap function  $\gamma$ and a
measure-length ratio cocycle function $J$
(see Definitions \ref{cocycledddd} and
~\ref{gapratiosds}). The domain of a stable
gap function is the set of all pairs
$(\xi_1, \xi_2)$ with the following
properties: (a) the stable leaf segments
$\xi_1$ and $\xi_2$ intersect the invariant
set $\L$ just in their end points;
(b) there is a stable leaf  $K$ of  a Markov rectangle
which contains $f \xi_1$ and $f \xi_2$.
The \emph{stable gap function} $\gamma$ is a
H\"older continuous function satisfying the
following algebraic equality $\gamma
(\xi_{1}:\xi_{2}) = \gamma
(\xi_{1}:\xi_{3})\gamma (\xi_{3}:\xi_{2})$.
(We use the notation $\xi_{1}:\xi_{2}$ rather than
$\xi_{1},\xi_{2}$ just to emphasise that
$\gamma$ will be measuring ratios.)
The domain of a stable measure-length ratio
cocycle is the set of all stable spanning
leaves of Markov rectangles. A \emph{stable
measure-length ratio cocycle} is a function
$J=\kappa/(\kappa\comp f)$ where  $\kappa$
is a positive H\"older continuous function
satisfying a set of algebraic inequallties
given by \eqref{aazzqq3}. The unstable gap
functions and measure-length
ratio cocycle are defined similarly.
In the case of Smale horseshoes the domains
of the gap  functions and  of the
measure-length ratio cocycle functions are
Cantor sets. In the case of codimension one
attractors the domains of  the  stable gap
functions and  of the stable measure-length
ratio cocycle functions are a finite union of
intervals, and the stable cocycle-gap pairs
have to satisfy the \emph{cocycle-gap
property} (see Definiton \ref{ghhhhgdd} and Lemma \ref{ggdgdhhbna}) which is due  to
the fact that the Markov rectangles have
common boundaries along the stable
laminations. The stable and unstable
cocycle-gap   pairs  give rise to an
infinite dimensional metric  space, that we
denote respectively by ${\mathcal
JG}^s(\nu,\delta_s)$ and ${\mathcal
JG}^u(\nu,\delta_u)$.
\begin{theorem}
\label{gthrthrtr}
(i) \emph{(Smale horseshoes)} There is a natural  map
$$g \to
(\gamma_s(g),J_s(g),\gamma_u(g),J_u(g))$$
which gives a one-to-one correspondence
between $C^{1+}$ conjugacy classes of
diffeomorphisms $g$ contained in
$[\nu]\subset \cTf(\delta_s ,\delta_u)$ and stable and unstable
cocycle-gap pairs contained in  ${\mathcal
JG}^s(\nu,\delta_s) \times {\mathcal
JG}^u(\nu,\delta_u)$.

(ii)
(Codimension one attractors)
There is a natural map $g \to
(\gamma_s(g),J_s(g))$ which gives a
one-to-one correspondence between $C^{1+}$
conjugacy classes of diffeomorphisms $g$
contained in $[\nu ]\subset \cTf(\delta_s ,1)$
and stable cocycle-gap pairs contained in
${\mathcal JG}^s(\nu,\delta_s)$.
\end{theorem}

\subsection{Self-renormalisable structures}

\begin{figure}[htbp]
\includegraphics[width=11cm]{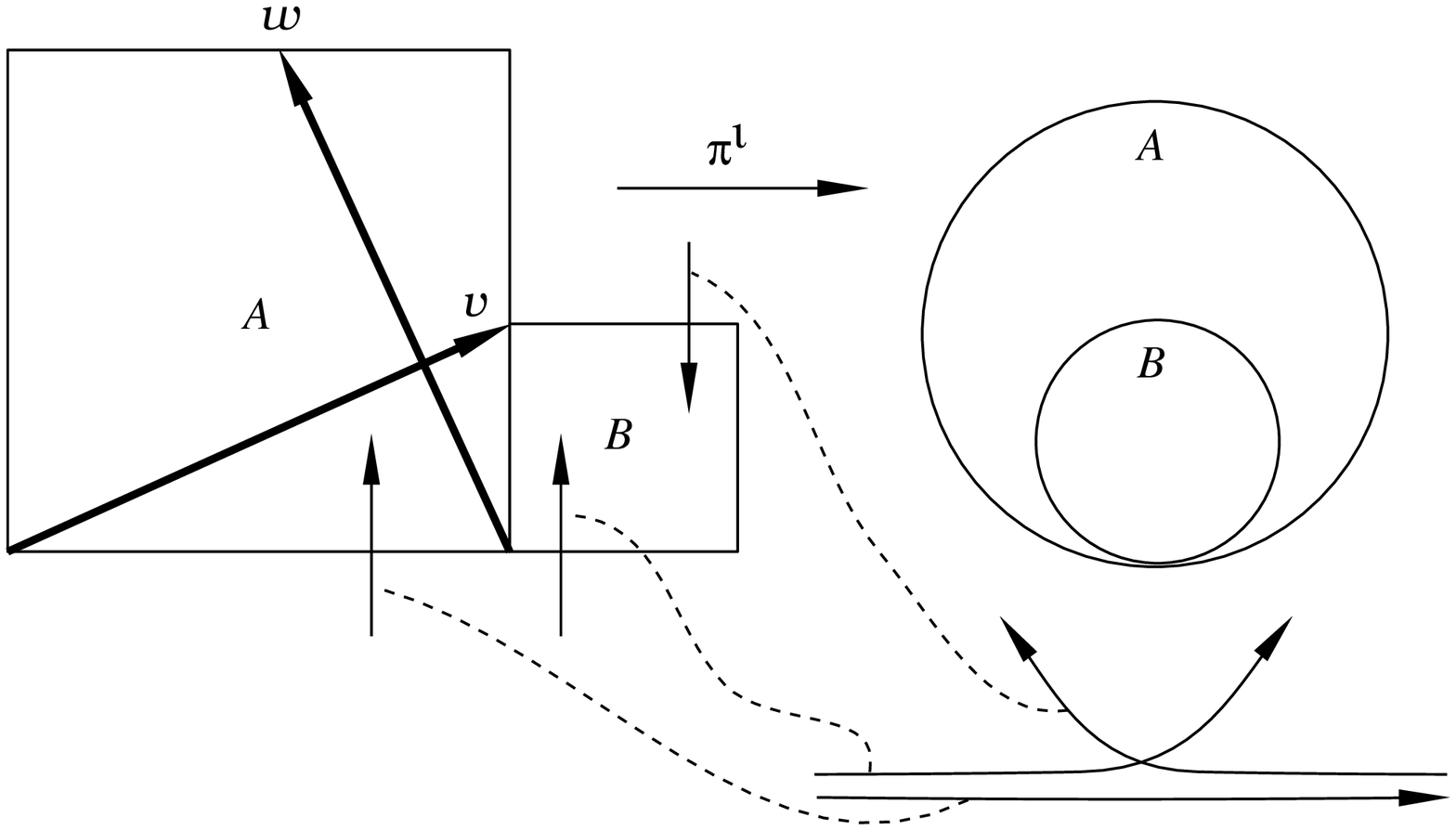}
\caption{This figure illustrates a (unstable)
train track for the Anosov map
$g: \reals^2 \setminus (\integers v \times \integers w) \to
\reals^2 \setminus (\integers v \times \integers w)$ defined by
$g(x,y) = (x+y,y)$.
The rectangles $A$ and $B$ are the Markov rectangles and the
vertical arrows show paths along unstable manifolds from
$A$ to $A$ and from $B$ to $A$.
The train track is represented by the pair of circles and the curves
below it show the smooth paths through the junction of the
two circles which arise from the smooth paths between the rectangles
$A$ and $B$ along unstable manifolds.
Note that there is no smooth path from $B$ to $B$ even though in this
representation of the train track it looks as though there ought to be.
This is because there is no unstable manifold running directly from
the rectangle $B$ to itself.}
\label{Train_tracksAnosov}
\end{figure}

\begin{figure}[tbp]
\includegraphics[width=11cm]{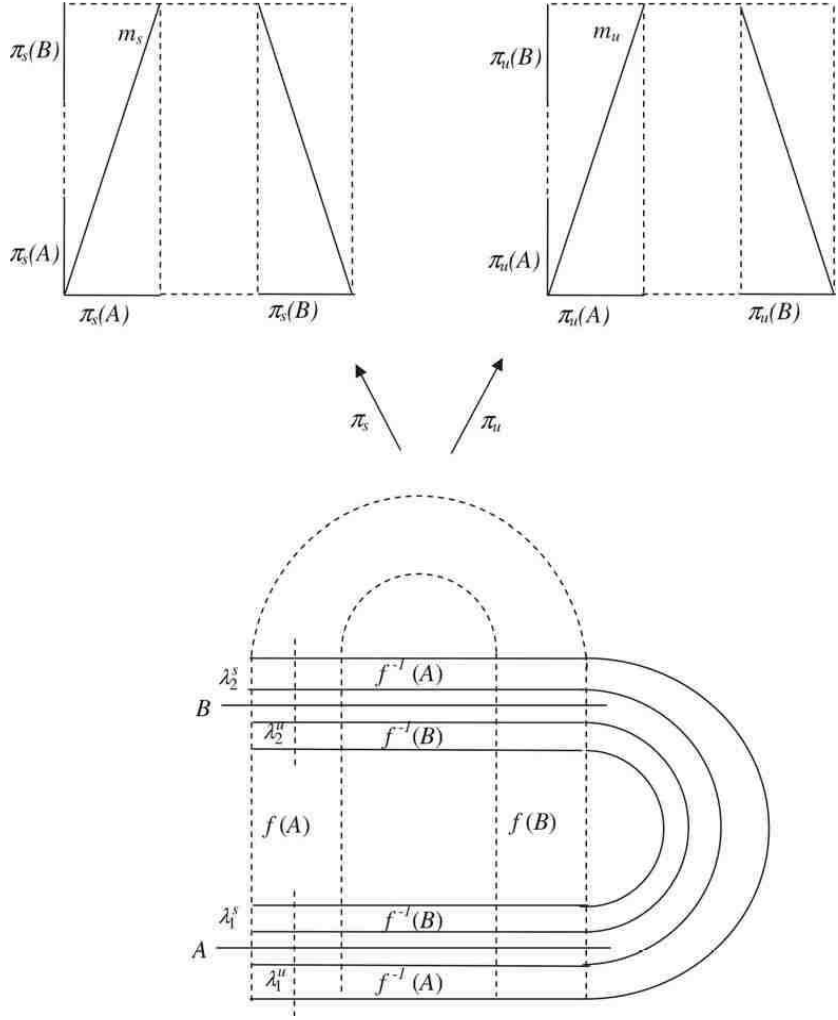}
\caption{\label{example0002}
A Markov partition for the Smale-shoe $f$ into two rectangles
$A$ and $B$.
A representation of the  Markov maps
$m_s:\Theta^s \to \Theta^s $  and
$m_u:\Theta^u \to \Theta^u $
for  Smale horseshoes.}
\label{ref:Smale}
\end{figure}

\begin{figure}[tbp]
\includegraphics[width=17cm]{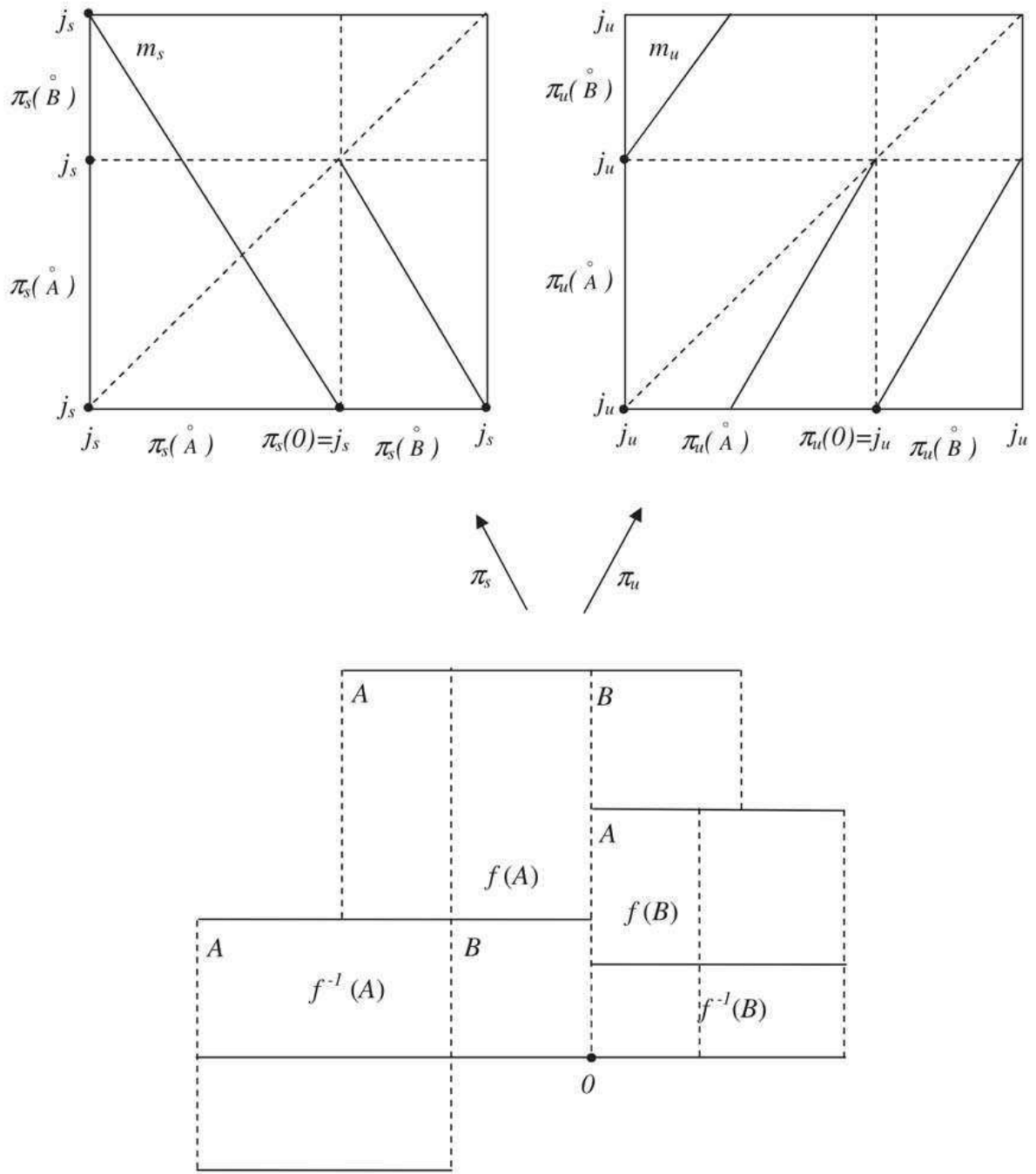}
\caption{\label{example03}
A representation of the  Markov maps
$m_s:\Theta^s \to \Theta^s $  and
$m_u:\Theta^u \to \Theta^u $
as  maps  of the interval for Anosov diffeomorphisms.}
\end{figure}

In \S ~\ref{bgtrewss}, we construct
\emph{$C^{1+}$  stable and unstable  self-renormalisable
structures on  train-tracks}. The train-tracks  are a form
of optimal local leaf-quotient
space of the stable and unstable
laminations of $\Lambda$.  Locally,
these train-tracks are just the
quotient space of stable or
unstable leaves within a Markov
rectangle, but globally the
identification of leaves common to
two more than one rectangle gives a
non-trivial structure and
introduces junctions.
They are
characterised by being the compact
quotient on which the Markov map
induced by the action of $f$ is
continuous with the minimal number
of identifications.  A smooth
structure on the stable or unstable
leaves of $\Lambda$ induces a smooth
structure on the corresponding
train-tracks and \emph{vice-versa}.
Now, we use that  the
holonomies of codimension one
hyperbolic systems are $C^{1+}$ (see \cite{Holonomy}),
and so the holonomies also  project in the train-tracks
and together with the Markov maps
give rise to what we call self-renormalisable $C^{1+}$ structures
(see Figures  \ref{Train_tracksAnosov}, \ref{example0002} and \ref{example03}).

\begin{remark}
These structures are called self-renormalisable
because the train track $\TT$ has defined on it both a Markov
map $m$ and a pseudo group of holonomy maps $\{ h\}$.
For the train tracks arising from Anosov diffeomorphisms
of the torus  $\{ h\}$ can be identified with a
$C^{1+}$ diffeomorphism $g$ of the circle and the Markov map $m$
defines a renormalisation of $g$ which agrees with the usual one
(\cite{orss,rand:survey}). The renormalised map $Rg$ is $C^{1+}$conjugate
to $g$. Thus these $C^{1+}$ self-renormalisable structures are in one-to-one
correspondance with $C^{1+}$ fixed points of the the circle
map renormalisation operator. It then follows from the results of
this paper that these fixed points are in one-to-one
correspondance with those Anosov diffeomorphisms of the the
torus which preserve an a measure that is absolutely continuous
with respect to Lebesgue measure. For train  tracks arising from non-Anosov diffeomorphisms one gets other interesting renormalisation structures, for example,
for interval exchange maps.
\end{remark}

In \S ~\ref{dfghnnvvfr}, we prove
the following  useful   equivalence between
2-dimensional dynamics and 1-dimensional dynamics.

\begin{theorem}
\label{pfiejd}
There is a natural map  $g \rightarrow
({\cS}_s(g),\cS_u(g))$ which  gives a
one-to-one correspondence between $C^{1+}$
conjugacy classes  in $\cT(\fL )$ and pairs
of    stable and unstable $C^{1+}$
self-renormalisable structures.
\end{theorem}

Hence,  for a pair $(\cS_s,\cS_u)$ of
$C^{1+}$ self-renormalisable structures to be
 realisable by a  $C^{1+}$ hyperbolic
diffeomorphism in $\cT(\fL )$, the unstable
$C^{1+}$ self-renormalisable structure does
not impose any restriction in the stable
$C^{1+}$ self-renormalisable structure, and
\emph{vice-versa}. The same is no longer true
 if we ask $g \in \cT(\fL )$ to be a
$C^{1+}$-Hausdorff realisation of a Gibbs
measure as we describe in the next section.

\subsection{Realisation of Gibbs measures on train-tracks}
\label{fdgsdgsg}
We are going to study the $C^{1+}$-Hausdorff
realisations of   Gibbs measures as
self-renormalisable structures. Then we use
this information to study the
$C^{1+}$-Hausdorff realisations of  Gibbs
measures as  $C^{1+}$ hyperbolic
diffeomorphisms, going like this from
one-dimensional dynamics to two-dimensional
dynamics.

%Possibe figure here figfig

Let $\pi_u: \Theta \to \Theta^u$ and  $\pi_s:
\Theta \to \Theta^s$ be respectively the
natural projections on the left and right
infinite words. Let  $\tau_u:\Theta^u \to
\Theta^u$ and $\tau_s:\Theta^s \to \Theta^s$
be the corresponding right and left shifts.
Let  $\i \in\{s,u\}$ where $s$ denotes stable
and $u$ denotes unstable. Let us denote by
$\TT^\i$ the $\i$ train-track and  by $i_\i:
\Theta^\i \to \TT^\i$ the natural marking
induced by $i:\Theta \to \L$ (see \S
~\ref{fhhhffgf}). A $\tau$-invariant measure
$\nu$ on $\Theta$ determines a unique
$\tau_\i$-invariant measure
$\nu_\i=(\pi_\i)_* \nu$ on $\Theta^\i$.
Conversely, a $\tau_\i$-invariant measure
$\nu_\i$ on $\Theta^\i$ has an unique
natural extension to a $\tau$-invariant
measure $\nu$ on $\Theta$. We say that a
$\tau_\i$-invariant measure $\nu_\i$ is a
\emph{Gibbs measure} if its natural extension
$\nu$ is a Gibbs measure on $\Theta$. A
$C^{1+}$ $\i$ self-renormalisable structure
$\cS_\i$ is a \emph{$C^{1+}$-Hausdorff
realisation of a   Gibbs measure $\nu$ on
$\Theta$} if, for  every chart $(e,U)$ of the
self-renormalisable structure $\cS_\i$, the
pusforward $ (e  \circ i_\i)_* \nu_\i$ of the
measure $\nu_\i$ is a measure absolutely
continuous with respect to the Hausdorff
measure of the set $e(U)$. Let us denote by
$\delta (\cS_\i)$ the Hausdorff dimension of
the set $e(U)$ for every chart $(e,U)$  of
the smooth structure of $\cS_\i$. Clearly,
this is independent of the chart $(e,U)$. We
denote by  $\cD^\i(\nu,\delta)$   the set of
all $C^{1+}$ $\i$ self-renormalisable
structures $\cS_{\i}$ with $\delta
(\cS_\i)=\delta$ and which are
$C^{1+}$-Hausdorff realisations of the  Gibbs
measure $\nu$ on $\Theta$. By Theorem
\ref{thm:SRB},  every $C^{1+}$ $\i$
self-renormalisable structure $\cS_\i$ is a
$C^{1+}$-Hausdorff realisation of an unique
Gibbs measure $\nu_{\cS_\i}$ on $\Theta$.

\begin{theorem}
\label{rfgrghhbvcde}
The  map $g \rightarrow
({\cS}_s(g),\cS_u(g))$ gives a 1-1
correspondence between $C^{1+}$    conjugacy
classes  in $[\nu]\subset\cTf ( \delta_s,\delta_u)$ and
pairs  in $\cD^s(\nu,\delta_s) \times
\cD^u(\nu,\delta_u)$.
\end{theorem}

Hence,  if $g \in \cTf (\delta_s,\delta_u)$
then $\delta (\cS_s(g))=\delta_s$ and $\delta
(\cS_u(g))=\delta_u$. Let $\cS_\i$ be a
$C^{1+}$ $\i$ self-renormalisable structure.
If $\delta (\cS_\i)=1$   we    call $\TT^\i$
a \emph{no-gap train-track}. If $0<\delta
(\cS_\i)<1$   we call $\TT^\i$  a \emph{gap
train-track}. Let  $\ip$  denote the element
of $\su$ which is not $\i \in \su$.

\begin{theorem}
Let $\TT^s$ and $\TT^u$ be the stable and unstable train-tracks
determined by a $C^{1+}$ hyperbolic diffeomorphism $(f,\L)$.
\label{self1}
The set $\cD^\i(\nu,\delta_\i)$ is  non-empty if, and ony if,
the    $\i$-measure solenoid function
$\sigma_\nu:\msol^\i \to \reals^+$ of the Gibbs measure $\nu$
has the following properties:
\begin{rlist}
\item If $\TT^\i$ and $\TT^\ip$ are no-gap
train-tracks then  $\sigma_{\nu}$ has a
non-vanishing H\"older continuous extension
to the closure of $\msol^\i$ satisfying the
boundary condition.
\item  If $\TT^\i$ is a no-gap train-track
and $\TT^\ip$ is a gap train-track  then
$\sigma_{\nu}$ has a   non-vanishing H\"older
continuous extension to the closure of
$\msol^\i$.
\item  If $\TT^\i$ is a gap train-track
and $\TT^\ip$ is a  no-gap train-track  then
$\sigma_{\nu}$ has a   non-vanishing H\"older
continuous extension to the closure of
$\msol^\i$ satisfying the   cylinder-cylinder
condition. \item  If $\TT^\i$ and $\TT^\ip$
are gap train-tracks then $\sigma_\nu$ does
not have to satisfy any extra-condition.
\end{rlist}
Furthermore, $\cD^\i(\nu,\delta_\i) \ne
\emptyset$ if, and only if,
$\cD^\ip(\nu,\delta_\ip) \ne \emptyset$
\end{theorem}

By Theorem \ref{83dfsf}, the set of all
$\i$-measure solenoid functions $\s_\nu$ with
the properties indicated in Theorem
\ref{self1} determine an  infinite
dimensional metric  space ${\mathcal SOL}^\i$
which gives a nice parametrization of all
Lipschitz conjugacy classes
$\cD^\i(\nu,\delta)$ of    $C^{1+}$
self-renormalisable structures $\cS_{\i}$
with a given Hausdorff dimension
$\delta(\cS_{\i})=\delta$. We have used the
same notation ${\mathcal SOL}^\i$ as in \S ~\ref{dfegrggggd} above
because we will show that they are effectively the same sets.

\begin{theorem}
\label{dfgrggrrd}
Let us suppose that $\cD^\i(\nu,\delta) \ne \emptyset$.
\begin{rlist}
\item
(Flexibility) If $\TT^\i$ is a gap
train-track then $\cD^\i(\nu,\delta)$ is   an
infinite dimensional space parametrized by
cocycle-gap pairs contained in ${\mathcal
JG}^\i(\nu,\delta)$.
\item
(Rigidity) If $\TT^\i$ is a no-gap train-track
then $\cD^\i(\nu,1)$ consists of a single
$C^{1+}$ self-renormalisable structure.
\end{rlist}
\end{theorem}

By Lemma \ref{akju1}, each set $\cD^\i(\nu,\delta)$ is either
empty or a Lipschitz conjugacy class. Hence, if
$\TT^\i$ is a no-gap train-track then the
Lipschitz conjugacy class consists of  a
single   $C^{1+}$ self-renormalisable structure.
Furthermore, by Lemma \ref{akju2},  the set of   eigenvalues of all
periodic orbits of $\cS_\i$
is a complete invariant of
each set    $\cD^\i(\nu,\delta)$ (see also \cite{sullivan:nested}).

\begin{theorem}
\label{fdgdrgrrqs}
(Rigidity) If $\delta_\i=1$, the mapping $g
\rightarrow \cS_\ip(g)$ gives a 1-1
correspondence between $C^{1+}$  conjugacy
classes  in $[\nu ]\subset\cTf ( \delta_s,\delta_u)$ and
 $C^{1+}$   self-renormalisable structures in
$\cD^\ip (\nu,\delta_\ip)$.
\end{theorem}

Hence, the map $g \rightarrow \cS_\ip(g)$
gives rise to the following one-to-one
correspondences. For Anosov diffeomorphisms
$(f,\L)$, there is an one-to-one
correspondence between (i) $C^{1+}$ conjugacy
classes of hyperbolic diffeomorphisms    $g
\in \cT(\fL )$ which are $C^{1+}$-Hausdorff
realisations of a Gibbs measure; (ii)
$C^{1+}$ unstable self-renormalisable
structures; and (iii) $C^{1+}$ stable
self-renormalisable structures. For
codimension one attractors  $(f,\L)$, there
is a one-to-one  correspondence between (i)
$C^{1+}$ conjugacy classes of hyperbolic
diffeomorphisms    $g \in \cT(\fL )$ which
are $C^{1+}$-Hausdorff realisations of a
Gibbs measure; and (ii)  $C^{1+}$ stable
self-renormalisable structures.

\subsection{Solenoid functions and  Gibbs measures}
\label{fggghhreew}

In \S ~\ref{fgvdghbvc}, we introduce
\emph{HR-structures} (HR for
H\"older-ratios).  These associate an affine
structure to each stable and unstable leaf
segment in such a way that these vary
H\"older continuously with the leaf and are
kept invariant by the dynamics of $f$. As we
will describe the HR-structures are in
one-to-one correspondence with  the $C^{1+}$
conjugacy classes of $C^{1+}$ hyperbolic
diffeomorphisms in $\cT(\fL )$. The
HR-structures are characterized by a pair
$(r_s,r_u)$ of ratio functions. The ratio
functions $r_s$ and $r_u$ are independent one
from another. However, if we ask that the
$C^{1+}$ hyperbolic diffeomorphism have an
associated geometric measure then we will see
later that the ratio function $r_u$ imposes
restrictions on the ratio function $r_s$ and
vice-versa.

For $\i\in \{ s,u\}$,
we construct  a   topological set  $\sol^\i$
which is either
isomorphic to a finite union of closed intervals on the real line
or to  an  embedded Cantor set on the real line.
The ratio functions $r_\i$ when restricted
to the set   $\sol^\i$ are H\"older continuous functions
$\sigma_\i=r_\i|\sol^\i$   completely
characterised by a finite set of  properties that we explain also in \S
~\ref{fgvdghbvc}.
We call the functions $\sigma_\i$  \emph{solenoid functions}. They
form an infinite dimensional metric space
which parametrizes the $C^{1+}$ conjugacy
classes of $C^{1+}$ hyperbolic diffeomorphisms contained in $\cT(\fL )$
(see Corollary \ref{cor:sol}).

We say that  $D$ is a \emph{stable leaf $n$-cylinder segment of $\L$} if $f^{-n}D$ is
a
stable spanning leaf segment of a Markov rectangle. Let $E_D$ be a
stable spanning leaf segment of a Markov rectangle such that $D \subset E_D$.
Let $\nu$ be a Gibbs measure and $\mu=i_* \nu$ the corresponding invariant
measure in $\L$.
Let $\rho(D\cap \L, E_D \cap \L)$
be the ratio between the measure of $D\cap \L$
and the measure of $E_D\cap\L$
with respect to the conditional measure of $\mu$ in $E_D$.
The \emph{$\delta$-stable bounded solenoid class of a
Gibbs measure $\nu$} is the set of
all stable solenoid functions $\s$
with the following property:
There is $C=C(\s)>0$ such that for
every  stable leaf $n$-cylinder segment
of $\L$
we have
$$
\left|\delta_\i \log r (D\cap \L, E_D \cap \L) -\log  \rho(D\cap \L, E_D \cap
\L)  \right|
< C   \  ,
$$
where $r$ is the stable ratio function
determined by the stable solenoid
function $\s$
(see also Definition \ref{dfgdrfgggx}).
The $\delta$-unstable bounded solenoid
class of a Gibbs measure $\nu$ is
defined in a similar way.

\begin{theorem}
\label{ddfsffe1}
\begin{rlist}
\item
There is a natural map $g\to
(\s_s(g),\s_u(g))$ which gives a one-to-one
correspondence between $C^{1+}$ conjugacy
classes of $C^{1+}$ hyperbolic
diffeomorphisms  $g \in \cT(\nu,
\delta_s,\delta_u)$ and pairs
$(\s_s(g),\s_u(g))$  of stable and unstable
solenoid functions such that,  for   $\i$
equal to $s$ and $u$, $\s_\i(g)$ is contained
in the $\delta_\i$-bounded solenoid
equivalence class of $\nu$.
\item
There is a natural map $\cS_\i  \to
\s_{\cS_\i}$ which gives a one-to-one
correspondence between $C^{1+}$
self-renormalisable  structures $\cS_\i$
contained in  $\cD^\i(\nu,\delta_\i)$ and
$\i$-solenoid functions $\s_{\cS_\i}$
contained in the $\delta_\i$-bounded
equivalence class of $\nu$.
\item
Let us suppose that $\cD^\i(\nu,\delta_\i) \ne\emptyset$.
\begin{rlist}
\item[(a)]
(Rigidity) If $\delta_\i=1$ then the
$\delta_\i$-bounded solenoid equivalence
class of   $\nu$ is a singleton consisting in
the  continuous extension of the $\i$ measure
solenoid function $\sigma_{\nu,\i}$ to
$\sol^\i$.
\item[(b)]
(Flexibility) If $0<\delta_\i<1$ then the
$\delta_\i$-bounded solenoid equivalence
class of  $\nu$  is an infinite dimensional
space of  solenoid functions.
\end{rlist}
\end{rlist}
\end{theorem}

In \S ~\ref{dfgggvvd}, we use the cocycle-gap
pairs to construct explicitly the solenoid
functions in the $\delta_\i$-bounded
equivalence classes of the $C^{1+}$-Hausdorff
realisable Gibbs measures $\nu$.

By Lemma \ref{fghfhttf}, given an
$\i$-solenoid function $\s_\i$ there is a
unique $C^{1+}$ self-renormalisable structure
$\cS_\i$ such that $\s_\i=\s_{\cS_\i}$ and,
by Theorem \ref{thm:SRB} and Lemma
\ref{cor_geommeas}, there is a unique
$C^{1+}$-Hausdorff realisable Gibbs measure
$\nu = \nu_{\s_\i}$ such that $\cS_\i \in
\cD^\i(\nu,\delta_\i)$ with $\delta_\i =
\delta(\cS_\i)$. Hence, by Theorem
\ref{ddfsffe1} (ii), given an $\i$-solenoid
function $\s_\i$ there is an unique
$C^{1+}$-Hausdorff realisable Gibbs measure
$\nu$ such that $\sigma_\i$ belongs to the
$\delta_\i$-bounded solenoid equivalence
class of  $\nu$.

\begin{theorem}
\label{ddfsffe3}
Given an $\i$-solenoid function $\s_\i$ and
$0<\delta_\ip \leq 1$, there is a unique Gibbs measure
$\nu$ and a unique $\delta_\ip$-bounded
equivalence class of $\nu$ consisting of
$\ip$-solenoid functions $\s_\ip$ such that
the $C^{1+}$ conjugacy class of hyperbolic
diffeomorphisms $g \in  \cTf (\delta_s,\delta_u)$
determined by the pair
$(\s_s,\s_u)$ have an
invariant measure $\mu=(i_g)_* \nu$ absolutely
continuous with respect to the Hausdorff
measure.
\end{theorem}

Putting together Theorem \ref{ddfsffe1} and
Theorem \ref{ddfsffe3}, we obtain the
following implications:
\begin{rlist}
\item
(Flexibility   for Smale horseshoes) For
$\i=s$ and $u$, given a $\i$-solenoid
function $\s_\i$ there is an infinite
dimensional space of solenoid functions
$\s_\ip$ such that the $C^{1+}$  hyperbolic
Smale horseshoes determined by the pairs
$(\s_s,\s_u)$ have  an
invariant measure $\mu$ absolutely continuous
with respect to the Hausdorff measure.
\item
(Rigidity     for Anosov diffeomorphisms) For
$\i=s$ and $u$,   given an $\i$-solenoid
function $\s_\i$ there is an unique
$\ip$-solenoid function such that the
$C^{1+}$ Anosov diffeomorphisms determined by
the pair $(\s_s,\s_u)$ has an invariant measure $\mu$
absolutely continuous with respect to
Lebesgue.
\item
(Flexibility for codimension one attractors)
Given an unstable solenoid function $\s_u$ there
is an infinite dimensional space of
stable solenoid functions $\s_s$ such that the
$C^{1+}$ hyperbolic codimension one
attractors determined by the pairs
$(\s_s,\s_u)$ have an
invariant measure $\mu$ absolutely continuous
with respect to the Hausdorff measure.
\item
(Rigidity for codimension one attractors)
Given an $s$-solenoid function $\s_s$ there
is an unique unstable solenoid function $\s_u$
such that the $C^{1+}$   hyperbolic
codimension one attractors determined by the
pair $(\s_s,\s_u)$ have  an invariant measure $\mu$ absolutely
continuous with respect to the Hausdorff
measure using non-zero stable and unstable
pressures.
\end{rlist}

In this paper we prove, in fact, a
generalized version of the results presented
in the introduction. The theory that we
develop here studies the properties and
classifies all the Gibbs measures which are
\emph{$C^{1+}$-realisable as  natural
geometric  measures of $C^{1+}$
self-renormalisable structures} (see
Definition \ref{fgtyjtyjyu1}) and of \emph{$C^{1+}$
hyperbolic diffeomorphisms}  (see  Definition
\ref{fgtyjtyjyu}). The set of natural geometric
measures contains the invariant measures which are absolutely
continuous with respect to Hausdorff
measure.

The proof of all  theorems stated in the
introduction are given in \S ~\ref{fghtghrr}.

\section{Hyperbolic diffeomorphisms}

In this section, we present some
basic facts on hyperbolic dynamics,
that we include for clarity of the
exposition.

\subsection{Interval notation.} We
also use the notation of interval
arithmetic for some inequalities
where: \begin{rlist} \item if $I$
and $J$ are intervals then $I+J$,
$I .  J$ and $I/J$ have the obvious
meaning as intervals, \item if
$I=\{ x\}$ then we often denote $I$
by $x$, and \item $I\pm \epsilon$
denotes the interval consisting of
those $x$ such that
$|x-y|<\epsilon$ for all $y\in I$.
\end{rlist}
By $\phi (n) \in
1\pm {\mathcal O}(\nu^n)$ we mean that
there exists a constant $c>0$
depending only upon explicitely
mentioned quantities such that for
all $n\ge 0$, $1-c\nu^n<\phi
(n)<1+c\nu^n$.
By  $\phi (n) = {\mathcal O}(\nu^n)$ we mean  that
there exists a constant $c \ge 1$
depending only upon explicitely
mentioned quantities such that for
all $n\ge 0$,  $c^{-1} \nu^n \le \phi
(n) \le c\nu^n$.

\subsection{Stable and unstable
superscripts}
\label{superscripts}

Throughout the paper
we will use the following notation:
we use $\i$ to denote an element of
the set $\su$ of the stable and
unstable superscripts and $\ip$ to
denote the element of $\su$ that is
not $\i$.  In the main discussion
we will often refer to objects
which are qualified by $\i$ such
as, for example, an $\i$-leaf:
This means a leaf which is a leaf
of the stable lamination if $\i=s$,
or a leaf of the unstable lamination if
$\i=u$.  In general the meaning
should be quite clear.

We define the map $f_\i=f$
if $\i=u$ or $f_\i=f^{-1}$ if
$\i=s$.

\subsection{Leaf segments}
\label{leafsegments}

Let $d$ be a metric on $M$.
For $\iota \in \{s,u\}$,
if $x\in\Lambda$ we denote the local
$\iota$-manifolds
through $x$ by
$$
W^{\iota }(x,\varepsilon )=\left\{ y\in M:d(f_{\iota
}^{-n}(x),f_{\iota
}^{-n}(y))\leq \varepsilon
,~\mbox{for all}~n\ge 0\right\} .
$$
By the Stable Manifold Theorem (see
\cite{HirschPugh}), these sets are
respectively contained in the stable and
unstable immersed manifolds
$$
W^{\iota }(x)=\bigcup_{n\ge 0}f_{\iota }^{n}\left( W^{\iota }\left(
f_{\iota }^{-n}(x),\varepsilon _{0}\right) \right)
$$
which are the image of a $C^{1+\gamma}$
immersion $\kappa_{\iota,x}:\reals \to M$. An
{\em open (resp. closed) full} $\iota $-{\em
leaf segment $I$} is defined as a subset of
$W^{\iota }(x)$ of the form $\kappa_{\iota,x}
(I_{1})$ where $I_{1}$ is a non-empty open (resp.
closed) subinterval in $\reals$.
An {\em open} (resp. {\em closed) $\iota$-leaf
segment} is the intersection with $\Lambda$
of an open (resp. closed) full $\iota$-leaf
segment such that the intersection contains
at least two distinct points. If the
intersection is exactly two points we call
this  closed $\iota$-leaf segment an  {\em
$\iota$-leaf gap}. A {\em full $\iota$-leaf
segment} is either an open or closed full
$\iota$-leaf segment. An {\em $\iota$-leaf
segment} is either an open or closed
$\iota$-leaf segment. The {\em endpoints} of
a full $\iota $-leaf segment are the points
$\kappa_{\iota,x}(u)$ and
$\kappa_{\iota,x}(v)$ where $u$ and $v$ are
the endpoints of $I_{1}$. The {\em endpoints}
of an $\iota $-leaf segment $I$ are the
points of the minimal closed full
$\iota$-leaf segment containing $I$. The {\it
interior} of an $\iota$-leaf segment $I$ is
the complement of its boundary. In
particular, an $\iota$-leaf segment $I$ has
empty interior if, and only if, it is an
$\iota$-leaf gap. A map $c:I \to \reals$ is
an {\it $\iota$-leaf chart} of an
$\iota$-leaf segment $I$ if has an extension
$c_E:I_E \to \reals$ to a full $\iota$-leaf
segment $I_E$ with the following properties:
$I \subset I_E$ and $c_E$ is a homeomorphism
onto its image.

\subsection{Smoothness}\label{sect:Smoothness}
In this
paper, when we say that a map,
atlas or structure is $C^r$ we
include the case $C^{k+}$ where $k$
is a positive integer.  For maps
$f$ this means that $f$ is
$C^{k+\alpha }$ for some $0<\alpha
<1$, i.e. $C^{k}$ with
$\alpha$-H\"older continuous
$k$th-order derivatives.  For an
atlas or structure by $C^{k+ }$ we mean that
each pair of charts in the atlas or
structure are $C^{k+\alpha }$
compatible for some $0<\alpha <1$
where the $\alpha$ might depend
upon the charts.  In the case of an
atlas, we suppose that (i) one can
choose $\alpha$ to be independent
of the charts and (ii) the overlap
maps have $C^{k+\alpha}$ norm
bounded independent of the charts
considered.  This is immediately
verified if the number of charts
contained in the $C^{k+}$ atlas is
finite.  Thus a $C^{k+}$ atlas is
$C^{k+\alpha }$, for some $0<\alpha
<1$.  This is not the case for
$C^{k+}$ structures.

\subsection{Topological and smooth conjugacies}
\label{dfgdbbccees}

Let  $(f,\Lambda)$  be a
$C^{1+}$  hyperbolic diffeomorphism.
Somewhat unusually we also desire to
highlight the $C^{1+}$ structure on
$M$ in which $f$ is a
diffeomorphism.  By a
\emph{$C^{1+}$ structure on
$M$} we mean a maximal set of
charts with open domains in $M$
such that the union of their
domains cover $M$ and
whenever $U$ is an open subset
contained in the domains of any two
of these charts $i$ and $j$ then
the overlap map $j \comp
i^{-1}:i(U) \to j(U)$ is
$C^{1+\alpha}$, where $\alpha >0$
depends on $i$, $j$ and $U$.  We
note that by compactness of
$M$, given such a $C^{1+}$
structure on $M$, there is an
atlas consisting of a finite set of
these charts which cover $M$
and for which the overlap maps are
$C^{1+\alpha}$ compatible and
uniformly bounded in the
$C^{1+\alpha}$ norm, where
$\alpha>0$ just depends upon the
atlas. We denote by
$\cC_{f}$ the $C^{1+}$ structure
on $M$ in which $f$ is a
diffeomorphism.  Usually one is not
concerned with this as, given two
such structures, there is a
homeomorphism of $M$ sending one
onto the other and thus, from this
point of view, all such structures
can be identified.  For our
discussion it will be important to
maintain the identity of the
different smooth structures on
$M$.

We say that a map $h:\Lambda_f \to \Lambda_g$
is a \emph{topological conjugacy} between two
 $C^{1+}$ hyperbolic diffeomorphisms
$(f,\Lambda_f)$ and $(g,\Lambda_g)$   if
there is a homeomorphism $h:\Lambda_f \to
\Lambda_g$ with the following properties:
\begin{rlist} \item $g \circ h(x)=h \circ
f(x)$ for every $x \in \Lambda_f$. \item The
pull-back of the $\iota $-leaf segments of
$g$ by $h$ are $\iota $-leaf segments of $f$.
\end{rlist}

We say that a topological conjugacy
$h:\Lambda_f \to \Lambda_g$ is a
\emph{Lipschitz conjugacy} if $h$ has a
bi-Lipschitz  homeomorphic extension to  an
open neighbourhood of $\Lambda_f$ in the
surface $M$ (with respect to the $C^{1+}$
structures $\cC_{f}$ and $\cC_{g}$,
respectively).

Similarly, we say that a topological
conjugacy $h:\Lambda_f \to \Lambda_g$ is a
\emph{$C^{1+}$ conjugacy} if $h$ has a
$C^{1+\alpha}$ diffeomorphic extension to  an
open neighbourhood of $\Lambda_f$ in the
surface $M$, for some $\alpha >0$.

Our approach is to fix a $C^{1+}$ hyperbolic
diffeomorphism $(f,\L)$ and consider
$C^{1+}$ hyperbolic diffeomorphism
$(g_1,\Lambda_{g_1})$ topologically conjugate
to $(f,\L)$. The topological conjugacy
$h:\Lambda \to \Lambda_{g_1}$ between $f$ and
$g_1$ extends to  a homeomorphism $H$ defined
on a neighbourhood of $\Lambda$. Then,   we
obtain the new $C^{1+}$-realization
$(g_2,\L_{g_2})$ of $f$ defined as follows:
(i) the map $g_2=H^{-1}  \circ g_1 \circ H$;
(ii) the basic set is $\Lambda_{g_2} =
H^{-1}| \Lambda_{g_1}$; (iii)  the $C^{1+}$
structure $\cC_{g_2}$ is given by the
pull-back $\left(H \right)_*\cC_{g_1}$ of the
$C^{1+}$ structure $\cC_{g_1}$. From (i) and
(ii), we get that $\Lambda_{g_2} = \Lambda$
and $g_2|\L=f$. From (iii), we get that $g_2$
is  $C^{1+}$ conjugated to $g_1$.  Hence, to
study the  conjugacy classes of $C^{1+}$
hyperbolic    diffeomorphisms $(f,\L)$ of
$f$, we can just consider the $C^{1+}$
hyperbolic diffeomorphisms $(g,\L_{g})$ with
$\L_{g} = \L$ and $g|\L = f|\L$, which we will do from now on
for simplicity of our exposition.

\subsection{Rectangles}
\label{sect:Markov_parts}
\label{spanning}

Since $\L$ is a hyperbolic invariant set of a diffeomorphism $f:M \to M$,
for $0<\eps<\eps_0$
there is $\delta =\delta (\eps )>
0$ such that, for all points $w,z
\in\L$ with $d(w,z)<\d$,  $W^u (w, \eps)$
and $W^s (z,\eps)$
intersect in a unique point
that we denote by
$[w,z]$.
Since we assume
that the hyperbolic set has a
\emph{local product structure}, we have that
$[w,z] \in \Lambda$.  Furthermore,
the following properties are satisfied:
(i)         $[w,z]$
varies continuously with
$w,z\in\L$;
(ii)        the bracket map is continuous
on a $\delta$-uniform
neighbourhood of the diagonal
in $\L\times \L$; and
(iii)        whenever both sides are defined
$f([w,z]) = [f(w),f(z)]$.
Note that the bracket map does not
really depend on $\delta$ provided
it is sufficiently small.

Let us underline that it is a
standing hypothesis that all the
hyperbolic sets considered here
have such a local product
structure.

A {\it rectangle} $R$ is a
subset of $\L$ which is (i)
closed under the bracket i.e.\
$x,y\in R\implies [x,y]\in R$,
and (ii) proper i.e.\ is the
closure of its interior in
$\L$.
This definition imposes that a
rectangle has always to be proper
which is more restrictive than the
usual one which only insists on the
closure condition.

\begin{figure}[tbp]
\centerline{\includegraphics[width=7cm]{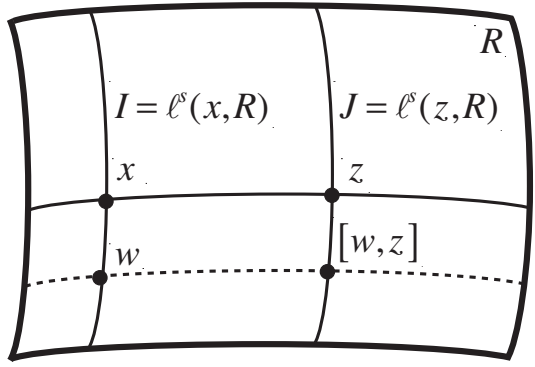} }
%\ \picill 170pt by 117pt (fig2_basicholos.eps) }
\caption{A rectangle.}
\label{Rectangle}
\end{figure}

If $\ell^s$ and $\ell^u$ are
respectively stable and unstable
leaf segments intersecting in a single point
then we denote by
$[\ell^s , \ell^u]$ the set
consisting of all points of the
form $[w,z]$ with $w\in\ell^s$ and
$z\in\ell^u$.  We note that if the
stable and unstable leaf segments
$\ell$ and $\ell'$ are closed then
the set $[\ell , \ell' ]$ is a
rectangle.  Conversely in this
$2$-dimensional situations, any
rectangle $R$ has a product
structure in the following sense:
for each $x \in R$ there are closed
stable and unstable leaf segments
of $\L$, $\ell^s(x,R)\subset
W^s(x)$ and $\ell^u(x,R)\subset
W^u(x)$ such that $R=
[\ell^s(x,R),\ell^u(x,R)]$.  The
leaf segments $\ell^s(x,R)$ and
$\ell^u(x,R)$ are called
\emph{stable and unstable spanning
leaf segments} for $R$ (see Figure \ref{Rectangle}).
For $\i \in \{s,u\}$, we denote by
$\partial\ell^{\i}(x,R)$ the set
consisting of the endpoints of
$\ell^{\i}(x,R)$, and we denote
by $\rint \ell^{\i}(x,R)$
the set
$\ell^{\i}(x,R) \setminus \partial\ell^{\i}(x,R)$.
The \emph{interior of} $R$ is given by $\rint
R= [\rint \ell^s(x,R) ,\rint
\ell^u(x,R) ]$, and the \emph{boundary of}
$R$ is given by $\partial R
=[\partial\ell^s(x,R),\ell^u(x,R) ]
\union
[\ell^s(x,R),\partial\ell^u(x,R)
]$.

\subsection{Markov partitions}
\label{markgorprp}

A {\it Markov partition of $f$}
is a collection $\cR = \{
R_1,\ldots , R_k \}$ of
rectangles such that
(i) $\Lambda \subset
\union_{i=1}^k R_i$;
(ii)         $R_i \inter R_j =\partial
R_i \inter \partial R_j$
for all $i$ and $j$;
(iii)        if $x\in \iint R_i$ and $f
x \in \iint R_j$ then
\begin{rlist}
 \item[(a)]
 $f(\ell^s(x,R_i))\subset
 \ell^s(fx,R_j)$ and
 $f^{-1}(\ell^u(fx,R_j))\subset
 \ell^u(x,R_i)$
 \item[(b)]
 $f(\ell^u(x,R_i))\inter
 R_j = \ell^u(fx,R_j)$
 and
 $f^{-1}(\ell^s(fx,R_j))\inter
 R_i = \ell^s(x,R_i)$.
\end{rlist}

The last condition means that
$f(R_i)$ goes across $R_j$ just
once.  In fact, it follows from
condition (a) providing the
rectangles $R_j$ are chosen
sufficiently small (see \cite{Mane}).
The rectangles making up the
Markov partition are called
\emph{Markov rectangles}.

We note that there is a Markov partition
$\cR$ of $f$ with the following
\emph{disjointness property} (see
\cite{Bowen,PalisNewhouse,sinais9}):
\begin{rlist}
\item
if    $0<\delta_{f,s}<1$ and
$0<\delta_{f,u}<1$ then the stable and
unstable leaf boundaries of any two Markov
rectangles do not intesect.
\item
if    $0<\delta_{f,\i}<1$ and
$\delta_{f,\ip}=1$ then the $\ip$-leaf
boundaries of any two Markov rectangles do
not intersect except, possibly, at their
endpoints.
\end{rlist}
If $\delta_{f,s}=\delta_{f,u}=1$, the
disjointness property does not apply and so
we consider that it is trivially satisfied
for every Markov partition. For simplicity of
our exposition, we will just  consider Markov
partitions   satisfying the  disjointness
property.

\subsection{Marking  the invariant set $\Lambda$}

\label{ggrdgdedx}

The properties of the  Markov partition $\cR
= \{ R_1,\ldots , R_k \}$ of $f$ imply the
existence of an unique two-sided  subshift
$\tau$ of finite type $\Theta=\Theta_A$ and a
continuous surjection $i:\Theta \to \Lambda$
such that (a) $f \circ i = i \circ \tau$ and
(b) $i(\Theta_{j})=R_j$ for every
$j=1,\ldots,k$. We call such a map $i:\Theta
\to \Lambda$     \emph{a  marking of a
$C^{1+}$ hyperbolic diffeomorphism $(f,\L)$}.

As we have explained before  a $C^{1+}$
hyperbolic diffeomorphism $(f,\L)$ admits
always a marking which is not necessarily
unique.

\subsection{Leaf $n$-cylinders and leaf $n$-gaps}

For $\i=s$ or $u$, an \emph{$\i$-leaf primary
cylinder of a Markov rectangle $R$}   is a
\label{page:primary} spanning $\i$-leaf
segment of $R$.  For $n \ge 1$, an
\emph{$\i$-leaf $n$-cylinder of $R$} is an
$\i$-leaf segment $I$ such that
\begin{rlist}
\item
$f_\i^n I$ is
an $\i$-leaf primary cylinder of a Markov rectangle $M$;
\item
$f_\i^n \left( \ell^\ip(x,R) \right) \subset M$ for  every $x \in I$.
\end{rlist}
For $n \ge 2$, an \emph{$\i$-leaf $n$-gap $G$
of $R$} is an $\i$-leaf gap $\{x,y\}$  in a
Markov rectangle $R$ such that  $n$ is the
smallest integer such that both leaves
$f_\i^{n-1} \ell^\ip (x,R)$  and  $f_\i^{n-1}
\ell^\ip(y,R)$ are contained in
$\ip$-boundaries of   Markov rectangles; An
\emph{$\i$-leaf primary  gap} $G$ is the
image $f_\i G'$ by $f_\i$ of an $\i$-leaf
$2$-gap $G'$.

We note that an $\i$-leaf segment $I$ of a
Markov rectangle $R$ can be simultaneously an
$n_1$-cylinder, $(n_1+1)$-cylinder, $\ldots$,
$n_2$-cylinder of $R$ if $f^{n_1}(I)$,
$f^{n_1+1}(I)$, $\ldots$, $f^{n_2}(I)$ are
all spanning $\i$-leaf segments. Furthermore,
if $I$ is an $\i$-leaf segment contained in
the common boundary of two Markov rectangles
$R_i$ and $R_j$ then $I$ can be an
$n_1$-cylinder of $R_i$ and an $n_2$-cylinder
of $R_j$ with $n_1$ distinct of $n_2$. If
$G=\{x,y\}$ is an $\i$-gap of $R$ contained
in the interior of $R$ then there is a unique
$n$ such that $G$ is an $n$-gap. However, if
$G=\{x,y\}$ is contained in the common
boundary of two Markov rectangles $R_i$ and
$R_j$ then $G$ can be an $n_1$-gap of $R_i$
and an $n_2$-gap of $R_j$ with $n_1$ distinct
of $n_2$. Since the number of Markov
rectangles $R_1, \ldots ,R_k$ is finite,
there is $C \ge 1$ such that, in all the
above cases for cylinders and gaps we have
$|n_2-n_1| \leq C$.

We say that a  leaf segment  $K$ is the
$i$-th \emph {mother} of an $n$-cylinder or
an $n$-gap $J$ of $R$ if $J \subset K$ and
$K$ is a leaf $(n-i)$-cylinder of $R$. We
denote $K$ by $m^iJ$.

By the properties of a Markov partition, the
smallest full $\i$-leaf ${\hat K}$ containing
a leaf $n$-cylinder $K$ of a Markov rectangle
$R$ is equal to the  union of all smallest
full $\i$-leaves containing either a leaf
$(n+j)$-cylinder or a leaf $(n+i)$-gap of
$R$, with $i\in \{1,\ldots,j\}$, contained in
$K$.

\subsection{Metric on $\Lambda$}
\label{sfddfs}

We say that a rectangle $R$ is an
$(n_s,n_u)$\emph{-rectangle} if
there is $x\in R$ such that, for
$\i=s$ and $u$, the spanning leaf
segments $\ell^\i(x,R)$ are either
an $\i$-leaf $n_\i$-cylinder or the
union of two such cylinders with a
common endpoint.

The reason for allowing the
possibility of the spanning leaf
segments being inside two touching
cylinders is to allow us to regard
geometrically very small rectangles
intersecting a common boundary of
two Markov rectangles to be small
in the sense of having $n_s$ and
$n_u$ large.

If $x,y\in\L$ and $x\neq y$ then
$d_\L (x,y) = 2^{-n}$ where $n$ is
the biggest integer such that both
$x$ and $y$ are contained in an
$(n_s,n_u)$-rectangle with $n_s\geq n$
and $n_u\geq n$.  Similarly if $I$
and $J$ are $\i$-leaf segments then
$d_\L (I,J) = 2^{-n_\ip}$ where $n_\i = 1$ and
$n_\ip$ is the biggest integer
such that both $I$ and $J$ are
contained in an $(n_s,n_u)$-rectangle.

\subsection{Basic holonomies}
\label{basicholonomies}

\begin{figure}[tbp]
\centerline{\includegraphics[width=7.5cm]{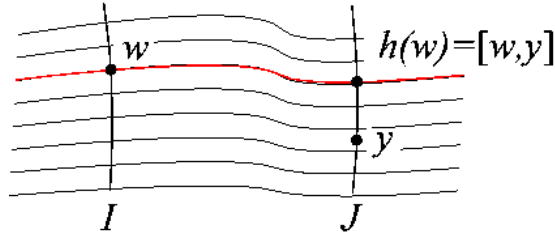} }
%\ \picill 141pt by 89pt (fig3_basich.eps) }
\caption{A basic stable holonomy from $I$ to $J$.}
\label{basic_holonomy}
\end{figure}

Suppose that $x$ and $y$ are two
points inside any rectangle $R$ of
$\L$.
Let $\ell(x,R)$ and $\ell(y,R)$ be two stable leaf
segments respectively containing
$x$ and $y$ and inside $R$.  Then
we define $\theta : \ell(x,R)\to \ell(y,R)$ by
$\theta (w) = [w,y]$.  Such maps
are called the \emph{basic stable
holonomies} (see Figure \ref{basic_holonomy}).  They generate the
pseudo-group of all stable
holonomies.  Similarly we define
the   basic unstable holonomies.

By Theorem 2.1 in \cite{Holonomy},
the holonomy $\theta: \ell^\i(x,R) \to
\ell^\i(y,R)$ has a $C^{1+\alpha}$ extension to the   leaves
containing $\ell^\i(x,R)$ and $\ell^\i(y,R)$, for some $\alpha  >
0$.

\subsection{Foliated lamination atlas}

In this section when we refer to a $C^r$ object
$r$ is allowed to take the values $k+\alpha$ where $k$ is a positive integer
and $0 < \alpha \le 1$.
Two $\i$-leaf charts $i$ and $j$  are $C^r$ compatible if
whenever $U$ is an open subset of an
$\i$-leaf segment contained in the
domains of $i$ and $j$ then $j\comp
i^{-1}:i(U)\to j(U)$ extends to a
$C^r$ diffeomorphism of the real
line.  Such maps are called
\emph{chart overlap
maps}.\label{sect:Crlamstr}\label{sect:Cr_lam_str}
A \emph{bounded $C^r$
$\i$-lamination atlas} $\cA^\i$ is
a set of such charts which (a)
cover $\L$, (b) are pairwise $C^r$
compatible, and (c) the chart
overlap maps are uniformly bounded
in the $C^r$ norm.

Let $\cA^\i$ be a bounded
$C^{1+\alpha}$ $\i$-lamination
atlas, with $0<\alpha \le 1$.  If
$i:I \to \reals$ is a chart in
$\cA^\i$ defined on the leaf
segment $I$ and $K$ is a leaf
segment in $I$ then we define
$|K|_i$ to be the length of the
minimal closed interval containing
$i(K)$.  Since the atlas is
bounded, if $j:J \to \reals$ is
another chart in $\cA^\i$ defined
on the leaf segment $J$ which
contains $K$ then the
ratio between the lengths $|K|_i$
and $|K|_j$ is universally bounded
away from $0$ and $\infty$.  If $K'
\subset I \inter J$ is another such
segment then we can define the
ratio $r_i(K:K')=|K|_i/|K'|_i$.
Although this ratio depends upon
$i$, the ratio is exponentially
determined in the sense that if $T$
is the smallest segment containing
both $K$ and $K'$ then
$$
r_j\left (K:K' \right ) \in \left( 1 \pm
\cO \left(|T|_i^\alpha
\right) \right) r_i \left (K:K' \right ) \ .
$$
This follows from the
$C^{1+\alpha}$ smoothness of the
overlap maps and Taylor's Theorem.

A $C^{r}$ lamination atlas $\cA^\i$ has
\emph{bounded geometry}
(i) if $f$ is a $C^r$ diffeomorphism with  $C^r$ norm
uniformly bounded in this
atlas;
(ii)         if for all pairs $I_1, I_2$
of $\i$-leaf $n$-cylinders
or $\i$-leaf $n$-gaps with
a common point, we have that
$r_i(I_1:I_2)$ is
uniformly bounded away from
$0$ and $\infty$ with the
bounds being independent of
$i$, $I_1$, $I_2$ and $n$;
and (iii) for all endpoints
$x$ and $y$ of an $\i$-leaf
$n$-cylinder or $\i$-leaf
$n$-gap $I$, we have that
$|I|_i \le
\cO\left((d_\L(x,y))^\beta
\right)$ and $d_\L(x,y) \le
\cO\left(|I|_i^\beta
\right)$, for some $0<\beta
<1$, independent of $i$, $I$ and $n$.

A $C^{r}$ bounded lamination atlas $\cA^\i$
is $C^{r}$\emph{foliated}
(i) if   $\cA^\i$ has bounded geometry; and
(ii) if         the basic holonomies are
$C^r$ and have a $C^r$ norm
uniformly bounded in this
atlas, except possibly for
the dependence upon the
rectangles defining the
basic holonomy.
A  bounded lamination atlas $\cA^\i$ is $C^{1+}$\emph{foliated}
if $\cA^\i$ is  $C^{r}$foliated for  some $r > 1$.

\subsection{Foliated atlas $\cA^\i(g,\rho)$}
\label{sfgvgbg}

Let $g \in  \cT(\fL )$
and $\rho=\rho_g$ be a   $C^{1+}$ Riemannian metric on
the manifold containing $\L$.
The  \emph{$\i$-lamination atlas $\cA^\i (g,\rho)$
determined by} $\rho$ is the
set of all maps $e:I\to \reals$
where $I=\L\cap \hat{I}$ with
$\hat{I}$  a full $\i$-leaf
segment, such that
$e$ extends to an isometry between
the induced Riemannian metric on
$\hat{I}$ and the Euclidean
metric on the reals.  We call the maps
$e\in\cA^\i (\rho) $ the \emph{$\i$-lamination
charts}.  If $I$ is an $\i$-leaf
segment (or a full $\i$-leaf
segment) then by $|I|_{\rho}$ we
mean the length in the Riemannian
metric $\rho$ of the minimal full
$\i$-leaf containing $I$.
By Theorem 2.2 in \cite{Holonomy}, the lamination atlas
$\cA^\i (g,\rho)$   is $C^{1+}$foliated for $\i=\{s,u\}$.

\section{Solenoid functions}
\label{fgvdghbvc}

In this section, we construct the
stable and unstable solenoid functions,
and we  prove an equivalence between $C^{1+}$ hyperbolic diffeomorphisms
and pairs of stable and unstable solenoid functions.

\subsection{HR-H\"older ratios}
A \emph{HR-structure} associates an
affine structure to each stable and
unstable leaf segment in such a way
that these vary H\"older
continuously with the leaf and are
invariant under $f$.

An affine structure on a stable or
unstable leaf is equivalent to a
\emph{ratio function} $r(I:J)$ which
can be thought of as prescribing
the ratio of the size of two leaf
segments $I$ and $J$ in the same
stable or unstable leaf.  A {\it
ratio function} $r(I:J)$ is
positive (we recall that each leaf segment has at least two distinct points)
and continuous in the
endpoints of $I$ and $J$.
Moreover,
\begin{equation}
%\label{eq1}
\label{eq:ratios}
r(I:J) = r(J:I)^{-1} ~~~{\rm
and}~~~ r(I_1\cup I_2 :K) = r(I_1:K)
+ r(I_2:K)
\end{equation}
provided $I_1$ and $I_2$ intersect
at most in one of their endpoints.

We say that $r$ is a {\it
$\i$-ratio
function}
if (i) for all $\i$-leaf
segments $K$, $r(I:J)$
($I,J\subset K$) defines a
ratio function on $K$; (ii) $r$
is invariant under $f$, i.e.
$r(I:J)=r(fI: fJ)$ for all
$\i$-leaf segments; and
(iii) for every basic
$\i$-holonomy map $\theta:I\to J$
between the leaf segment
$I$ and the leaf segment
$J$  defined  with respect to a rectangle $R$
and for every
$\i$-leaf segment $I_0\subset I$ and
every $\i$-leaf segment or gap
$I_1\subset I$,
\begin{equation}
     \label{1eq:Holderness_of_ratios}
     \left| \log \frac{r(\theta
     I_0 :\theta
     I_1)}{r(I_0:I_1)} \right|
     \le \cO \left(\left( d_\Lambda
     (I,J
     )\right)^\epsilon \right)
\end{equation}
where $\epsilon\in (0,1)$ depends upon $r$ and
the constant of proportionality
also depends upon   $R$, but not on
the segments considered.

A {\it HR-structure} is a pair
$(r_s,r_u)$ consisting of a stable
and an unstable ratio function.

\subsection{Realised ratio functions}

Let $(g,\L) \in \cT(\fL )$ and let $\cA(g,\rho)$ be
an  $\i$-lamination atlas which is $C^{1+}$ foliated.
Let $|I|=|I|_\rho$ for every $\i$-leaf segment $I$.
By hyperbolicity of $g$ in
$\Lambda$, there are $0<\nu<1$ and
$C>0$ such that for all $\i$-leaf
segments $I$ and all $m \ge 0$ we
get
$|g_\ip^{m}I| \le C \nu^m |I|$.
Thus, using the
mean value theorem and the
fact that $g_\i$
is $C^r$, for all short leaf
segments $K$ and all leaf segments
$I$ and $J$ contained in it, the
$\i$-realised ratio function $r_{g,\i}$
given by
\begin{eqnarray*}\label{eqn:12131}
r_{g,\i} (I:J) & = &
\lim_{n \to \infty}
\frac {|g_\ip^{n}I|}{|g_\ip^{n} J|} \nonumber \\
& = &
\frac{|g_\ip^{m}I|}{|g_\ip^{m} J|}
\prod_{n=m}^{\infty} \left(
\frac
{|g_\ip^{n+1}I|}{|g_\ip^{n+1} J|}
\frac {|g_\ip^{n} J|}{|g_\ip^{n} I|} \right) \nonumber \\
& \in &
\frac{|g_\ip^{m}I|}{|g_\ip^{m} J|}
\prod_{n=m}^{\infty}
\left(1 \pm \cO \left(\nu^n |K|^\alpha \right) \right) \nonumber \\
& \subset &
\frac{|g_\ip^{m}I|}{|g_\ip^{m} J|}
\left(1 \pm \cO \left( \nu^n |K|^\alpha \right) \right)  \nonumber
\end{eqnarray*}
is well-defined, where $\alpha = \min\{1, r-1\}$.
This construction gives the
HR-structure on $\L$ determined by
$g$.�By \cite{HR}, we get the following equivalence:

\begin{theorem}\label{thm:HR}
The map $g\to (r_{g,s},r_{g,u})$ determines a one-to-one correspondence
between $C^{1+}$ conjugacy classes in $\cT(\fL )$
and HR-structures.
\end{theorem}

\subsection{Lamination atlas}

\begin{figure}[tbp]
 \centerline{\includegraphics[width=7cm]{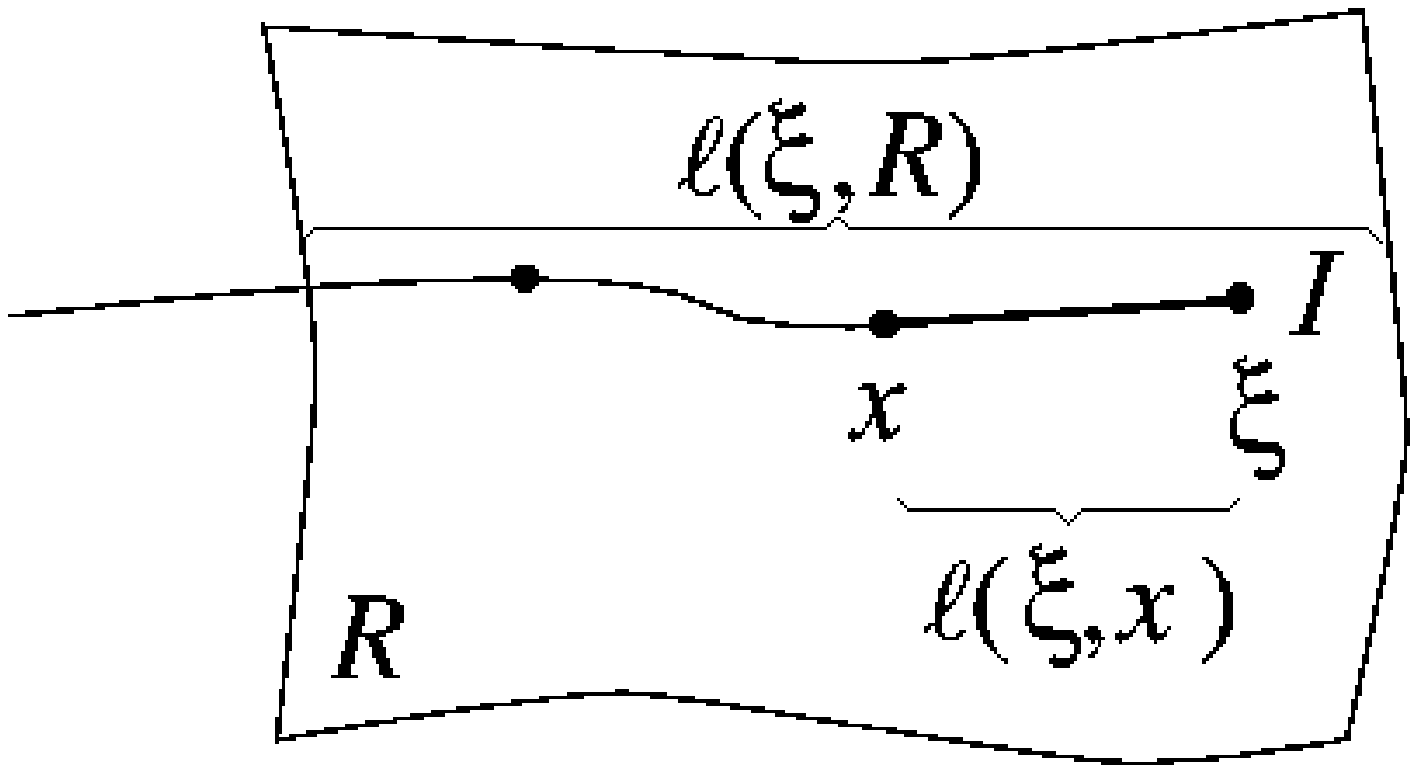} }
 %\ \picill 170pt by 117pt (fig2_basicholos.eps) }
 \caption{The embedding $e:I \to \reals$.}
 \label{embedding}
 \end{figure}

Given an $\i$-ratio function $r$, we
define the embeddings $e:I \to
\reals$ by
\begin{equation}\label{eqn:e_charts}
e(x)=r(\ell (\xi , x),        \ell(\xi,R))
\end{equation}
where $\xi$ is an endpoint of the
$\i$-leaf segment $I$ and $R$ is a
Markov rectangle containing $\xi$ (see Figure \ref{embedding}).
For this definition it is not
necessary that $R$ contains $I$.  We
denote the set of all these
embeddings $e$ by $\cA(r)$.

The embeddings $e$ of $\cA (r)$
have overlap maps with affine
extensions, therefore the atlas
$\cA (r)$ extends to a $C^{1+\a}$
lamination structure $\cL (r)$.
By Proposition 4.2
in \cite{Holonomy}, we obtain
that $\cA (r)$ is a $C^{1+}$foliated atlas.

Let $g \in \cT(\fL )$ and  $\cA(g,\rho)$
a   $C^{1+}$ foliated $\i$-lamination atlas
determined by a Riemmanian metric $\rho$.
Putting together Proposition 2.5 and Proposition 3.5 of \cite{HR}, we get
that the overlap map
$e_1 \circ e_2^{-1}$ between a chart $e_1 \in \cA(g,\rho)$
and a chart $e_2 \in  \cA (r_{g,\i})$ has a $C^{1+}$
diffeomorphic extension to the reals.
Therefore, the atlasses $\cA(g,\rho)$ and
$\cA (r_{g,\i})$ determine the same $C^{1+}$ foliated
$\i$-lamination. In particular,
for all short leaf
segments $K$ and all leaf segments
$I$ and $J$ contained in it, we obtain that
$$
r_{g,\i} (I:J) =
\lim_{n \to \infty}
\frac {|g_\ip^{n}I|_{\rho}}{|g_\ip^{n} J|_{\rho}}
=\lim_{n \to \infty}
\frac {|g_\ip^{n}I|_{i_n}}{|g_\ip^{n} J|_{i_n}}
$$
where $i_n$ is any chart in $ \cA (r_{g,\i})$
containing the segment $g_\ip^{n}K$ in its domain.

\subsection{Realised solenoid functions}

For $\i=s$ and $u$, let $\sol^\i$
denote the set of all ordered pairs
$(I,J)$ of $\i$-leaf segments
with the following properties:
\begin{rlist}
\item
The intersection of $I$ and $J$
consists of a single endpoint.
\item
if $\delta_{f,\i}=1$ then $I$ and $J$ are primary $\i$-leaf cylinders.
\item
if $0 < \delta_{f,\i}<1$ then
$f_\ip I$  is an $\i$-leaf  $2$-cylinder of a Markov rectangle $R$
and $f_\ip J$
is an $\i$-leaf $2$-gap also of the same Markov rectangle $R$.
\end{rlist}
(See \S
~\ref{sect:Markov_parts} for the
definitions of leaf cylinders and
gaps).  Pairs $(I,J)$ where both
are primary cylinders are called
\emph{leaf-leaf pairs}.
Pairs  $(I,J)$  where $J$ is a gap
are called  \emph{leaf-gap pairs}
and in this case we refer to   $J$ as a \emph{primary gap}.
The set $\sol^\i$ has a very nice topological structure.
If $\delta_{f,\ip}=1$ then the set $\sol^\i$ is
isomorphic to a finite union
of intervals, and if
$\delta_{f,\ip} < 1$ then the set $\sol^\i$ is isomorphic
to an embedded Cantor set.

We define a pseudo-metric
$d_{\sol^\i}:\sol^\i \times \sol^\i
\to \reals^+$ on the set $\sol^\i$
by
$$
d_{\sol^\i}\left (\left
(I,J\right),\left
(I',J'\right)\right)= \max \left \{
d_\L \left (I,I'\right), d_\L \left
(J,J'\right)\right\} \ .
$$
Let $g \in \cT(\fL )$.
For $\i=s$ and $u$, we call the restriction
of an $\i$-ratio function $r_{g,\i}$ to $\sol^i$
a \emph{realised solenoid function} $\s_{g,\i}$.
By construction, for $\i=s$ and $u$, the restriction
of an
$\i$-ratio function to
$\sol^i$ gives an H\"older continuous function
satisfying the  matching condition, the boundary
condition and the cylinder-gap condition
as we now proceeed to describe.

\subsection{H\"older continuity
of solenoid functions} This means
that for $t=(I,J)$ and $t'=(I',J')$
in $\sol^\i$, $\left
|\s_\i(t)-\s_\i(t')\right | \leq \cO
\left(\left(d_{\sol^\i} \left
(t,t'\right )\right)^\alpha \right)
.$ The H\"older continuity of  $\s_{g,\i}$ and the compactness of its domain
imply that  $\s_{g,\i}$ is  bounded away from zero and
infinity.

\subsection{Matching
condition}

\begin{figure}[tbp]
\includegraphics[width=4in]{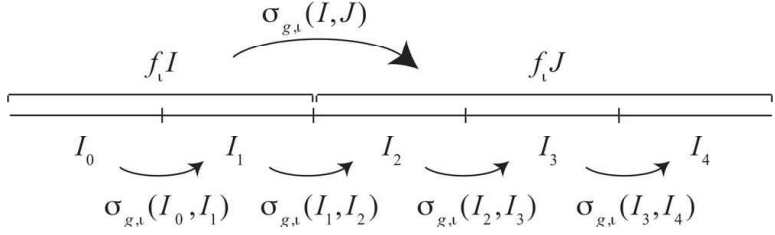}
\caption{The $f$-matching condition for $\i$-leaf segments.}
\label{ref:matching}
\end{figure}

Let $(I,J) \in \sol^\i$ be  a pair of primary cylinders
and  suppose that we
have pairs $$(I_0,I_1),
(I_1,I_2),\ldots ,
(I_{n-2},I_{n-1})\in \sol^\i$$
of primary cylinders such that
$f_\i I=\union_{j=0}^{k-1}  I_j$ and
$f_\i J=\union_{j=k}^{n-1}  I_j$.  Then
$$\frac{|f_\i I|}{|f_\i J|}
=
\frac{\sum_{j=0}^{k-1} |I_j|}
{\sum_{j=k}^{n-1} |I_j|}
=
\frac{1+
\sum_{j=1}^{k-1}
\prod_{i=1}^{j}
| I_i|/| I_{i-1}|}
{\sum_{j=k}^{n-1}
\prod_{i=1}^{j}|  I_i|/|  I_{i-1}|
}
\ .
$$
Hence, noting that $g|\Lambda=f|\Lambda$, the realised solenoid function
$\s_{g,\i}$  must
satisfy the   \emph{matching
condition} (see Figure
\ref{ref:matching})
for all such leaf segments:
\begin{equation}
%\label{eqn:mathching}
\label{eqn:solenoid_matching}
\s_{g,\i} (I:J) =
\frac{1+\sum_{j=1}^{k-1}\prod_{i=1}^{j}\s_{g,\i}
(I_i:I_{i-1})}{\sum_{j=k}^{n-1}
\prod_{i=1}^{j}\s_{g,\i}
(I_i:I_{i-1})} \ .
\end{equation}

\subsection{Boundary condition}
\label{uyggrersx}

If the stable and unstable leaf segments
have Hausdorff dimension   equal to $1$,
then leaf segments $I$ in the boundaries
of Markov rectangles can sometimes
be written as the union of primary
cylinders in more than one way.
This gives rise to the existence of a boundary condition
that
the realised solenoid functions have to satisfy
as we pass to explain.

\begin{figure}[tbp]
\includegraphics[width=9cm]{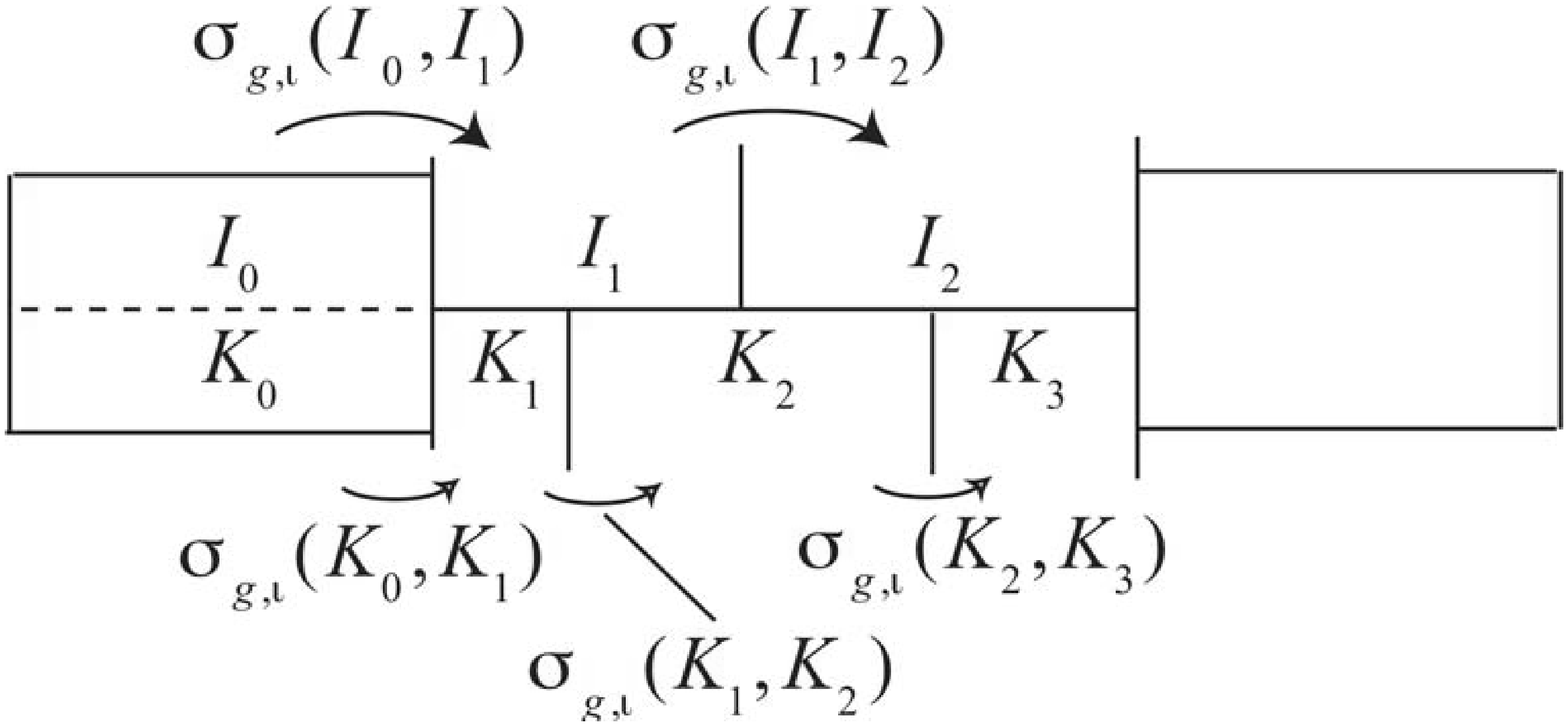}
\caption{The boundary condition for  $\i$-leaf segments.}
\label{ref:boundary}
\end{figure}

If
$J$ is another leaf segment
adjacent to the leaf segment  $I$ then the value of
$|I|/|J|$ must be the same whichever
decomposition we use.   If we write
$J=I_0=K_0$ and $I$ as
$\union_{i=1}^mI_i$ and
$\union_{j=1}^n K_j$ where
the $I_i$ and $K_j$ are
primary cylinders with $I_i\neq
K_j$ for all $i$ and $j$,
then the above two ratios are
$$
\sum_{i=1}^{m} \prod_{j=1}^{i}
\frac{|  I_j|}{|  I_{j-1}|}  =
\frac{|  I|}{|  J|} =
\sum_{i=1}^{n} \prod_{j=1}^{i}
\frac{|  K_j^|}{|  K_{j-1}|} \ .
$$
Thus,  noting that $g|\Lambda=f|\Lambda$,  a realised solenoid function
$\s_{g,\i}$ must satisfy the following
\emph{boundary condition} (see
Figure \ref{ref:boundary}) for all such leaf segments:
\begin{equation}
\label{eqn:solenoid_boundary}
\sum_{i=1}^{m}\prod_{j=1}^{i} \s_{g,\i}
\left (I_j:I_{j-1} \right ) =
\sum_{i=1}^{n}\prod_{j=1}^{i} \s_{g,\i}
\left (K_j:K_{j-1} \right ) \ .
\end{equation}

\subsection{Scaling function}

If the $\i$-leaf segments
have Hausdorff dimension less than one and the $\ip$-leaf
segments have Hausdorff dimension equal to $1$,
then a primary cylinder $I$ in the $\i$-boundary
of a Markov rectangle can also
be written as the union of gaps and
cylinders of other Markov rectangles.
This gives rise to the existence of a cylinder-gap condition
that
the $\i$-realised solenoid functions have to satisfy.

Before defining the cylinder-gap condition, we will introduce
the scaling function that will be useful to express the cylinder-gap
condition, and also, in Definitions \ref{fgrdegerdg}, the bounded
equivalence classes of solenoid functions
and, in Definition \ref{dfgdrfgggx}, the $(\delta,P)$-bounded solenoid equivalence classes of
a Gibbs measure.

Let $\Sol^\i$ be the set of all  pairs  $(K,J)$
of $\i$-leaf segments with the following properties:
\begin{rlist}
\item $K$ is a  leaf $n_1$-cylinder or an $n_1$-gap  segment for some $n_1
> 1$;
\item $J$ is a  leaf $n_2$-cylinder or an $n_2$-gap  segment for some
$n_2
> 1$;
\item $m^{n_1-1}K$ and $m^{n_2-1}J$ are the same  primary cylinder.
\end{rlist}

\begin{lemma}
\label{llkk} Every fuction $\sigma_\i : \sol^\i \to \reals^+$ has a
canonical extension $s_\i$ to $\Sol^\i$. Furthermore, if $\sigma_\i$
is the restriction of a ratio function
$r_{\i}|\sol^\i$ to $\sol^\i$ then  $s_\i= r_{\i}|\Sol^\i$.
\end{lemma}

\begin{remark}
The above map $s_\i : \Sol^\i \to \reals^+$ is the {\it scaling
function} determined by the solenoid function
$\sigma_\i : \sol^\i \to \reals^+$.
\end{remark}

\noindent
\emph{Proof of Lemma \ref{llkk}}.
We are going to give an explicit construction of a
realised scaling function $s_{g,\i}$
from a realised solenoid function  $\s_{g,\i}$
with the   property that  $s_{g,\i}= r_{g,\i}|\Sol^i$
where $r_{g,\i}$ is a ratio function,
i.e  for every
$(K,J) \in  \Sol^\i$ we have
$$
s_{g,\i}(K,J)= \lim_{k \to \infty} \frac{|g_\ip^k K|_\rho}{|g_\ip^k J|_\rho}
$$
where $\cA(g,\rho)$ is an $\i$-lamination atlas.

This construction is canonical and applies to every
function $\sigma_\i : \sol^\i \to \reals^+$ determining
a canonical extension $s_\i : \Sol^\i \to \reals^+$
of $\sigma_\i$.

Let us proceed to construct the   $\i$-scaling function
$s: S^\i \to \reals^+$
from an $\i$-solenoid function $\sigma$.
Suppose that $J$ is an $\i$-leaf $n$-cylinder
or   $n$-gap. Then there are   pairs $$(I_0,I_1),
(I_1,I_2),\ldots ,
(I_{l-1},I_l)\in \sol^\i$$
such that
$mJ=\union_{j=0}^l f_\ip^{n-1} I_j$ and
$J= f_\ip^{n-1} I_s$ for some $0 \le s \le l$.
Let us denote    $f_\ip^{n-1} I_j$ by $I_j'$  for every $0 \le s \le l$.
Then
$$
\frac{|mJ|}{|J|}
= \sum_{j=0}^l
\frac{|I_j'|} {|I_s'|}
= 1+
\sum_{j=0}^{s-1}    \prod_{i=s}^{j+1}
\frac{| I_{i-1}'|} {|  I_i'|}
+ \sum_{j=s+1}^{l}    \prod_{i=s}^{j-1}
\frac{|  I_{i+1}'|} {|   I_i'|} \ ,
$$
where  the first sum above is empty   if $s=0$,
and  the second sum above is empty if $s=1$.
Therefore, we define the extension $s$ from  $\sigma$ to
the pairs $(mJ,J)$   by
$$
s(mJ,J)= 1+
\sum_{j=0}^{s-1}
\prod_{i=s}^{j+1}
\sigma(I_{i-1},I_i)
+ \sum_{j=s+1}^{l}
\prod_{i=s}^{j-1}  \sigma(I_{i+1},I_i) \ ,
$$
where  the first sum above is empty   if $s=0$,
and  the second sum above is empty  if $s=1$.
For every  $(K,J) \in \Sol^\i$ there is a primary leaf segment $I$
such that  $m^{m_1}K= I = m^{m_2}J$ for some $m_1 \ge 1 $ and $m_2\ge  1$.
Then,
$$
\frac{| K|}{|  J|} =
\prod_{j=1}^{m_1} \frac{| m^{j} J|}{|  m^{j-1} J|}
\prod_{j=1}^{m_2} \frac{| m^{j-1} K|}{|  m^j K|}   \ .
$$
Therefore, we define the extension $s$  to
$(K,J)$  by
$$
s(K,   J) =    \prod_{j=1}^{m_1}  s(m^jJ, m^{j-1} J)
\prod_{j=1}^{m_2}  s(m^{j-1} K, m^jK) \ .
$$
Hence, we have constructed a scaling function $s$ from $\sigma$ on the set
$\Sol^\i$ such that if $\sigma$ is the restriction of a ratio function
$r_{g,\i}|\sol^\i$ to $\sol^\i$ then  $s= r_{g,\i}|\Sol^\i$.\qed

\subsection{Cylinder-gap condition}
\label{dfsfbbvv}

\begin{figure}[tbp]
\includegraphics[width=11cm]{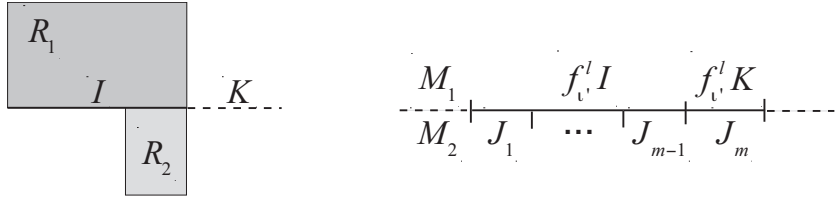}
\caption{The cylinder-gap condition for  $\i$-leaf segments.}
\label{ref:gap}
\end{figure}

Let $(I,K)$ be  a leaf-gap pair such that the primary cylinder
$I$ is the $\i$-boundary of a Markov rectangle $R_1$.
Then the  primary cylinder
$I$ intersects another
Markov rectangle $R_2$ giving rise to the existence of a cylinder-gap condition
that
the realised solenoid functions have to satisfy
as we proceed to explain.
Take the smallest $l \ge 0$ such that $f^l_\ip I \cup f^l_\ip K$
is contained in the intersection of the boundaries of two Markov rectangles
$M_1$ and $M_2$. Let $M_1$ be the Markov rectangle with the property
that $M_1 \cap f^l_\ip R_1$ is a rectangle with non-empty interior
(and so $M_2 \cap f^l_\ip R_2$ also has non-empty interior).
Then, for some positive $n$, there are distinct $n$-cylinder and gap leaf segments
$J_1,\ldots,J_m$ contained in a primary cylinder of $M_2$ such that
$f^l_\ip K=J_m$ and the smallest full $\i$-leaf segment
containing $f^l_\ip I$ is equal to the union
$\cup_{i=1}^{m-1} {\hat J}_i$, where ${\hat J}_i$ is the smallest
full $\i$-leaf segment containing $J_i$.
Hence,
$$
\frac{|f^{l}_\ip I|}{|f^{l}_\ip K|} =
\sum_{i=1}^{m-1}
\frac{|J_i|}{|J_m|} \ .
$$
Hence,  noting that $g|\Lambda=f|\Lambda$,  a realised  solenoid function $\s_{g,\i}$
must satisfy the
\emph{cylinder-gap condition} (see
Figure \ref{ref:gap}) for all such leaf segments:
$$
\sigma_{g,\i}(I, K) =   \sum_{i=1}^{m-1} s_{g,\i}(J_i,J_m)
$$
where  $s_{g,\i}$ is the    scaling
function   determined by the solenoid function $\s_{g,\i}$.

\subsection{Solenoid functions}

\label{SOLLLLLL}

Now, we are ready to present the definition of an $\i$-solenoid
function.

\begin{definition}
An H\"older continuous function $\s_\i:\sol^\i \to \reals^+$
is an \emph{$\i$-solenoid function}
if $\s_\i$ satisfies the matching condition, the boundary
condition and the cylinder-gap condition.
\end{definition}

We denote by ${\mathcal PS}(f)$ the set of pairs $(\sigma_s,\sigma_u)$
of stable and unstable solenoid functions.

\begin{remark}
Let $\s_\i:\sol^\i \to \reals^+$ be an $\i$-solenoid function.
The matching, the boundary
and the cylinder-gap  conditions are trivially satisfied
except in the following cases:
\begin{rlist}
\item
The matching condition
if $\delta_{f,\i}=1$.
\item
The boundary condition
if $\delta_{f,s}=\delta_{f,u}=1$.
\item
The cylinder-gap condition
if $\delta_{f,\i}<1$ and $\delta_{f,\ip} = 1$.
\end{rlist}
\end{remark}

\begin{lemma}
\label{gdfrrn}
The map $r_\i \to  r_\i|\sol^\i$ gives a one-to-one correspondence
between $\i$-ratio functions and $\i$-solenoid functions.
\end{lemma}

\begin{proof}
Every  $\i$-ratio function restricted to
the set $\sol^\i$ determines an $\i$-solenoid function
$r_\i|\sol^\i$.
Now we prove the converse.  Since
the solenoid functions are
continuous and their domains are
compact they are bounded away from
$0$ and $\infty$.  By this
boundedness and the $f$-matching
condition of the solenoid functions
and by iterating the domains
$\sol^s$ and $\sol^u$ of the
solenoid functions backward and
forward by $f$, we determine the
ratio functions $r^s$ and $r^u$ at
very small (and large) scales, such
that $f$ leaves the ratios
invariant.  Then, using the
boundedness again, we extend the
ratio functions to all pairs of
small adjacent leaf segments by
continuity.
By the boundary condition and the cylinder-gap condition of the
solenoid functions, the ratio
functions are well determined at
the boundaries of the Markov
rectangles.   Using the H\"older
continuity of the
solenoid function, we   deduce inequality
\eqref{1eq:Holderness_of_ratios}.
\end{proof}

The set ${\mathcal PS}(f)$ of all pairs $(\s_s,\s_u)$ has a natural
metric.
Combining Theorem \ref{thm:HR} with Theorem \ref{gdfrrn},
we obtain that the set ${\mathcal PS}(f)$  forms a moduli space
for the $C^{1+}$ conjugacy classes of $C^{1+}$ hyperbolic diffeomorphisms
$g \in \cT(\fL )$:

\begin{corollary}
\label{cor:sol}
The map $g \to (r_{g,s}|\sol^s,r_{g,u}|\sol^u)$
determines a one-to-one correspondence
between $C^{1+}$ conjugacy classes of $g \in \cT(\fL )$
and pairs of solenoid functions in ${\mathcal PS}(f)$.
\end{corollary}

\begin{definition}
\label{fgrdegerdg}
We say that any two $\i$-solenoid functions $\s_1:\sol^\i \to \reals^+$
and $\s_2:\sol^\i \to \reals^+$
are in the same \emph{bounded equivalence class} if
the corresponding   scaling functions
$s_1:\Sol^\i \to \reals^+$ and $s_2:\Sol^\i \to \reals^+$
satisfy the following property:
There is $C>0$ such that for every  $\i$-leaf $(i+1)$-cylinder or $(i+1)$-gap
$J$
\begin{equation}
\label{dsdfgreww}
\left| \log s_1(J,m^iJ) - \log s_2(J,m^iJ) \right|
< C   \  .
\end{equation}
\end{definition}

Later, in Lemma \ref{asdro},
we prove that two $C^{1+}$ hyperbolic diffeomorphisms
$g_1$ and $g_2$ are Lipschitz conjugate if, and only if,
the solenoid functions $s_{g_1,\i}$ and $s_{g_2,\i}$
are in the same bounded equivalence class
for $\i$ equal to $s$ and $u$.

\section{Self-renormalisable structures}
\label{bgtrewss}

In this section, we construct the
stable and unstable self-renormalisable structures living in $1$-dimensional spaces,
and we  prove an equivalence between $C^{1+}$ hyperbolic diffeomorphisms
and pairs of stable and unstable self-renormalisable structures.

\subsection{Train-tracks}\label{sect:train_tracks}

Roughly
speaking train-tracks are the optimal
leaf-quotient spaces on which the
stable and unstable Markov maps
induced by the action of $f$ on
leaf segments are local
homeomorphisms.

For each Markov rectangle $R$ let
$t_R^\i$ be the set of
$\ip$-segments of $R$.  Thus by the
local product structure one can
identify $t_R^\i$ with any spanning
$\i$-leaf segment $\ell^\i (x,R)$
of $R$.

We form the space $\TT^\i$ by
taking the disjoint union
$\bigsqcup_{R \in \cR} t_R^\i$
(union over all Markov rectangles $R$ of the Markov partition $\cR$) and
identifying two points $I\in
t_R^\i$ and $J\in t_{R'}^\i$ if
either (i) the $\ip$-leaf segments
$I$ and $J$ are $\ip$-boundaries of Markov rectangles
and their intersection contains at least a point  which is not an
endpoint of $I$ or $J$ or (ii)
there is a sequence
$I=I_1,\ldots,I_n=J$ such that all
$I_i,I_{i+1}$ are both identified
in the sense of (i).  This space is
called the $\i$-\emph{train-track}
and is denoted $\TT^\i$.

Let $\pi_{\TT^\i}:\bigsqcup_{R \in \cR} R \to \TT^\i$ be the natural
projection sending $x\in R$ to
the point in $\TT^\i$ represented
by $\ell^\ip (x,R)$.  A
\emph{topologically regular point} $I$
in $\TT^\i$ is a point with an
unique preimage under $\pi_{\TT^\i}$ (i.e. the pre-image of $I$ is not a union of distinct
$\ip$-boundaries of Markov rectangles).  If
a point has more than one preimage
by $\pi_{\TT^\i}$ then we call it a
\emph{junction}.
Since there are only a finite number of
$\ip$-boundaries of Markov rectangles
there are only finitely many junctions (see Figures \ref{Train_tracksAnosov}).

% and \ref{Train_tracksAnosov3}).
%\begin{figure}
%\includegraphics[width=10cm]{Anosov.eps}
%\caption{This figure illustrates the stable and unstable (horocycle)
%train tracks for the Anosov map}
%\label{Train_tracksAnosov3}
%\end{figure}

Let $d_\i$\label{page:metric} be the
metric on $\TT^\i$ defined as
follows: if $\xi,\eta \in \TT^\i$,
$d_{\TT^\i}(\xi,\eta ) = d_\L(\xi,\eta)$.

\subsection{Train-track segments and charts.}
\label{fddhhhuy}

We say that $I_T$ is a \emph{train-track segment}
if there is an $\i$-leaf segment $I$,
not intersecting $\i$-boundaries of Markov rectangles,
such that $\pi_\i|I$ is an injection and
$\pi_\i(I) = I_T$.
Let $\cA$ be an $\i$-lamination atlas
(take for instance $\cA$ equal to $\cA^\i(f,\rho)$ or to $\cA(r_{f,\i})$).
The chart $i:I \to \reals$
in $\cA$ determines a
\emph{train-track chart} $i_T:I_T \to \reals$ for $I_T$
given by $i_T = i \comp \pi_\i^{-1}$.
We denote by $\cB$ the set of all
train-track charts  for all train-track segments determined by $\cA$.

Two train-track  charts
$(i_T,I_T)$ and $(j_T,J_T)$ on the train-track $\TT^\i$ are
\emph{$C^{1+}$-compatible} if the
\emph{overlap map}
$j_T \comp i_T^{-1}: i_T(I_T \cap J_T) \to j_T(I_T \cap J_T)$
has a
$C^{1+}$ extension.
A $C^{1+}$ atlas $\cB$ is a set of $C^{1+}$-compatible charts
with the following property: For  every short train-track segment
$K_T$  there is a chart $(i_T,I_T) \in \cB$
such that $K_T \subset I_T$.
A \emph{
$C^{1+}$ structure $\cS$ on
$\TT^\i$} is a maximal set of
$C^{1+}$-compatible charts with a given atlas $\cB$ on
$\TT^\i$. We say that two $C^{1+}$ structures
$\cS$ and $\cS'$ are in the same
\emph{Lipschitz equivalence class} if
for every chart   in $\cS$ and
every chart  in $\cS'$
the overlap map $ e_1 \circ e_2^{-1}$
has a  bi-Lipschitz extension.

Given any train-track charts $i_T:I_T \to \reals$
and $j_T:J_T \to \reals$ in $\cB$,
the overlap map
$j_T \comp i_T^{-1}: i_T(I_T \cap J_T) \to j_T(I_T \cap J_T)$
is equal to $j_T \comp i_T^{-1} = j \comp h \comp i^{-1}$
where $i=i_T \comp \pi_\i : I \to \reals$
and $j=j_T \comp \pi_\i : J \to \reals$
are charts in $\cA$, and
$$h : i^{-1}(i_T(I_T \cap J_T)) \to j^{-1}(j_T(I_T \cap J_T))$$
is a basic $\i$-holonomy.
Let us denote by $\cB^\i(g,\rho)$ and $\cB(r_{g,\i})$
the train-track atlasses determined respectively by
$\cA^\i(g,\rho)$ and $\cA(r_{g,\i})$ with $g \in \cT(\fL )$.
Since $\cA^\i(g,\rho)$ and $\cA(r_{g,\i})$ are $C^{1+}$foliated atlases,
there is $\eta > 0$
such that, for all  train-track charts $i_T$ and $j_T$
in  $\cB^\i(g,\rho)$ (or in $\cB(r_{g,\i})$),
the overlap maps
$j_T \comp i_T^{-1} = j \comp h \comp i^{-1}$
have $C^{1+\eta}$ diffeomorphic
extensions with a uniform bound
for the $C^{1+\eta}$ norm.
Hence, $\cB^\i(g,\rho)$ and $\cB(r_{g,\i})$ are $C^{1+\eta}$ atlas.

\subsection{Markov maps}
\label{gfgrtr}

The \emph{Markov map}
$\tau_\i :\TT^\i\to\TT^\i$
is the
mapping induced by the action of
$f$ on leaf segments i.e.\ it is
defined as
follows:\label{page:markov} if
$I \in \TT^\i$,
$\tau_{\i} I=\pi_{\TT^\i}(f_\i I)$ is the $\ip$-leaf
segment containing the
$f_\i$-image of the $\ip$-leaf segment
$I$ (for simplicity of notation we use
the same symbols for the Markov maps
as for the shift maps). This map
$\tau_\i$
is a local
homeomorphism because $f_\i$ sends a
short $\i$-leaf segment
homeomorphically onto a short
$\i$-leaf segment.
For simplicity of notation, we will denote  $\tau_\i$ by $f_\i$ through the paper.

Given a topological  chart $(e,U)$ on the train-track $\TT^\i$
and a train-track segment $C \subset U$,
we denote by $|C|_e$ the lenght of the smallest interval containing $e(C)$.
We say that $f_\i$ has
\emph{bounded geometry} in   a $C^{1+}$ atlas $\cB$
if there is $\kappa_1 > 0$ such
that, for every     $n$-cylinder  $C_1$  and
$n$-cylinder or
$n$-gap $C_2$ with a common endpoint with $C_1$,
we have  $\kappa_1^{-1} < |C_1|_e / |C_2|_e <
\kappa_1$, where the lengths are
measured in any chart $(e,U)$ of the
atlas such that $C_1\cup C_2 \subset U$.
Hence there is $\kappa_2 > 0$ and $0<\nu <1$ such that
$|C|_e \le \kappa_2 \nu^n$ for every $n$-cylinder or $n$-gap $C$.
This property is equivalent to
the Markov map $f_\i$ being
uniformly expanding in $\TT^\i$.

Since $f$ on
$\L$ along leaves has affine
extensions with respect to the
charts in $\cA(r^\i)$ and the basic
$\i$-bolonomies have $C^{1+\eta}$
extensions we get that the Markov
maps $\tau_\i$ also have $C^{1+\eta}$
extensions with respect to the
charts in $\cB(r^\i)$ for some
$\eta > 0$. Since $\cA(r^\i)$ has bounded geometry,
we obtain that $f_\i$ also has
\emph{bounded geometry} in  $\cB(r^\i)$.
Since, for every $g \in \cT(\fL )$,
the $C^{1+}$ lamination atlas $\cA(g,\rho_g)$
has bounded geometry we obtain that
the Markov
map $f_\i$   has $C^{1+\eta}$
extensions with respect to the
charts in   $\cB(g,\rho_g)$, for some $\eta >0$,
and has bounded geometry.

\subsection{Holonomy pseudo-groups on $\TT^\i$}
The elements $\theta_\i=\theta_{f,\i}$ of the
\emph{holonomy pseudo-group on} $\TT^\i$
are the mappings defined as
follows.  Suppose that $I$ and $J$
are $\i$-leaf segments and $h:I\to
J$ a holonomy.  Then it follows
from the definition of the
train-track $\TT^\i$ that the map
$\theta :\pi_{\TT^\i}(I) \to \pi_{\TT^\i} (J)$
given by $\theta
(\pi_{\TT^\i}(x))=\pi_{\TT^\i}(h(x))$ is
well-defined.  The collection of
all such local mappings forms the
\emph{basic holonomy pseudo-group of} $\TT^\i$.
Note   that if $x$ is a junction
of $\TT^\i$ then there may be
segments $I$ and $J$ containing $x$
such that $I\cap J = \{ x \}$.  In this case the
image of $I$ and $J$ under the
holonomies will not agree in that
they will map $x$ differently.

\subsection{Markings on train-tracks}
\label{fhhhffgf}

For $\i=s$ and $u$, the Markov partition $\cR = \{
R_1,\ldots , R_m \}$ for $(f,\L)$ induces a Markov partition
$\cR^\i = \{ R_1^\i,\ldots , R_m^\i
\}$ for the Markov map $\tau=\tau_\i$ on
the train-track $\TT^\i$.
The marking $i:\Theta \to \Lambda$ determines   unique
\emph{markings} $i_u:\Theta^u \to \TT^u$
and  $i_s:\Theta^s \to \TT^s$
such that
$ i_u(w_0 w_1 \ldots)=\cap_{i\geq 0} R_{w_i}^u$
and
$i_s(\ldots w_{-1} w_0) =\cap_{i\geq 0} R_{w_i}^u$.
We note that
$
\pi_{\TT^\i} \circ i = i_\i \circ \pi_\i \ .
$
The map $i_\i$
is continuous, onto
$\TT^\i$ and semiconjugates the
shift map   on $\Theta^\i$
to the Markov map on $\TT^\i$.
Defining $\eps,\eps' \in \Theta^\i$
to be equivalent $(\eps
\overequiv{} \eps')$ if the point
$i^\i (\eps)=i^\i (\eps')$, we get that the
space $\Theta^\i / \overequiv{}$ is
homeomorphic to the train-track
$\TT^\i$.

Consider the Markov map $f_\i$
on  $\TT^\i$ induced by the action of
$f$ on $\ip$-leaves and described
above.  For $n \ge 1$, an   $n$-cylinder is the
projection into $\TT^\i$ of an
$\i$-leaf $n$-cylinder segment in $\L$.
Thus, each Markov rectangle in
$\Lambda$ projects in an unique
primary $\i$-leaf segment in
$\TT^\i$.
For $n \ge 1$,
an $n$-gap of $f_\i$ is the
projection into $\TT^\i$ of a
$\i$-leaf $n$-gap in $\L$.

We say that $\TT^\i$ is a \emph{no-gap train-track}
if $\TT^\i$ does not have gaps. Otherwise,
we call $\TT^\i$ a \emph{gap train-track}.

\subsection{Self-renormalisable structures}

\label{sect:ttsmooth}\label{foldjgjlke}

The $C^{1+}$ structure $\cS_\i$
on $\TT^\i$ is an
\emph{$\i$ self-renormalisable} if it has the following properties:
\begin{rlist}
\item
in this structure the
Markov mapping $m_\i$ is a
local diffeomorphism and has
bounded geometry in some $C^{1+}$ atlas of this structure; and
\item the
elements of the basic holonomy pseudo-group
are local diffeomorphisms
in $\cS_\i$.
\end{rlist}

We say that $\cB$ is a  $C^{1+}$ \emph{self-renormalisable atlas}
if $\cB$ has bounded geometry
and   extends to a $C^{1+}$ self-renormalisable structure.
By definition, a $C^{1+}$ self-renormalisable  structure
contains a  $C^{1+}$ self-renormalisable atlas.

A $C^{1+}$foliated $\i$-lamination atlas $\cA$ induces
a $C^{1+}$ $\i$ self-renormalisable atlas $\cB$
on $\TT^\i$ (and vice-versa) as follows:
The   holonomies are $C^{1+}$
with respect to the   atlas $\cA$
and so the  charts in
$\cB$
are $C^{1+}$ compatible, and the
basic holonomy pseudo-group of $\TT^\i$
are local diffeomorphisms.
Since $\cA$ has bounded geometry
the Markov mapping  $\tau_\i$ is a
local diffeomorphism and also has
bounded geometry in $\cB$.
Therefore,  $\cB$ is a  $C^{1+}$ self-renormalisable atlas
and extends to a   $C^{1+}$ self-renormalisable
structure $\cS(\cB)$  on $\TT^\i$.
Since  $\cA(r_\i)$ and $\cA^\i(g,\rho_g)$
are $C^{1+}$foliated $\i$-lamination atlas
we obtain that the atlases $\cB(r_\i)$ and $\cB^\i(g,\rho_g)$
determine respectively $C^{1+}$ self-renormalisable structures
$\cS(r_\i)$ and $\cS(g,\i)$.

\begin{lemma}
\label{fghfhttf}
The map $r_\i  \to \cS(r_\i)$ determines a one-to-one correspondence
between   $\i$-ratio functions
(or equivalently, $\i$-solenoid functions $r_\i|\sol^\i$)
and $C^{1+}$ self-renormalisable
structures on $\TT^\i$.
\end{lemma}

\begin{proof} Every ratio function $r_{\i}$
determines an unique $C^{1+}$ self-renormalisable $\Str$.
Conversely, let us prove that a given
$C^{1+}$ self-renormalisable
structure $\Str$ on $\TT^\i$ also determines an unique
ratio function $r_{\Str,\i}$.
Let $\cB$ be a bounded atlas for $\Str$.
Consider a small leaf segment $K$
and two leaf segments $I$ and $J$
contained in $K$.
Since the elements of the basic holonomy pseudo-group
on $\TT^\i$ are $C^{1+}$ and the Markov map
is also $C^{1+}$ and has bounded geometry we obtain
by Taylor's Theorem that the
following limit exists
\begin{eqnarray}\label{eqn:12133}
r_{\Str,\i} (I:J) & = & \lim_{n \to        \infty}
\frac {| \pi_\i f_\ip^{n}I|_{i_n}}{|\pi_\i f_\ip^{n} J|_{i_n}} \nonumber \\
& \in &
\frac{|\pi_\i I|_{i_0}}{|\pi_\i J|_{i_0}} \left(1 \pm \cO(|\pi_\i K|_{i_0}^\gamma) \right)
\ ,
\end{eqnarray}
where the size of the  leaf segments are
measured in   charts  of the
bounded atlas $\cB$. Furthermore, by   \cite{HR}  and \eqref{eqn:12133},
the charts in
$\cB(r_\i)$ and the charts in $\cB$ are $C^{1+}$ equivalent and so determine the
same  $C^{1+}$ self-renormalisable structure.
\end{proof}

\subsection{Hyperbolic diffeomorphisms}
\label{dfghnnvvfr}

Let  $g \in \cT(\fL )$ and
$\cA(g,\rho_g)$ be the
$C^{1+}$foliated $\i$-lamination atlas
determined by the Riemannian metric $\rho_g$.
As shown in  \S ~\ref{foldjgjlke}, the atlas  $\cA(g,\rho_g)$ induces a
$C^{1+}$ self-renormalisable atlas $\cB(g,\rho_g)$  on  $\TT^\i$
which generates a $C^{1+}$ self-renormalisable
structure $\cS(g,\i)$.

\begin{lemma}
\label{ghjjjyuu}
The mapping $g \rightarrow (\cS (g,s),\cS (g,u))$
gives  a 1-1 correspondence between
$C^{1+}$  conjugacy classes  in
$\cT(\fL )$
and pairs  $(\cS (g,s),\cS (g,u))$ of
$C^{1+}$ self-renormalisable
structures.
Furthermore, $r_{g,s} = r_{\cS (g,s),s}$ and
$r_{g,u} = r_{\cS (g,u),u}$.
\end{lemma}

\begin{proof}
By Lemma \ref{fghfhttf},
the pair $(\cS_{s},\cS_{u})$
determines a pair
$(r_{s,\Str}|\sol^s,r_{u,\Str}|\sol^u)$ of  solenoid  functions
and vice-versa.
By Corollary \ref{cor:sol}, the pair
$(r_{s,\Str}|\sol^s,r_{u,\Str}|\sol^u)$ determines an unique
$C^{1+}$ conjugacy class of diffeomorphisms $g \in \cT(\fL )$
which realise the pair
$(r_{s,\Str}|\sol^s,r_{u,\Str}|\sol^u)$
and vice-versa (and so  $(\cS (g,s),\cS (g,u))=(\cS_{s},\cS_{u})$).
Furthermore, by Lemma \ref{gdfrrn}, we get $r_{g,s} = r_{\cS (g,s),s}$ and
$r_{g,u} = r_{\cS (g,u),u}$.
\end{proof}

\section{Measure solenoid functions}

In this section,
we introduce   the following new concepts:
stable and unstable measure solenoid functions
and stable and unstable measure ratio functions.
Later, we will use the measure solenoid functions
and the measure ratio functions to determine
which Gibbs measures are $C^{1+}$-realisable by $C^{1+}$ hyperbolic
diffeomorphisms and by $C^{1+}$ self-renormalisable structures.

\label{gdghhnnhgr}

\subsection{Gibbs measures}
\label{ghhhnvcfer}

Let us give the  definition of an infinite two-sided  subshift of finite
type
$\Theta=\Theta(A)$. The elements of $\Theta$
are all   infinite two-sided
words $w=\ldots w_{-1}w_{0}w_{1}\ldots$
in the symbols $1,\ldots ,k$ such that $A_{w_{i}w_{i+1}}=1$, for all
$i \in \Bbb{Z}$.
Here $A=\left( A_{ij}\right) $ is any matrix
with entries $0$ and $1$ such that $A^{n}$ has all entries positive for some
$n\geq 1$.  We write $w\overequiv{n_1,n_2}w^{\prime }$ if
$w_j=w_j'$ for every $j=-n_1,\ldots ,n_2$. The metric $d$
on $\Theta$ is given by $d(w,w^{\prime })=2^{-n}$ if $n\geq 0$ is the
largest such that $w\overequiv{n,n}w^{\prime }$. Together with this metric
$\Theta$ is a compact metric space.
The two-sided shift map
$\tau:\Theta \to \Theta$   is the
mapping which sends $w=\ldots w_{-1}w_{0}w_{1}\ldots$
to $v=\ldots v_{-1}v_{0}v_{1}\ldots$ where $v_j=w_{j+1}$
for every $j \in \Bbb{Z}$.
An $(n_1,n_2)$-\emph{cylinder} $\Theta _{w_{-n_1}\ldots w_{n_2}}$, where
$w\in \Theta$, consists of all
those words $w^{\prime }$ in $\Theta$ such that
$w\overequiv{n_1,n_2}w^{\prime }$.
Let $\Theta^u$ be the set of all right-handed words $w_0 w_1 \ldots$
which extend to words $\ldots w_0 w_1 \ldots$ in $\Theta$,
and, similarly, let
$\Theta^s$ be the set of all left-handed words $\ldots w_{-1} w_0$
which extend to words $\ldots w_{-1} w_0   \ldots$ in $\Theta$.
Then $\pi_u:\Theta \to \Theta^u$
and $\pi_s:\Theta \to \Theta^s$ are
the natural projection given, respectively, by
$$\pi_u(\ldots w_{-1} w_0 w_1 \ldots) = w_0 w_1 \ldots
~~~~~{\rm and}~~~~~
\pi_s(\ldots w_{-1} w_0 w_1 \ldots) = \ldots w_{-1} w_0 \ .$$
An $n$-cylinder $\Theta^u _{w_{0}\ldots w_{n-1}}$
is equal to $\pi_u(\Theta _{w_{0}\ldots w_{n-1}})$
and an $n$-cylinder $\Theta^s _{w_{-(n-1)}\ldots w_{0}}$
is equal to $\pi_s(\Theta _{w_{-(n-1)}\ldots w_{0}})$.
Let  $\tau_u:\Theta^u \to \Theta^u$ and
$\tau_s:\Theta^s \to \Theta^s$ be
the corresponding one-sided shifts.

\begin{definition}
For $\i=s$ and $u$, we say that
$s_\i:\Theta^\i \to \Bbb{R}^+$ is an
\emph{$\i$-measure scaling function}
if $s_\i$ is a H\"older continuous function,
and for every $\xi \in \Theta^\i$
$$
\sum_{\tau_\i \eta= \xi} s_\i(\eta) = 1 \ ,
$$
where the sum is upon all $\xi \in \Theta^\i$ such that
$\tau_\i \eta= \xi$.
\end{definition}

For $\i \in \{s,u\}$, a $\tau$-invariant measure  $\nu$ on $\Theta$ determines a
unique  $\tau_\i$-invariant measure $\nu_\i=(\pi_\i)_* \nu$ on $\Theta^\i$.
We note that a $\tau_\i$-invariant measure $\nu_\i$ on $\Theta^\i$
has an unique $\tau$-invariant natural extension to an invariant measure $\nu$ on $\Theta$
such that $\nu(\Theta_{w_{0} \ldots w_{n_2}}) =
\nu_\i(\Theta^\i_{w_{0} \ldots w_{n_2}})$.

\begin{definition}
\label{fgdfgbbcc7}
A $\tau$-invariant measure   $\nu$ on $\Theta$ is a
\emph{Gibbs measure}:
\begin{rlist}
\item
if  the function
$s_{\nu,u}:\Theta^u \to \Bbb{R}^+$
given by
$$s_{\nu,u}(w_0 w_1 \ldots)=
\lim_{n \to \infty} \frac{\nu(\Theta_{w_0\dots w_n})}
{\nu(\Theta_{w_1\ldots w_n})} \ ,$$
is well-defined and it is an   $u$-measure scaling function; and
\item
if
the function
$s_{\nu,s}:\Theta^s \to \Bbb{R}^+$
given by
$$s_{\nu,s}(\ldots w_1 w_0)=
\lim_{n \to \infty}
\frac{\nu(\Theta_{w_n\dots w_0})}
{\nu(\Theta_{w_n\ldots w_1})} \ ,$$
is well-defined and it is a $s$-measure scaling function.
\end{rlist}
\end{definition}

The following theorem follows from
the results proved in
\cite{RPmeas}. It can also be
deduced from standard results about
Gibbs states such as those in
\cite{Bowen}.

\begin{theorem}
(Moduli space for Gibbs measures)
\label{fgdfgbbcc3}
Let $s_\i:\Theta^\i \to \Bbb{R}^+$ be an
\emph{$\i$-measure scaling function} for $\i=s$ or $u$.
Then there is an unique $\tau$-invariant Gibbs  measure $\nu$
such that $s_{\nu,\i}=s_\i$.
\end{theorem}

\subsection{Extended measure scaling function}
\label{rgergegr}

To present a   classification
of Gibbs measures  $C^{1+}$-Hausdorff which are realisable
by   codimension one attractors, we have to
define the cylinder-cylinder condition.
We will express the cylinder-cylinder condition,
in \S ~\ref{cylllll}, using the extended measure scaling
functions. These extended measure scaling
functions are  also used to present, in \S
~\ref{rgergegr}, the $\delta_\i$-bounded solenoid equivalence class  of
a Gibbs measure.

Throughout the paper, if $\xi \in\Theta^\ip$, we denote by $\xi_\L$
the  leaf primary cylinder segment $i(\pi_\ip^{-1} \xi) \subset \L$.
Similarly,
if $C$ is an $n_\i$-cylinder of $\Theta^\i$ then we denote
by $C_\L$ the      $(1,n_\i)$-rectangle $i(\pi_\i^{-1} C) \subset \L$.

\begin{figure}[tbp]
\centerline{\includegraphics[width=4.5cm]{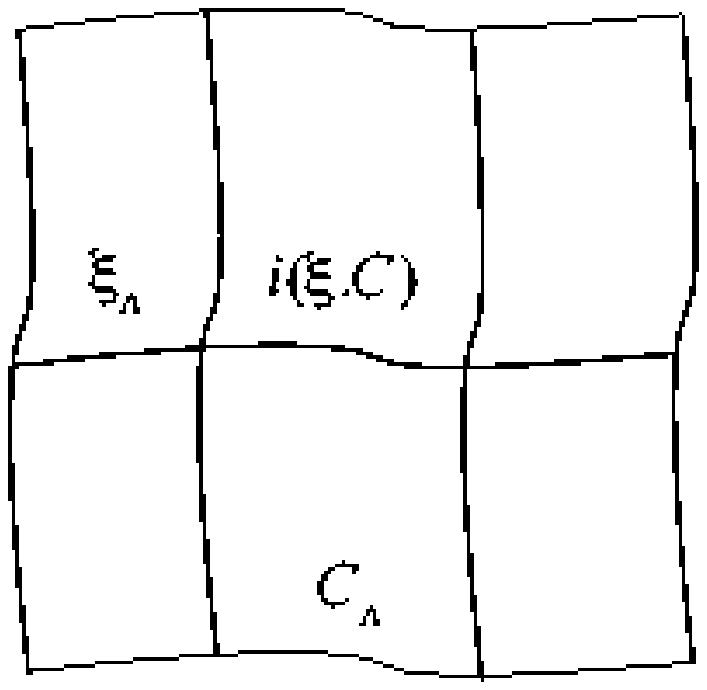} }
\caption{An $\i$-admissible pair $(\xi,C)$
where $\xi_\L = i(\pi_\ip^{-1} \xi)$ and
$C_\L = i(\pi_\i^{-1} C)$.}
\label{admissible_pair}
\end{figure}

We say that a $(n_1,n_2)$-cylinder  $\theta _{w_{-n_1}\ldots w_{n_2}}$
of $\Theta$
is an  \emph{$u$-symbolic leaf $n_2$-cylinder}
if $n_1=-\infty$.
Similarly, we say that a
$(n_1,n_2)$-cylinder  $\theta _{w_{-n_1}\ldots w_{n_2}}$
is a   \emph{$s$-symbolic leaf $n_1$-cylinder}
if $n_2=+\infty$.
Let $\xi \in \Theta^\ip$ and $C$ be  a
$n$-cylinder of $\Theta^\i$.
We say that the pair  $(\xi,C)$
is \emph{$\i$-admissible} if the set
$$\xi.C= \pi_\i^{-1} C \cap \pi_\ip^{-1} \xi $$
is non-empty (see Figure \ref{admissible_pair}).
We note that if the pair  $(\xi,C)$
is $\i$-admissible  then
$\xi.C$ is an $\i$-symbolic leaf $n$-cylinder,
and, conversely, any $\i$-symbolic leaf $n$-cylinder can be
expressed as $\xi.C$ where the pair  $(\xi,C)$
is $\i$-admissible.
The    set of all   {$\i$-admissible}  pairs  $(\xi,C)$
is  the $\i$-\emph{measure scaling set} $\msc^\i$.

Let
$C^u=\Theta^u_{w_0\ldots w_{n-1}}$ and
$C^s=\Theta^s_{w_{-(n-1)}\ldots w_0}$
be the $n$-cylinders  of
$\Theta^u$ and
of $\Theta^s$, respectively.
For $i<n$,  we denote
by $m^i C^u$ the \emph{$i$-th mother}
$\Theta^u_{w_0\ldots w_{(n-i-1)}}$ of $C^u$.
Similarly, we denote by $m^i C^s$ the
\emph{$i$-th mother} $\Theta^s_{w_{-(n-i-1)}\ldots w_0}$
of $C^s$.

\begin{figure}[tbp]
\centerline{\includegraphics[width=12cm]{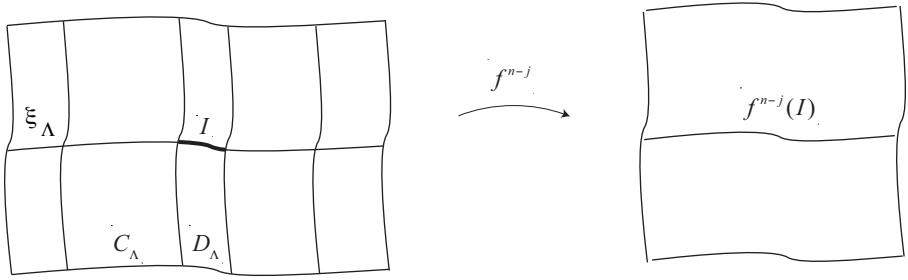} }
%\ \picill 141pt by 89pt (fig3_basich.eps) }
\caption{The
$(n-j+1)$-cylinder leaf segment $I = \xi_\L \cap D_\L$ and
the primary leaf segment
$f^{n-j}(I) = i(\pi_{\i'} \tau_{\i}^{n-j}(\xi . D))$,
where $D = m_{\i}^{j-1} C$.}
\label{alalp}
\end{figure}

Given an $\i$-measure scaling function  $s$,
we  construct the \emph{extended  $\i$-measure scaling function}
$\rho:\msc_\i \to  \Bbb{R}^+$ of $s$ as follows:
If $C$ is a $1$-cylinder then we define
$\rho_\xi(C)=1$.
If $C$ is a $n$-cylinder with $n \ge 2$, then we define
$$
\rho_\xi(C)=
\prod_{j=1}^{n-1}
s( \pi_\ip \tau_\i^{n-j} (\xi. m_\i^{j-1} C) )
$$
(see Figure \ref{alalp}).
We note that, if $s_{\nu,s}$ is the stable measure scaling function of a Gibbs measure $\nu$,
then
$$
\rho_\xi(C)=
\lim_{m \to \infty }\frac{\nu(\pi_s \circ \tau^m (\xi.C))}
{\nu(\pi_s \circ \tau^m (\xi))} \ .
$$
The unstable case is similar to the one above by taking $-m$ instead of $m$.
Hence, $\rho_\xi(C)$ is the ratio between the measure of $\xi. C$ and the measure of $\xi$  with respect to the
conditional measure determined by the Gibbs measure $\nu$ in $\xi$.

Recall that a $\tau$-invariant measure  $\nu$ on $\Theta$ determines a
unique  $\tau_u$-invariant measure $\nu_u=(\pi_u)_* \nu$ on $\Theta^u$
and a
unique  $\tau_s$-invariant measure $\nu_s=(\pi_s)_* \nu$ on $\Theta^s$.
The following theorem follows from
\cite{RPmeas}.

\begin{theorem}  (Ratio decomposition of a Gibbs measure)
\label{frgvnhgdfd}
Let
$\rho:\msc_\i \to  \Bbb{R}^+$
the an extended  $\i$-measure scaling function
and $\nu$ the corresponding  $\tau$-invariant Gibbs measure.
If $C$ is an $n$-cylinder  of $\Theta^\i$
     then
     $$                        \nu_\i (C) = \int_{\xi\in M}
     \rho_{\xi} (C)\nu_\ip (d\xi).
     $$
where $M= \pi_\ip \circ \pi_\i^{-1} C$ a $1$-cylinder  of $\Theta^\ip$.
Furthermore,  the ratios       $\nu_\i(C)/\rho_{\xi}(C)$ are
uniformly bounded away from $0$
and $\infty$.
\end{theorem}

\subsection{Cylinder-cylinder condition}
\label{cylllll}

We introduce the cylinder-cylinder condition that we will use to
classify all Gibbs measures that are $C^{1+}$-Hausdorff realizable by
codimension one attractors.

\begin{figure}[tbp]
\includegraphics[width=5in]{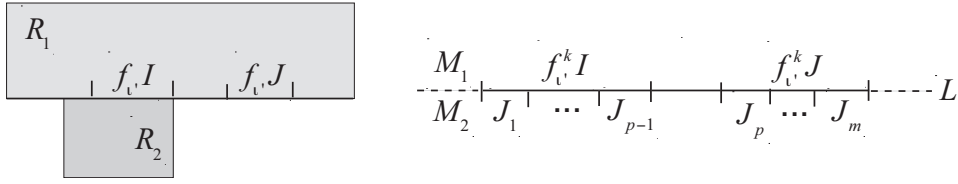}
\caption{The cylinder-cylinder condition for  $\i$-leaf segments.}
\label{ref:cylcyl}
\end{figure}

Similarly to the cylinder-gap condition given in \S ~\ref{dfsfbbvv}
for a given solenoid function, we are going to construct the
cylinder-cylinder condition for a given
measure solenoid function $\s_{\nu,\i}$.
Let $\delta_\i< 1$ and $\delta_\ip = 1$. Let  $(I,J) \in \Msol^\i$ be   such that
the $\i$-leaf segment
$f_\ip I \cup f_\ip J$  is contained in  an $\i$-boundary $K$ of a Markov rectangle $R_1$.
Then $f_\ip I \cup f_\ip J$ intersects another  Markov rectangle $R_2$.
Take the smallest $k \ge 0$ such that $f^k_\ip I \cup f^k_\ip J$
is contained in the intersection of the boundaries of two Markov rectangles
$M_1$ and $M_2$. Let $M_1$ be the Markov rectangle with the property
that $M_1 \cap f^k_\ip R_1$ is a rectangle with non empty interior,
and so $M_2 \cap f^k_\ip R_2$  has also non-empty interior.
Then, for some positive $n$, there are distinct  $\i$-leaf  $n$-cylinders
$J_1,\ldots,J_m$ contained in a primary cylinder $L$ of $M_2$ such that
$f^k_\ip I= \cup_{i=1}^{p-1} J_i$ and $f^k_\ip J= \cup_{i=p}^{m} J_i$.
Let $\eta \in \Theta^\ip$ be such that $\eta_\Lambda = L$
and,
for every $i=1,\ldots,m$, let
$D_i$ be a  cylinder  of $ \Theta^\i$ such that $i(\eta.D_i)=J_i$.
Let $\xi \in \Theta^\ip$ be such that $\xi_\Lambda = K$ and
$C_1$ and $C_2$ cylinders of $ \Theta^\i$ such that $i(\xi.C_1)=f_\ip I$
and $i(\xi.C_2)=f_\ip J$. We say that an $\i$-extended scaling function $\rho$
satisfies the  \emph{cylinder-cylinder  condition} (see Figure \ref{ref:cylcyl})  if
for all such leaf segments:
$$
\frac{\rho_{\xi}(C_2)}{\rho_{\xi}(C_1)}
=  \frac{ \sum_{i=p}^{m}  \rho_{\eta}(D_i)} { \sum_{i=1}^{p-1}  \rho_{\eta}(D_i)} \ .
$$

\subsection{Measure solenoid functions}
\label{fghbjhnnjw}

Let ${\Msol^\i}$ be the  set of all
pairs $(I,J)$ with the following properties:
(a)  If $\delta_\i=1$  then  $\Msol^\i=\sol^\i$. (b)
If $\delta_\i < 1$ then
$f_\ip I$ and   $f_\ip J$  are $\i$-leaf $2$-cylinders of a Markov rectangle $R$
such that $f_\ip I \cup f_\ip J$  is an $\i$-leaf segment,
i.e. there is an unique $\i$-leaf $2$-gap between them.
Let ${\msol^\i}$ be the  set of all
pairs $(I,J) \in \Msol^\i$ such that
the leaf segments $I$ and $J$  are not contained in an $\i$-global leaf
containing an $\i$-boundary of a Markov rectangle.
By construction, the set ${\msol^\i}$ is dense in $\Msol^\i$,
and   for every pair $(C,D) \subset {\msol^\i}$
there is an unique $\psi \in \Theta^\ip$
and an unique $\xi \in \Theta^\ip$
such that
$i(\pi_\ip^{-1}(\psi))=C$ and
$i(\pi_\ip^{-1}(\xi))=D$.
We will denote, in what follows, $i(\pi_\ip^{-1}(\psi))$ by $\psi_\L$
and $i(\pi_\ip^{-1}(\xi))$ by $\xi_\L$.

\begin{lemma}
\label{fghbhbfw}
Let $\nu$  be a Gibbs measure on $\Theta$.
The \emph{$s$-measure solenoid function  }
$\s_{\nu,s}:{\msol^s} \to \Bbb{R}^+$ of $\nu$
and the \emph{$u$-measure solenoid function  }
$\s_{\nu,u}:{\msol^u} \to \Bbb{R}^+$ of $\nu$
given by
$$
\s_{\nu,s}(\psi_\L,\xi_\L)=\lim_{n \to \infty}
\frac{\nu  ( \Theta_{\psi_0 \ldots \psi_n})}
{\nu ( \Theta_{\xi_0 \ldots  \xi_n})}
$$
and
$$
\s_{\nu,u}(\psi_\L,\xi_\L)=\lim_{n \to \infty}
\frac{\nu  ( \Theta_{\psi_n \ldots \psi_0})}
{\nu ( \Theta_{\xi_n  \ldots  \xi_0 })}
$$
are both well-defined.
\end{lemma}

\begin{proof}  Let $(I,J) \in \msol^\i$. By Property (iii) of $\msol^\i$, there is
$k=k(I,J)$ such that $f_\ip^k I$ and $f_\ip^k J$
are cylinders with the same mother
$m f _\ip^k I=m f_\ip^k J$.
Let $(\xi:C)$ and $(\xi:D)$ be the admissible pairs in
$\msc_\i$ such that $i(\xi.C)=f_\ip^k I$ and
$i(\xi.D)=f_\ip^k J$.
Since  the measure $\nu$ is $\tau$-invariant, we obtain that
$$
\s_{\nu,\i}(I,J) = \rho_{\xi} (C) \rho_{\xi} (D)^{-1}  \ ,
$$
where $\rho$ is the extended scaling function determined by the Gibbs measure $\nu$.
Therefore,   the $\i$-measure solenoid function $\s_{\nu,\i}$
is well-defined for $\i\in\{s,u\}$.
\end{proof}

\begin{lemma}
\label{fghbhbfw3}
Suppose $\delta_{f,\i}=1$.
If   an $\i$-measure solenoid function
$\s_{\nu,\i}:{\msol^\i} \to \Bbb{R}^+$ has
a continuous extension to $\sol^\i$
then its extension
satisfies the matching condition.
\end{lemma}

\begin{proof}
Let $(J_0,J_1) \in \sol^\i$ be  a pair of primary cylinders
and  suppose that we
have pairs $$(I_0,I_1),
(I_1,I_2),\ldots ,
(I_{n-2},I_{n-1})\in \sol^\i$$
of primary cylinders such that
$f_\i J_0=\union_{j=0}^{k-1}  I_j$ and
$f_\i J_1=\union_{j=k}^{n-1}  I_j$.
Since the set  $\msol^\i$ is dense in $\sol^\i$
there are pairs
$(J^l_0 ,J^l_1)  \in \msol^\i$
and pairs     $(I^l_j ,I^l_{j+1})$
with the following properties:
\begin{rlist}
\item
$f_\i J^l_0 =\union_{j=0}^{k-1}  I^l_j$ and
$f_\i J^l_1=\union_{j=k}^{n-1}  I^i_j$.
\item
The  pair $(J^l_0,J^l_1)$
converges to $(J_0,J_1)$ when $i$ tends to infinity.
\end{rlist}
Therefore, for every $j=0,\ldots,n-2$ the pair
$(I^l_j,I^l_{j+1})$
converges to $(I_j,I_{j+1})$ when $i$ tends to infinity.
Since $\nu$ is a $\tau$-invariant measure, we get that
the matching condition
$$
\s_{\nu,\i} (J^l_0:J^l_1) =
\frac{1+\sum_{j=1}^{k-1}\prod_{i=1}^{j}\s_{\nu,\i}
(I^l_j:I^l_{i-1})}{\sum_{j=k}^{n-1}
\prod_{i=1}^{j}\s_{\nu,\i}
(I^l_j:I^l_{i-1})}
$$
is satisfied for every $l \ge 1$.
Since the extension of $\s_{\nu,\i}:\msol^\i   \to \Bbb{R}^+$
to the set $\sol^\i$ is continuous,
we get that the matching condition also holds
for the   pairs
$(J_0,J_1)$ and  $(I_0,I_1), \ldots, (I_{n-2},I_{n-1})$.
\end{proof}

\begin{remark}
We say that an $\i$-measure solenoid function
$\s_{\nu,\i}$ of a Gibbs measure $\nu$ satisfies the cylinder-cylinder condition
if the extended scaling function of the Gibbs measure $\nu$ satisfies the cylinder-cylinder condition.
\end{remark}

\subsection{Measure ratio functions}
\label{fsefsefff}

We say that $\rho$ is a {\it
$\i$-measure ratio
function}        if
\begin{rlist}
\item $\rho(I:J)$ is well-defined for every pair of $\i$-leaf
segments $I$ and $J$ such that (a) there is
an $\i$-leaf segment $K$ such that  $I,J\subset K$,
and (b) $I$ or $J$ has non-empty interior;
\item
if $I$ is an $\i$-leaf gap then $\rho(I:J)=0$ (and $\rho(J:I)=+\infty$);
\item
if $I$ and $J$ have non-empty interiors then $\rho(I:J)$ is strictly positive;
\item
$\rho(I:J) = \rho(J:I)^{-1}$;
\item  if  $I_1$ and $I_2$ intersect
at most in one of their endpoints then
$\rho(I_1\cup I_2 :K) = \rho(I_1:K)+\rho(I_2:K)$;
\item            $\rho$ is invariant under $f$, i.e.
$\rho(I:J)=r(f I: f J)$ for all $\i$-leaf segments;
\item  for every basic
$\i$-holonomy map $\theta:I\to J$
between the leaf segment
$I$ and the leaf segment
$J$  defined  with respect to a rectangle $R$
and for every
$\i$-leaf segment $I_0\subset I$ and
every $\i$-leaf segment or gap
$I_1\subset I$,
\begin{equation}
     \label{3eq:Holderness_of_ratios}
     \left| \log \frac{\rho(\theta
     I_0 :\theta I_1)}
{\rho(I_0:I_1)} \right|
     \le \cO \left(\left( d_\Lambda
     (I,J
     )\right)^\epsilon \right)
\end{equation}
where $\epsilon\in (0,1)$ depends upon $\rho$ and
the constant of proportionality
also depends upon   $R$, but not on
the segments considered.
\end{rlist}

We note that if $\TT^\i$ is a no-gap train-track then
an  $\i$-measure
ratio function is  an $\i$-ratio function.

\begin{remark}
A function $\sigma: \msol^\i \to \reals^+$ that has an H\"older
continuous extension to $\Msol^\i$ determines a unique extended scaling function $\rho$,
and so we say that $\sigma$ satisfies
the   cylinder-cylinder  condition if the extended scaling function $\rho$
satisfies the cylinder-cylinder condition.
\end{remark}

Let ${\mathcal SOL}^\i$ be the space  of all
H\"older continuous functions $\sigma_\i: \Msol^\i \to \reals^+$
with the following properties:
\begin{rlist}
\item
If  $\TT^\i$ is a no-gap train-track  then $\sigma_\i$ is an $\i$-solenoid function.
\item
If  $\TT^\i$ is a gap train-track  and
$\TT^\ip$ is a no-gap train-track  then  $\sigma_\i$
satisfies the   cylinder-cylinder condition.
\item  If   $\TT^\i$ and
$\TT^\ip$ are no-gap train-tracks  then  $\sigma_\i$
does not have to satisfy any extra property.
\end{rlist}

\begin{lemma}
\label{altum}
The map $\rat \to \rat|\Msol^\i$ determines an one-to-one
correspondence between $\i$-measure ratio functions
and functions contained in ${\mathcal SOL}^\i$.
\end{lemma}

\begin{proof}
The proof follows similarly to the proof of Lemma \ref{gdfrrn}.
\end{proof}

\begin{remark}
\label{prator}
\begin{enumerate}
\item[(i)]
By Lemma \ref{altum},
a Gibbs measure $\nu$ with an $\i$-measure solenoid function with
an extension $\hat{\sigma}$ to
$\Msol^\i$ such that $\hat{\sigma} \in {\mathcal SOL}^\i$
determines an unique $\i$-measure ratio function $\rho_\nu$.
\item[(ii)]
A measure ratio function $\rho$ determines naturally a measure scaling function,
and so, by Lemma \ref{fgdfgbbcc3}, a Gibbs measure $\nu_\rho$.
\item[(iii)]
By Lemma \ref{altum},
a  function $\sigma: \msol^\i \to \reals^+$  with
an extension $\hat{\sigma}$ to
$\Msol^\i$ such that $\hat{\sigma} \in {\mathcal SOL}^\i$
determines an $\i$-measure ratio function,  and, by (ii), a unique
Gibbs measure $\nu$ such that $\sigma = \sigma_{\nu}$.
\end{enumerate}
\end{remark}

\section{Natural geometric measures}
\label{dfgdgdbbbw}

In this section, we define
the natural geometric
measures $\mu_{\cS,\delta}$ associated with a
self-renormalisable structure
$\cS$ and $\delta >0$.
The natural geometric
measures are
measures determined by
the length scaling structure of
the cylinders. We will prove that
every  natural geometric measure is a pushforward of a
Gibbs measure with the property that the
measure solenoid function determines a  measure ratio function.
In \S ~\ref{fgrrffghtrd}, we will show that a Gibbs measure
with the property that its
measure solenoid function determines a measure ratio function
is $C^{1+}$-realisable by a self-renormalisable structure.

\begin{definition}
\label{fgtyjtyjyu1}
Let $\cS$ be a $C^{1+}$ self-renormalisable structure
on $\TT^{\i}$.
If $\TT^{\i}$ is a gap train-track let $0<\delta <1$,
and if $\TT^{\i}$ is a no-gap train-track let $\delta =1$.
\begin{rlist}
\item
We say that $\cS$ has a \emph{natural geometric measure
$\mu_\i=\mu_{\cS,\delta}$ with pressure} $P=P(\cS,\delta)$
if (a) $\mu_\i$ is a $f_\i$-invariant measure;
(b), there exists $\kappa > 1$ such that
for all $n \ge 1$ and all $n$-cylinders $I$  of  $\TT^\i$,
we have
\begin{equation}
\label{dcsfere}
\kappa^{-1} <\frac{\mu_\i(I)} {|I|_i^\delta        e^{-nP}} < \kappa \ ,
\end{equation}
where $i$ is a chart
containing $I$ of a bounded atlas $\cB$ of $\cS$;
\item
We say that  $\cS$
is a $C^{1+}$
\emph{realisation of a Gibbs measure} $\nu=\nu_{\cS,\delta}$ if
$\mu_\i= (i_\i)_* \nu_\i$ where  $\nu_\i=(\pi_\i)_*\nu$ and
$\mu_\i=\mu_{\cS,\delta}$  is a natural geometric measure
of $\cS$.
\end{rlist}
\end{definition}

Suppose that we have a $C^{1+}$
self-renormalisable  structure $\Str$ on $\TT^\i$ and
that $\cB$ is a bounded atlas for
it.  Let $\d > 0$.  If $I$ is a
segment in $\TT^\i$, let $|I|=|I|_i$ be
its length in any chart $i$ of this
atlas which contains it. If $C$
is a $m$-cylinder,  let us denote $m$ by $n(C)$
and  $i_\i(C)$ by $I_{C}$.
For $m_1 \ge 1$ and $m_2 \ge 1$,
let $C$ be an $m_1$-cylinder and $D$ an $m_2$-cylinder
contained in the same $1$-cylinder.
Let
\begin{equation}
\label{dfeeegh}
L_{\delta,s}(C:D) =
     \frac { \sum_{C' \subset C}
     |I_{C'}|^\delta  e^{-n(C')s } } {
     \sum_{D' \subset D} |I_{D'}|^\delta
     e^{-n(D')s }}
\end{equation}
where the sums are respectively
over all cylinders contained in
$C$ and $D$ and the values $|I_{C'}|$ and $|I_{D'}|$
are determined using  the same chart in $\cB$.
Let the \emph{pressure} $P=P(\cS,\delta)$ be the infimum  value of $s$ for which
the numerator (and the denominator) are finite.

If $\xi \in\Theta^\ip$, then
the  leaf $1$-cylinder segment $\xi_\L=i(\pi_\ip^{-1} \xi) \subset \L$
is also   regarded, without ambiguity, as a point in the train-track
$\TT^\ip$. Similarly,
if $C$ is an $n$-cylinder of $\Theta^\i$ then
the $(1,n)$-rectangle $C_\L=i(\pi_\i^{-1} \xi) \subset  \L$
is also   regarded, without ambiguity, as an $n$-cylinder
of  the train-track $\TT^\i$.

The following theorem follows from
the results proved in
\cite{RPmeas}. It can also be
deduced from standard results about
Gibbs states such as those in
\cite{Bowen}.

\begin{theorem}\label{thm:SRB}
Let $\cS$ be a $C^{1+}$ self-renormalisable structure
on $\TT^{\i}$. For every $\delta > 0$,
there is a unique geometric natural
measure $\mu_\i=\mu_{\cS,\delta}$
with pressure $P=P(\cS,\delta) \in \reals$,
and there is an   unique $\tau$-invariant Gibbs measure $\nu=\nu_{\cS,\delta}$
on $\Theta$ such that $\mu_\i=(i_\i)_* \nu_\i$ where $\nu_\i=(\pi_\i)_* \nu$.
Furthermore, the measure  $\mu_\i$
has the following properties:
\begin{rlist}
      \item There is $0<\a<1$
     such that if $C$ and $D$
     are any two $n$-cylinders in
$\Theta^{\i}$
such that $I_C$ and $I_D$ are contained
in a common segment $K$
       then
     $$\frac{\mu_\i(I_{C})}{\mu_\i(I_{D})}
     \in
     \left(1\pm\cO(|K|^\a)\right)
     L_{\delta,P} (C:D) \ .$$
\item  If $\rho:\msol^\i \to \reals$
is the extended measure scaling function of $\nu_\i$,
then
$$
\rho_\xi(C)=
\lim_{m\to \infty} L_{\delta,P}  (C_m:\xi_m) \ ,
$$
where $C_m$ and $\xi_m$ are the cylinders given by
$I_{C_m} = f^{m}_\ip (C_\L \cap \xi_\L)$ and
$I_{\xi_m}=f^{m-1}_\ip \xi_\L$.
\item ({\it ratio decomposition})  if $C$ is an
     $n$-cylinder in $\Theta^{\i}$
     and $C_{p}$ is the primary
     cylinder containing $C$
     then
     \begin{equation} \label{eqn:ratio_decomp3}
             \mu_\i (I_{C}) =
           \int_{\xi \in \pi_\ip(C) }
     \rat_{\xi} (C) \nu_\ip (d\xi)  \ .
     \end{equation}
\end{rlist}
\end{theorem}

\begin{figure}[tbp]
\centerline{\includegraphics[width=9.5cm]{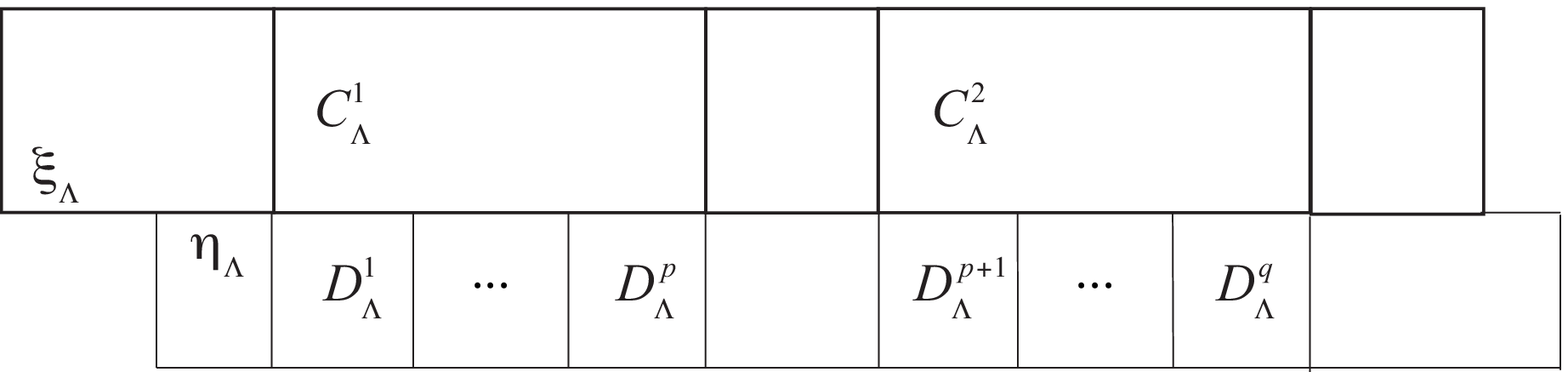} }
%\ \picill 141pt by 89pt (fig3_basich.eps) }
\caption{The rectangles $C_\L^1,C_\L^2$ and  $D_\L^1, \ldots, D_\L^q$}
\label{primary_leaf}
\end{figure}

\begin{lemma}\label{rem:close_ratios}
Let $\cS$ be a $C^{1+}$ self-renormalisable structure
on $\TT^{\i}$ and let $\rho$ be
the extended measure scaling function of
the Gibbs measure $\nu_{\cS,\delta}$.
\begin{rlist}
\item
If  $C$ and $D$ are two
cylinders contained in an
$n$-cylinder $E$  of $\Theta^\i$
then, for all
$\xi,\eta$  contained in the $1$-cylinder
$\pi_\ip (\pi_\i^{-1} E)$ of $\Theta^\ip$,
\begin{equation}
\label{fdfdgghbnn}
\frac{\rat_{\eta}(C)}{\rat_{\eta}(D)}
\in \left( 1\pm
\cO(\theta^{n})\right)
\frac{\rat_{\xi}(C)}{\rat_{\xi}(D)} .
\end{equation}
\item
Let $\TT^\ip$ be a no-gap train-track. Let $\xi, \eta \in \Theta^\ip$
be such that the corresponding leaf segments in $\L$
have a common intersection $K$ (or coincide).
Let   $(\xi:C^1), (\xi:C^2), (\eta:D^1),\ldots,(\eta:D^q)$
be admissible pairwise distinct pairs in $\msc^\i$ such that
(a) $\xi_\L \cap C^1_\L = \xi_\L \cap (\cup_{i=1}^p D^i_\L) \subset K$, and
(b) $\xi_\L \cap C^2_\L = \xi_\L \cap (\cup_{i=p+1}^q D^i_\L) \subset K$
(see Figure \ref{primary_leaf}).
Then
\begin{equation}
\label{fdgeghnn}
\frac{\rat_{\xi}(C^1)}{\rat_{\xi}(C^2)}=
\frac{\sum_{i=1}^p  \rat_{\eta}(D^i)}
{\sum_{i=p+1}^q \rat_{\eta}(D^i)} \ .
\end{equation}
\item
Let $\TT^\i$ be a no-gap train-track (and $\delta=1$).
Then for every admissible pair $(C:\xi) \in \msc^\i$
we get
\begin{equation}
\label{ewqw}
\rat_{\xi}(C)=r^{\i}
(C_\L \cap    \xi_\L : \xi_\L)
\end{equation}
where $r^{\i}$  is the
$\i$-ratio function determined by
the $C^{1+}$
self-renormalisable structure.
\end{rlist}
\end{lemma}

\begin{proof}
\emph{Proof of (i) and (ii)}.
Suppose that $C$ and $D$ are two
cylinders contained in an
$n$-cylinder $E$. Let $E_{1}$ be
a $(n+1)$-cylinder whose image
under the shift map $\tau$ is
   $E$  and let $C_{1}$ and $D_{1}$
be the cylinders in $E_{1}$ such
that $\tau C_{1}=C$ and $\tau
D_{1}=D$. Then
$$
L_{\delta,P }(C_{1}:D_{1})\in \left(
1\pm \cO(\theta^{n})\right)
L_{\delta,P}(C:D)
$$
where (i) $0<\theta <1$ is independent
of $C$, $D$, $E$ and $E_{1}$, and
$P=P(\cS,\delta)$ is the pressure. This
follows directly from the
definition of $L_{\delta,P }$ together
with the fact that for all
cylinders $C'$, $D'$ in $E_1$,
$$
\frac{|I_{D'}|}{|I_{C'}|} \in \left(
1\pm \cO(\theta^{n})\right)
\frac{|I_{\tau D'}|}{|I_{\tau C'}|}.
$$
As a corollary of this we deduce
\eqref{fdfdgghbnn}.
Then, equality \eqref{fdgeghnn} follows from using
that the   local holonomies are local diffeomorphisms in the
self-renormalisable structure of $\TT^\i$.

\emph{Proof of (iii)}. In this case  the
self-renormalisable structure
$\cS$ is a local manifold
structure as defined in \S
~\ref{sect:ttsmooth} (i.e.\ the
charts are homeomorphisms onto
a subinterval of $\reals$),
and $\delta=1$.
Using \eqref{dfeeegh}, we get
$P(\cS,\delta)=0$ and so  the ratios
$\mu(I)/|I|$ are uniformly
bounded away from $0$ and
$\infty$ for all segments $I$
in $\TT^{\i}$.  Moreover, in
this case, the length system
$\l$ matches in the sense that,
if $C$
is an $n$-cylinder then
$\sum_{C'} |I_{C'}| =|I_C|$
where the sum is over all $m$-cylinders
$C'$ contained in $C$ and $|I_C|$ and $|I_{C'}|$ are
obtained using the
same chart in $\cB$.  Thus,
if $C$ and $D$ are $n$-cylinders and
$I_C \cup I_D$ is a segment of $\TT^\i$
then
$$                \frac{\mu_\i(I_{C})}{\mu_\i(I_{D})}
     \in(1\pm\cO(\theta^{n}))
     \frac{|I_{C}|}{|I_{D}|}.
$$
Hence,
$$
     \rat_{\xi}(C) =
     \lim_{m\to\infty}
     \frac{|f^{m}_\ip (C_\L \cap \xi_\L)|}
       {|f^{m}_\ip \xi_\L|}
$$
which implies \eqref{ewqw}.
\end{proof}

\begin{lemma}\label{cor_geommeas}
If $\d$ is the Hausdorff
dimension of $\TT^{\i}$ then
the ratios
$\mu_\i(I_{C})/|I_{C}|^\d$ are
uniformly bounded away from $0$
and $\infty$.  It follows from
this that the Hausdorff
$\delta$-measure $\cH^\d$ is
finite and positive on
$\TT^{\i}$ and such that
$\mu_\i$ is absolutely continuous
with respect to $\cH^\d$.
\end{lemma}

\begin{proof}
Since in this case $\d$ is the
Hausdorff dimension of $\TT^{\i}$,
$P_{\d}=0$.
Now, let us prove that  the
ratios $\mu_\i(I)/|I|^\d$ are uniformly
bounded away from $0$ and $\infty$
for all segments $I$.
Suppose that I
is any segment in $\TT^{\i}$.
Then either there exists a
cylinder $C$ with $I_{C}\subset
I$ such that $I\subset I_{mC}$
or there exist cylinders
$I_{C}, I_{D}\subset I$ with a common endpoint
such         that $I\subset I_{mC}\cup I_{mD}$.
In the first case let
       $\II=I_C$ and $m\II=I_{mC}$,  otherwise let
$\II=I_C\cup I_D$ and
$m\II=I_{mC}\cup I_{mD}$.
Since by bounded geometry there exists $\sigma >0$ such
that for all cylinders
$|I_C|/|I_{mC}|\in
[\sigma,\sigma^{-1}]$, $$\sigma |m\II |\leq
|\II |\leq |I| \leq |m\II |\leq
\sigma^{-1} |\II |.$$ By Theorem
\ref{thm:SRB}, the ratios $\mu_\i(I_C)/|I_C|^\delta$
are uniformly bounded away from
$0$ and $\infty$ for cylinders
$C$, and so  the
ratios   $\mu_\i(I)/|I|^\delta$ for segments $I$ are
also uniformly bounded away from
$0$ and $\infty$.

Suppose that $A$ is a  subset of
$\TT^\i$.  Recall the definition of
$\cH^\d_\eps (A)$ as the infimum of
the sums $\sum r_i^\d$ where the
$r_i$ are the lengths of the
segments of an $\eps$-cover of $A$.
Then the Hausdorff $\d$-dimensional
outer measure of $A$ is $\cH^\d
(A)=\lim_{\eps\to 0}\cH^\d_\eps (A)
= \sup_{\eps > 0} \cH^\d_\eps (A)$.
Now suppose that $A$ is a subset of
$\TT^\i$ and $\{ B_i \}$ is a cover
of $A$ by segments of length $r_i$.
Then, $\cH^\d (A)\geq\cO (\mu_\i(A))$
because
$$
\cH^\d (A) \geq        \sum r_i^\d \geq \cO \left ( \sum
\mu_\i(B_i) \right ) \geq \cO \left ( \mu_\i \left (\union
B_i \right ) \right )\geq \cO ( \mu_\i(A)).
$$
Now suppose that $A$ is a  Borel
subset of $\TT^i$ with $\mu_\i (A) >0$.
The set of segments in
$\TT^\i$ is a Vitali class for $A$.
Therefore, by the Vitali Covering
Theorem (e.g.\
\cite{Falconer,Federer}), given
$\eps > 0$, there is a countable
disjoint sequence $B_{j}$ of
segments such that either $\sum_{j}
|B_{j}|^{\delta} = \infty$ or $
\cH^\d (A) \leq
\sum_{j}|B_{j}|^{\delta} + \eps $.
But since the $B_{j}$'s are
disjoint we get
$$\sum_{j}|B_{j}|^{\delta} \leq \cO
\left( \sum_{j} \mu_\i (B_{j}) \right) \leq \cO
(\mu_\i (A) ) \ .$$
Thus $ \cH^\d (A) \leq
\cO (\mu_\i (A) )+\eps$. Letting $\eps \to
0$ gives $\cH^\d (A) \leq
\cO (\mu_\i (A) ).$
This proves that the measure $\mu_\i$
is proportional to $\cH^\d$. It
follows immediately that $\cH^\d$ is
positive and finite.
\end{proof}

\subsection{Measure ratio functions}
\label{sect:2dmeasures}

In this subsection,
we prove
that,  for  every $\delta >0$,
a given $C^{1+}$ self-renormalisable structure  $\cS$ on $\TT^\ip$
determines an $\i$-measure ratio function $\rat_{\cS,\delta}$
such that
the Gibbs measure $\nu_\rho$ determined by $\rat_{\cS,\delta}$
(see Remark \ref{prator})
is the same as the Gibbs measure
$\nu_{\cS,\delta}$ which is $C^{1+}$ realisable by the self-renormalisable structure  $\cS$.

\begin{figure}[tbp]
\centerline{\includegraphics[width=8.5cm]{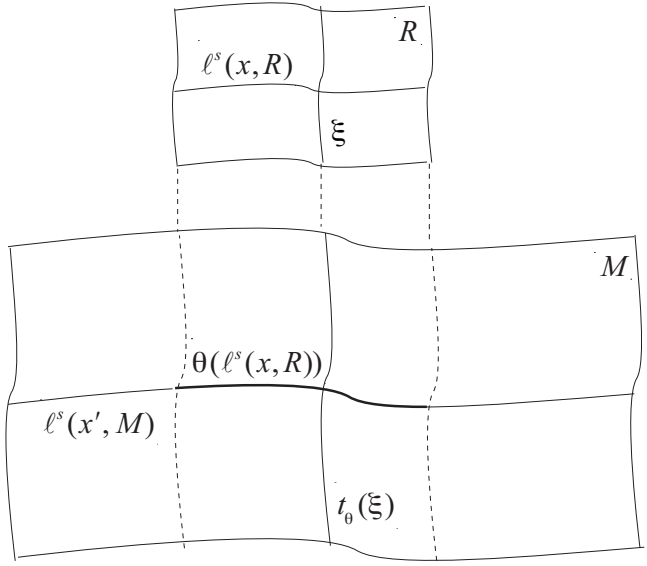} }
%\ \picill 170pt by 117pt (fig2_basicholos.eps) }
\caption{The holonomy injection $t_{\theta}$.}
\label{holonomy_injection}
\end{figure}

Let $\xi \in  \TT^\ip$ be an  $\i$-leaf
segment spanning of a Markov rectangle
$M$. Let $R$ be a
rectangle inside $M$.  There are cylinders
$C_j \in \Theta^\i$ such that
$\pi_{\i}R$ is  the
countable (or finite) union  $\cup_{j \in {\rm Ind}} I_{C_{j}}$ of cylinders
$I_{C_{j}}=i_\i(C_j)$ of
$\TT^{\i}$, and  any two of which
intersect at most in a point of
their boundary (see Figure \ref{holonomy_injection}).
Suppose that $\xi  \cap R \neq \emptyset$. Let  $\xi' \in \Theta^\ip$
be  such  that  $i_\i (\xi')=\xi$ (we note that $\xi'$ might not be uniquely determined).
Using  \eqref{fdgeghnn}, the following ratio is well-defined
\begin{equation}
\label{eqn:rhos}
\rho_{\i,\xi}(R:M)  =
\sum_{j \in {\rm Ind}}\rho_{\xi'}(C_{j}) \ ,
\end{equation}
where $\rho_{\xi'}(C_{j})$ is the $\i$-extended measure scaling function
of $\nu_{\cS,\delta}$.
If $\xi  \cap R = \emptyset$ then
we define $\rho_{\xi'}(R:M) =0$.
More generally, suppose
that $R_0$ and $R_1$ are $\ip$-spanning
rectangles contained in $R$.
Then we define
$$\rho_{\i,\xi}(R_0:R_1) = \rho_{\i,\xi}(R_0:M)
\rho_{\i,\xi}(R_1:M)^{-1} \ .$$

\begin{theorem}
(2-dimensional ratio decomposition)
Let $\cS$ be a $C^{1+}$ self-renormalisable structure
and $\mu_\i=\mu_{\cS,\delta}$  a   natural geometric  measure
for some $\delta >0$.
\label{thm:ratiodecomp_general}
Suppose that $R$ is a rectangle
contained in a Markov rectangle
$M$.
Then
\begin{equation}
     \label{eqn:gen_thm}
\mu (R) =
\int_{\pi_{\TT^\ip}(R)}
\rat_{\i,\xi}
(R:M) \mu_\ip (d\xi ) \  .
\end{equation}
\end{theorem}

We now consider the case where
$\TT^\i$ is a no-gap train-track.
Let $\cS$ be a $C^{1+}$ self-renormalisable   structure
and  $\mu_\i=\mu_{\cS,1}$   the natural measure
(with pressure $P=0$).
Recall the definition of
$t_{R}^{\ip}$ as the set of
spanning $\i$-leaf segments of the
rectangle $R$ (not necessarily a Markov rectangle). By
the local product structure, one
can identify $t_{R}^{\ip}$ with
any spanning $\ip$-leaf segment
$\l^{\ip}(x,R)$ of $R$.
Suppose that $R$ is a rectangle  and $M$ is a Markov
rectangle and that $\theta : \l=
\l^{\ip}(x,R)\to
\l' \subset\l^{\ip}(x',M)$ is a basic
holonomy defined on the spanning
$\ip$-leaf segment $\l$.  This
defines an injection
$t_{\theta}:t_{R}^{\ip}\to
t_{M}^{\ip}$ which we call the
\emph{holonomy injection} induced
by $\theta$ (see Figure \ref{holonomy_injection}).
The measure
$\mu_\ip$ on $\TT^{\ip}$ induces a
measure on
$t_{M}^{\ip}$ which we
can pull back to
$t_{R}^{\ip}$ using
$t_{\theta}$ to obtain a measure
$\mu_{R,M}^{\theta}$ i.e.
$\mu_{R,M}^{\theta} (E) =
\mu_\ip (\pi_{\TT^\ip}(t_{\theta}(E)))$.

\begin{theorem}(2-dimensional
ratio decomposition for SRB
measures)
\label{thm:ratiodecomp_SRB}
Let
$\TT^\i$ be a no-gap train-track.
Let $\cS$ be a $C^{1+}$ self-renormalisable   structure
and  $\mu_\i=\mu_{\cS,1}$   the natural measure
(with pressure $P=0$).
If $t_{\theta}:t_{R}^{\ip}\to
t_{M}^{\ip}$ is a holonomy
injection as above with $P$ a
Markov rectangle then
\begin{equation}
     \label{eqn:gen_SRB}
\mu (R) = \int_{t_{R}^{\ip}}
r_\i (\xi:t_{\theta}(\xi))
\mu_{R,M}^{\theta} (d\xi ) \ ,
\end{equation}
where $r_{\i}$ is the
$\i$-ratio function determined by $\cS$.
\end{theorem}

\begin{remark}
Note that if $R\subset M$ then
$t_{\theta}(\xi)$ is just the
$M$-spanning $\i$-leaf
containing $\xi$ and
$\mu_{R,M}^{\theta}=\mu_\ip$.
\end{remark}

Since any rectangle can be
written as the union of
rectangles $R$ with the
property hypothesised in the
theorem for some Markov
rectangle, the above theorem gives an
explicit formula for the
measure of any rectangle in
terms of a ratio decomposition using the
ratio function which characterises the smooth structure of
the train-track.

\noindent
\emph{Proof of Theorems
\ref{thm:ratiodecomp_general} and
\ref{thm:ratiodecomp_SRB}.} Suppose
that $R$ is any rectangle, $M$ is a
Markov rectangle and
$t_{\theta}:t_{R}^{\ip}\to
t_{M}^{\ip}$ is a holonomy
injection as above (in the case of
Theorem \ref{thm:ratiodecomp_general}
$t_{\theta}$ is the identity map).
Then  we note that there is
$0<\nu<1$ such that for all $n>0$
we can write $R=R_{0}\cup\dots\cup R_{N(n)}$  where
\begin{rlist}
\item
$R_{0},\dots ,R_{N(n)}$
are rectangles  which intersect
at most in their boundary
leaves and their spanning
$\i$-leaf segments are also $R$-spanning $\i$-leaf segments;
\item $P_i= t_{\theta} R_i$ and $\pi_\ip(P_{i})$ is an $n$-cylinder
of $\TT^\ip$ for every $0 \leq i \leq
N(n)$;
\item $R_0$ is the empty set, or  $\pi_\ip(P_0)$ is
strictly contained in an $n$-cylinder
of $\TT^\ip$, and so, using the bounded geometry of the Markov map
(see \S ~\ref{gfgrtr}) and \eqref{dcsfere},
$\mu(R_{0})<\cO(\epsilon_0^{n})$ for some $0<\epsilon_0<1$;
\end{rlist}
Let  $S_{i}= f_\ip^{n}R_{i}$
and $Q_{i}= f_\ip^{n}P_{i}$
for         $1\leq i\leq N(n)$, and note that
the rectangles $Q_{i}$ are
$\i$-spanning $(1,n)$-rectangles
of some   Markov rectangle
$M_{i}$. We note that if $t_{\theta}$ is not the identity
there might be a non-empty set $V_n$ of values of $i$
such that   $S_i$ is not be contained in the Markov rectangle $M_i$.
However, since there are a finite number of Markov rectangles,
the  cardinality of the set $V_n$ is bounded away from infinity,
independently of $n \ge 0$. Hence, we desregard in what follows
these values of $i \in V_n$, since the measure of the corresponding sets $S_i$
converges to $0$ when $n$ tends to infinity.
To prove the theorems we firstly
note that by Lemma
\ref{rem:close_ratios} and
by \eqref{eqn:rhos} we obtain that, if
$Q_{i}$, $P_{i}$ and $M_{i}$ are as
above, for all $\xi,\eta\in M_{i}$,
$$
\rho_{\xi}(S_{i}:Q_{i})\in\left(
1\pm\cO(\epsilon^{n})\right)
\rho_{\eta}(S_{i}:Q_{i}),$$
for some $0<\epsilon < 1$.
Thus, since $\mu
(S_{i})=\mu(\pi_\i(S_{i}))$ and
$\mu (Q_{i})=\mu(\pi_\i(Q_{i}))$ and by \eqref{eqn:ratio_decomp3},
if $\xi\in t_{M_{i}}^\ip$, $$\mu
(S_{i})\in\left(
1\pm\cO(\epsilon^{n})\right)
\rho_{\xi}(S_{i}:Q_{i})
\mu(Q_{i}).
$$
Now consider the case of Theorem
\ref{thm:ratiodecomp_general}.
Then, since $R_{i}$ and $P_{i}$ are
contained in the same Markov
rectangle,
$\rho_{\xi}(S_{i}:Q_{i})$ equals
$\rho_{\xi_{i}}(R:M)$ for some
$\xi_{i}\in t_{R_{i}}^\ip$ and
$\rho(Q_{i})=\rho(P_{i})$ which is equal to
$\mu_\ip  P_{i}  $
since $P_{i}$ is an
$\i$-spanning rectangle of the Markov rectangle $M$.
Thus we
have deduced that up to addition of
a term that is $\cO(\nu^{n})$,
$$\mu(R)\in\left(
1\pm\cO(\epsilon^{n})\right)
\sum_{i=1}^{N(n)} \rho_{\xi_{i}}
(R:M)
\mu_\ip ( P_i ).$$
Equation (\ref{eqn:gen_thm})
follows on taking the limit
$n\to\infty$.

Now   consider the case of Theorem
\ref{thm:ratiodecomp_SRB}.  Under
its hypotheses we have that
$\rho_{\xi}(S_{i}:Q_{i}) =
r^\i(\xi \cap S_{i}:\xi\cap Q_{i})$
by   \eqref{ewqw} and  \eqref{eqn:rhos}.
By the $f$-invariance of $r^\i$
there is $\xi_{i}\in t_{R_{i}}^\ip$
such that
$r^\i(\xi_{i}:t_{\theta}(\xi_{i}))
= r^\i(\xi \cap S_{i}: \xi \cap Q_{i})$.
Thus, as above, we deduce that
$$
\mu(R)\in\left(
1\pm\cO(\epsilon^{n})\right)
\sum_{i=1}^{N(n)}
r^\i(\xi_{i}:t_{\theta}(\xi_{i}))
\mu_\ip(t_\theta R_{i}) \ .
$$
Equation (\ref{eqn:gen_SRB})
follows on taking the limit
$n\to\infty$.
\qed

\begin{lemma}\label{lemma:ratio}
Let $\cS$ be a $C^{1+}$ self-renormalisable
structure on $\TT^\i$ with natural measure  $\mu_\i=\mu_{\cS,\delta}$
for some $\delta >0$.
Suppose that $R$ is contained in a
$(n_{s},n_{u})$-rectangle and
that $R'$ and $R''$ are $\ip$-spanning
rectangles contained in $R$.
Suppose in addition that either
(i) $R$ is contained in a  Markov
rectangle or (ii) $\TT^\i$ does not have gaps and
there is a holonomy injection
of $R$ into a Markov rectangle (in this case $\delta=1$ and $P=0$).
Then for every $\i$-leaf segment $\xi \in t^{\ip}_{R}$
we have that
\begin{equation}
\label{hjyjyyyhrr}
\frac{\mu(R')}{\mu(R'')}\in\left(
1\pm\cO(\epsilon^{n_{s}+n_{u}})\right)
\rho_{\xi}(R':R'')
\end{equation}
for some constant $0<\epsilon<1$
independent of $R$, $R'$, $R''$, $n_{s}$
and $n_{u}$, (and in case (ii)
$\rho_{\xi}(R:R')= r^\i(\xi \cap R:\xi \cap R')$).
\end{lemma}

\begin{proof}  We give the proof for the second
case since that for the first is
similar.
By Theorem \ref{thm:ratiodecomp_SRB},
we have that
$$
\frac{\mu (R)}{\mu
(R')}=\frac{\int_{t_{R}^{\ip}}
r^\i (\xi:t_{\theta}(\xi))
\mu_{R,M}^{\theta} (d\xi
)}
{\int_{t_{R'}^{\ip}} r^\i
(\xi:t_{\theta}(\xi))
\mu_{R',M}^{\theta} (d\xi )}
=
\frac{\int_{t_{R}^{\ip}}
r_{\xi}(R:R') r^\i
(\xi:t_{\theta}(\xi))
\mu_{R,M}^{\theta} (d\xi
)}
{\int_{t_{R'}^{\ip}} r^\i
(\xi:t_{\theta}(\xi))
\mu_{R',M}^{\theta} (d\xi )} .
$$
where
$r_{\xi}^\i (R:R') =r^\i (R \cap \xi:R' \cap \xi)$.
Let $F=f_\ip^{n_s+n_u}$.
By inequality \eqref{1eq:Holderness_of_ratios}
(or  inequality  \eqref{fdfdgghbnn} in case (i)),
there is $0<\epsilon<1$ such that
$$
r_{F\eta}^\i(F R: F R')\in
(1\pm\cO(\epsilon^{n_{s}+n_{u}}))
r_{F \xi}^\i(F R:F R')
$$
for all $\xi,\eta\in t_{R}^{\ip}$.
Thus,
$$
r_{\eta}^\i(R:R')\in
(1\pm\cO(\epsilon^{n_{s}+n_{u}}))r_{\xi}^\i(R:R')
$$
and so
$$
\label{rggggga1}
\frac{\mu (R)}{\mu (R')}\in
(1\pm\cO(\epsilon^{n_{s}+n_{u}}))
r_{\xi}^\i(R:R') \ .
$$
Similarly,
$$
\label{rggggga3}
\frac{\mu (R) }{ \mu (R'')} \in
(1\pm\cO(\epsilon^{n_{s}+n_{u}}))
r_{\xi}^\i(R:R') \ .
$$
Putting  together the previous two equations
we obtain \eqref{hjyjyyyhrr}.
\end{proof}

\begin{lemma}
\label{dsgdgnnww}
Let $\cS$ be a $C^{1+}$ self-renormalisable
structure on $\TT^\i$ with natural measure  $\mu_\i=\mu_{\cS,\delta}$
for some $\delta >0$.
The values   $$
\rat_{\cS,\delta}(\xi \cap R':\xi \cap R'')=\rho_{\xi}(R':R'')
$$
(as in Lemma \ref{lemma:ratio})
determine   an $\i$-measure ratio function  $\rat_{\cS,\delta}$
with the following properties:
(i) The Gibbs measure $\nu_\rho$ determined by the $\i$-measure ratio function $\rat_{\cS,\delta}$
(see Remark \ref{prator})
is the same as the Gibbs measure
$\nu_{\cS,\delta}$ which is $C^{1+}$ realisable by the self-renormalisable structure  $\cS$;
(ii)
If  $\TT^\i$ is a no-gap train-track then $\rat_{\cS,1} = r$,
where $r$  is the ratio function determined by the
$C^{1+}$ self-renormalisable
structure $\cS$.
\end{lemma}

\begin{proof}
Let us prove this lemma first in the case where
the train-track $\TT^\i$ does not have gaps and then
in the case where
the train-track $\TT^\i$ has gaps.

\noindent
\emph{(i)  $\TT^\i$ does not have gaps}.  Then $\delta=1$ and, by
Lemma \ref{rem:close_ratios} (iii),
we have  $\rho_{\xi}(R':R'')= r (\xi \cap R':\xi \cap R'')$
where $r$ is the ratio function determined by the
$C^{1+}$ self-renormalisable
structure $\cS$.
Hence $\rat_{\cS,\delta}=r$ is
an $\i$-measure ratio function. Using \eqref{hjyjyyyhrr},
we get that the Gibbs measure,
which is a $C^{1+}$ realisation of the natural geometric measure  $\mu_{\cS,\delta}$, determines
an   $\i$-measure solenoid function which induces the
$\i$-measure ratio function $\rat_{\cS,\delta}$.

\noindent
\emph{(ii) $\TT^\i$   has gaps}.
Let $\xi,\eta \in \TT^\ip$ be boundaries of Markov rectangles
such that $\xi \cap \eta \ne \emptyset$.
Let $M$, $M'$, $R$ and $R'$ be rectangles such that
$\xi \cap R=\eta \cap M$ and $\xi \cap R'=\eta \cap M'$.
Using \eqref{fdgeghnn} and \eqref{eqn:rhos}, we get that
$$
\rho_{\xi}(R:R')=\rho_{\eta}(M: M') \ .
$$
Therefore, the $\i$-measure ratio function $\rat_{\cS,\delta}$
is well-defined for every $\i$-leaf
segments $I$ and $J$ contained in
a spanning  leaf segment of a Markov rectangle.
By $f$-invariance of $\mu$ and \eqref{hjyjyyyhrr},
we get that
\begin{equation}
\label{gththhh}
\rho_{\xi}(R:R')=\rho_{f_{\ip}\xi}(f_{\ip}R:f_{\ip}R')
\end{equation}
is invariant under $f$. Let  $I$ and $J$ be   $\i$-leaf
segments such that (a) there is
an $\i$-leaf segment $K$ such that  $I,J\subset K$,
and (b) $I$ or $J$ has non-empty interior.
Then there is  $n >0$, $\xi \in \TT^\ip$, $R$ and  $R'$
such that
$f^n_\ip I=\xi \cap R$ and $f^n_\ip J=\xi \cap R'$.
Hence, using \eqref{gththhh}, the ratio
$$
\rat_{\cS,\delta}(I:J)=\rho_{\xi}(R:R')
$$
is well defined independently of $n$.
Using \eqref{hjyjyyyhrr},
we get that \eqref{3eq:Holderness_of_ratios}
is satisfied and
the Gibbs measure,
which is a $C^{1+}$ realisation of
the natural geometric measure  $\mu_{\cS,\delta}$, determines
an   $\i$-measure solenoid function which induces the
$\i$-measure  ratio function $\rat_{\cS,\delta}$.
\end{proof}

\subsection{Dual measure ratio function}

We will show that an $\i$-measure ratio function
$\rho_\i$ determines an unique
dual function $\rho_\ip$
which is an $\ip$-measure ratio function.

\begin{definition}
We say that  the
$\i$-measure
ratio function $\rho_\i$
and the
$\ip$-measure
ratio function  $\rho_\ip$
are \emph{dual} if both
determine the same Gibbs measure
$\nu=\nu_{\rho_\i}=\nu_{\rho_\ip}$ on $\Theta$
(see Remark \ref{prator}).
\end{definition}

\begin{lemma}
\label{ploky}
Let $\cS$ be a $C^{1+}$ self-renormalisable
structure on $\TT^\i$ with $\i$-measure
ratio function $\rho_\i=\rho_{\cS,\delta}$
corresponding to the Gibbs measure $\nu=\nu_{\cS,\delta}$.
Then there is an unique
$\ip$-measure ratio function
$\rho_\ip$
determining the same Gibbs measure
$\nu_{\rho_\ip}=\nu$ on $\Theta$.
\end{lemma}

\begin{figure}[tbp]
\centerline{\includegraphics[width=9cm]{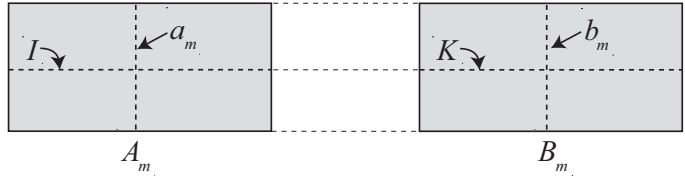} }
%\ \picill 170pt by 117pt (fig2_basicholos.eps) }
\caption{The rectangles $A_{m}=[I,a_m]$ and
$B_{m}=[K,b_m]$.}
\label{AB_m}
\end{figure}

\begin{proof}
Let $\mu=i_*\nu$.
The dual  $\rho_\ip$ of  $\rho_\i$ is constructed as follows:
Let  $I$ and $K$ be
(i) two $\ip$-leaf segments
contained in a common $n$-cylinder $\ip$-leaf,
or also (ii)
two $\ip$-leaf segments contained in a union of two
$n$-cylinders   with a
common endpoint in the case of a local manifold structure.
Choose
$p\in I$ and $p'\in K$. Let
$a_{m}$ be
the $\i$-leaf
$N$-cylinders containing $p$, and
$b_m$
the $\i$-leaf containing
$p'$ and holonomic to $a_{m}$.
Let
$A_{m}=[I,a_m]$ and
$B_{m}=[K,b_m]$ (see Figure \ref{AB_m}).
Now, let us prove that
\begin{rlist}
\item
\begin{equation}
\label{eqn:to_prove}
\frac{\mu(A_{m+1})}{\mu(B_{m+1})} \in
(1\pm\cO(\epsilon^{n+m}))
\frac{\mu(A_{m})}{\mu(B_{m})}
\end{equation}
for some $0<\epsilon < 1$;
\item the dual measure raio function is given by
\begin{equation}
\label{fdfdg}
\rho_\ip (I:K)= \lim_{m \to \infty}
\frac{\mu(A_{m})}{\mu(B_{m})} \  ;
\end{equation}
\end{rlist}
By Lemma \ref{lemma:ratio}, there is $0<\epsilon < 1$
such that
$$ \mu(A_{m+1})/\mu(A_{m})
\in (1\pm
\cO(\epsilon^{n+m}))\rho_\i (a_{m+1}:a_{m}) \ ,$$
and, similarly,
$$\mu(B_{m+1}) / \mu(B_{m})
\in (1\pm
\cO(\epsilon^{n+m}))\rho_\i (b_{m+1}:b_{m})   \ .
$$
Since $\rho_\i$ is an $\i$-ratio function,
$$
\rho_\i (a_{m+1}:a_{m}) \in
\cO(\epsilon^{n+m}))\rho_\i (b_{m+1}:b_{m}) \ .
$$
Therefore, \eqref{eqn:to_prove} follows.
Furthermore,  \eqref{eqn:to_prove} implies  \eqref{fdfdg}.

Using \eqref{eqn:to_prove},
we obtain that   $\rho_\ip$ is an $\ip$-measure ratio function:
$\rho_\ip$ is $f$-invariant,
$\rho_\ip(I:K) = \rho_\ip(K:I)^{-1}$ and
$$\rho_\ip(I:K) = \rho_\ip(I_1:K)+\rho_\ip(I_2:K)$$ for   $\i$-leaf segments
$I_1$ and $I_2$ with at most one common point and
such that $I=I_1 \cup I_2$.
Again using  \eqref{eqn:to_prove},
$\rho_\ip$ satisfies inequality \eqref{3eq:Holderness_of_ratios}.
\end{proof}

\section{Cocycle-gap pairs}
\label{dfgggvvd}

In this section,
we are going to introduce
the gap-cocyle pairs $(\gamma,J)$
consisting of a  gap ratio function $\gamma$
and of a measure-length ratio cocycle $J$ satisfying
a gap-cocyle property.
In \S ~\ref{fgrrffghtrd}, we will  apply this description to give
an explicit geometric construction of all $C^{1+}$
self-renormalisable structures and all $C^{1+}$
hyperbolic diffeomorphisms
which have  natural geometric measures.

\subsection{Measure-length
ratio cocycle}\label{sect:measure-length}

Let $\TT^\i$ be a gap train-track.
For each Markov rectangle $R$ let
$t_R^\ip$ be the set of
$\i$-segments of $R$. Let us denote by $\TT^\ip_o$
the   disjoint union $\sqcup_{i=1}^m t_{R_i}^\ip$ over
all Markov rectangles $R_1,\ldots,R_m$
(without doing any extra-identification).
In this section,
for every $\xi \in \TT^\ip_o$ and $n\ge 1$, we denote by $\xi_n$
the $n$-cylinder $\pi_\i f^{n-1}_\ip \xi$ of $\TT^\i$.

\begin{definition}
\label{cocycledddd}
Let $\TT^\i$ be a gap train-track and
$\rho$   be a $\i$-measure ratio function.
We say that
$J :\TT^\ip_o \to \reals^+$
is a $(\rho,\delta,P)$ \emph{$\i$-measure-length ratio
cocycle}
if $J =\kappa/(\kappa\comp f_\ip)$
where $\kappa$ is a
positive H\"older continuous
function on $\TT^\ip_o$ and is
bounded away from $0$, and
\begin{equation}
\label{aazzqq3}
\sum_{f_\ip  \eta = \xi }
J (\eta) \rho(f_\ip \eta:m(f_\ip \eta))^{1/\delta}  e^{P/\d} < 1 \ ,
\end{equation}
for every $\eta \in \TT^\ip_o$.
\end{definition}

We note that in \eqref{aazzqq3},
the  mother of $\eta$ is not defined
because $\eta$ is a leaf primary cylinder segment,
and so we used instead the mother of the
leaf $2$-cylinder $f_\ip \eta$.

Let us consider a $C^{1+}$
self-renormalisable structure $\cS$
on $\TT^{\i}$, and fix a bounded
atlas $\cB$ for $\cS$. Let $\delta > 0$.
By Theorem \ref{thm:SRB},
the $C^{1+}$
self-renormalisable structure $\cS$ $C^{1+}$-realises a Gibbs measure
$\nu=\nu_{\cS,\delta}$
as a natural invariant measure
$\mu=\mu_{\cS,\delta}= i_* \nu$ with pressure $P=P(\cS,\delta)$.
Let $\rho=\rho_{\cS,\delta}$ be the corresponding
$\i$-measure ratio function (see Lemma \ref{dsgdgnnww}).
Since $\mu$ is a natural geometric measure,
for every $\xi \in \TT^\ip_o$,
the ratios
$|\xi_n|_i e^{-nP /\d}/\mu(\xi_n)^{1/\d}$
are
uniformly bounded away from $0$ and
$\infty$, where the length $|\xi_n|_i$ is measured
in any chart  $i \in \cB$ containing $\xi_n$ in its domain.
Therefore, $$
\kappa_i (\xi_n) = \frac{ |\xi_n|_i e^{-nP / \d}}{\mu(\xi_n)^{1/\d}}
$$
is well-defined.
By Lemma \ref{lemma:ratio},  we get
$$
\frac{\mu_\i( \xi_{n})}{\mu_\i(m \xi_{n})}
\in
(1\pm\cO(\epsilon^{n})) \rho(f_\ip \xi:m(f_\ip \xi))
$$
for some $0 <\epsilon < 1$. Hence, the ratios
$\mu_\i(\xi_n)/\mu_\i (m \xi_n)$
converge exponentially fast along
backward orbits $\xi$ of
cylinders.
By \eqref{eqn:12133},
we get that
$|\xi_n|/|m \xi_n|$ also converge exponentially fast along
backward orbits $\xi$ of
cylinders.
Therefore,  it follows that there is
a H\"older function $J_{\cS,\delta}:\TT^\ip_o \to \reals$
such that
\begin{equation}
\label{ggsgghbvb}
\frac{\kappa_i (\xi_n)}{\kappa_i (m\xi_n)}
\in        (1\pm\cO(\epsilon^{n})) J_{\cS,\delta} (\xi)
\end{equation}
for some $0 <\epsilon < 1$.

\begin{lemma}
\label{gfthttt}
Let $\TT^\i$ be a gap train-track.
Let $\cS$ be a $C^{1+}$
self-renormalisable structure,
and $\delta > 0$. Let  $\mu_{\cS,\delta}$ be the natural geometric measure
with pressure $P=P(\cS,\delta)$, and
$\rho=\rho_{\cS,\delta}$  the corresponding
$\i$-measure ratio function.
The function $J_{\cS,\delta}:\TT^\ip_o \to \reals^+$
given by \eqref{ggsgghbvb}
is a $(\rho, \delta, P)$
$\i$-measure-length ratio
cocycle.
\end{lemma}

\begin{proof}
If $I$ is an $n$-cylinder in $\TT^\i$,
then $\sum_{mI'=I} |I'| < |I|$, where  the lengths
are measured in the same chart. Thus,
since
$|I'| =\kappa_i (I') \mu_\i (I')^{1/\delta}
e^{(n+1)P/\d}$ we deduce that
$$
\sum_{mI'=I} \frac {\kappa_i (I')}
{\kappa_i (I)}
~
\left(\frac{\mu_\i(I')}{\mu_\i(I)} \right)^{1/\d} ~
e^{P/\d} < 1.
$$
For every  $\xi \in \TT^\ip_o$, we have that
$\tau_\ip \eta=\xi$ if, and only if,   $\eta_{n+1} \subset \xi_n$
for every $n \ge 1$.
Hence, the H\"older continuous function $J=J_{\cS,\delta}$  satisfies
\eqref{aazzqq3}.

Now, suppose that        $\xi  \in \TT^\ip_o$
is such that there exists  $p \ge 1$
with  the   property  that $\xi_{n_p}\subset \xi_{n}$
for every $n \ge 1$.
By \eqref{ggsgghbvb}, we get
\begin{eqnarray*}
\frac{\kappa_{i_0} (\xi_{jp})}{\kappa_{i_0}
(\xi_{(j-1)p})} & = &
\prod_{l=0}^{p-1}
     \frac{\kappa_{i_{l+1}}
(\xi_{(j-1)p+l+1})}
{\kappa_{i_l}
(\xi_{(j-1)p+l})} \\
& \in &
(1\pm\cO(\nu^{(j-1)p}) )
\prod_{l=0}^{p-1}
J (f_\ip^l(\xi))
\end{eqnarray*}
where $i_0,\ldots ,i_{p-1}$ are charts contained in a bounded atlas
of $\cS$.
Thus, for all  $1 <m < M$, we have
\begin{eqnarray*}
\frac{\kappa_{i_0} (\xi_{Mp})}{\kappa_{i_0}
(\xi_{mp})} & = &\prod_{n=0}^{M-m-1}
\frac{\kappa_{i_0}
(\xi_{(n+m+1)p})}{\kappa_{i_0}
(\xi_{(n+m)p})} \\
& \in &
(1\pm\cO(\nu^{mp})) \left[
\prod_{l=0}^{p-1}
J (f_\ip^l(\xi))
\right]^{M-m}.
\end{eqnarray*}
Since the term on the left of
this equation is uniformly
bounded away from $0$ and
$\infty$, it follows that $
\prod_{l=0}^{p-1}
J (f_\ip^l(\xi)) =1$.
From Liv\v sic's
theorem (e.g.\ see
\cite{KatokHassBook}) we get that
$J_{\cS,\delta}=\kappa/(\kappa\comp f_\ip)$
where $\kappa$ is a
positive H\"older continuous
function on $\TT^\ip_o$ and is
bounded away from $0$.
\end{proof}

\subsection{Gap ratio function}

\label{gapratioeee}

Let $\TT^\i$ be a gap train-track.
Let $\cG^\ip$ be the set of all pairs  $(\xi_1:\xi_2) \in \TT^\ip_o \times \TT^\ip_o$
such that $m f_\ip \xi_1 = m f_\ip \xi_2$.
The metric $d_\L$ induces a natural metric
$d_{\cG^\ip}$ on $\cG^\ip$ given by
$$
d_{\cG^\ip}((\xi_{1}:\xi_{2}),(\eta_{1}:\eta_{2})) =
\max  \{d_\L(\xi_{1},\eta_{1}) ,d_\L(\xi_{2},\eta_{2}) \} \ .
$$

\begin{definition}
\label{gapratiosds}
A
function $\gamma:\cG^\ip \to \reals^+$ is an   \emph{$\i$-gap ratio
function} if it satisfies
the following conditions:
\begin{rlist}
     \item $\gamma
     (\xi_{1}:\xi_{2})$ is
     uniformly bounded away from
     $0$ and $\infty$;
     \item $\gamma (\xi_{1}:\xi_{2}) =
     \gamma
     (\xi_{1}:\xi_{3})\gamma
     (\xi_{3}:\xi_{2})$;
     \item there are $0 < \theta  < 1$ and $C > 1$ such that
\begin{equation}\label{eqn:crf_Holder}
|\gamma (\xi_{1}:\xi_{2})  - \gamma(\eta_{1}:\eta_{2})| \le
C \left ( d_{\cG^\ip} ((\xi_{1}:\xi_{2}) , \gamma(\eta_{1}:\eta_{2})) \right)^\theta \ .
     \end{equation}
\end{rlist}
\end{definition}
We note that part (ii) of this
definition implies that $\gamma
(\xi_{1}:\xi_{2})=\gamma
(\xi_{2}:\xi_{1})^{-1}$.

Let $\cS$ be
a $C^{1+}$ self-renormalisable
structure on $\TT^{\i}$ and
$\cB$ a bounded atlas for
$\cS$.  Then  the
gap ratio function $\gamma_{\cS}$
is well-defined by
\begin{equation}
\label{ffbgnnnn}
\gamma_{\cS}  (\xi:\eta) =
\lim_{n \to \infty}  \frac{|\pi_{\TT^\i} f_\i^n\xi|_{i_n}}
{|\pi_{\TT^\i} f_\i^n\eta|_{i_n}}
\end{equation}
where $i_n \in \cB$  contains in its domain
the $n$-cylinder $m f_\i^n\xi$ (we note that $m f_\i^n\xi=m f_\i^n\eta)$.

\subsection{Ratio functions}

We are going to construct the ratio function  of a
$C^{1+}$ self-renormalisable structure from the
gap ratio function  and measure-length ratio cocycle.

\begin{lemma}
\label{dfgdgdeewkuy}
Let $\TT^\i$ be a gap train-track.
Let $\cS$ be a $C^{1+}$ self-renormalisable structure.
Let $r_{\cS}$ be the corresponding $\i$-ratio function.
Let $\delta >0$ and  let
$\mu=\mu_{\cS,\delta}$ be the natural geometric measure
with pressure
$P=P(\cS,\delta)$ and $\rho_{\cS,\delta}$ the corresponding $\i$-measure
ratio function. Let $J_{\cS,\delta}$ and $\gamma_\cS$ be
the corresponding $\i$-gap ratio function and $\i$-measure-length ratio cocycle.
Then the following   equalities are satisfied:
\begin{rlist}
\item
Let $I$ be an $\i$-leaf $n$-cylinder
contained in the $\i$-leaf $(n-1)$-cylinder
$L$.  Then
\begin{equation}
\label{fdfdbng}
r_{\cS}(I:L) = J_{\cS,\delta} (\xi_I)\,
\rho_{\cS,\delta}(I:L)^{1/\d}\, e^{P/\d}
\end{equation}
where $\xi_I=   f_\i^{n-1} I \in \TT^\i_o$.
\item
Let $I$  be  an
$n$-cylinder and  $K$   an $n$-gap
and   both contained in a
$(n-1)$-cylinder $L$. Then
\begin{equation}
\label{fdfdbng3}
r_{\cS}(I:K) = r_{\cS}(I:L)
\frac{\sum_{G \subset L} \gamma_{\cS} (G:K)}
{1-\sum_{D \subset L} r_{\cS}        (D:L)}.
\end{equation}
where the   sum in the numerator  is over all  $n$-gaps $G \subset L$
and the sum in the denominator is over all  $n$-cylinders $D \subset L$.
\end{rlist}
\end{lemma}

\begin{proof}
For every $n$-cylinder $I \subset \TT^\i$, define
$\kappa_i (I) =
|I|_i e^{-nP/\d} / \mu_\i(I)^{1/\delta}$ and
let $J_{\cS,\delta}$ be the associate
measure-length cocycle.
Let  $I$ be an  $\i$-leaf $n$-cylinder,
$L$ the $\i$-leaf $(n-1)$-cylinder containing $I$.
Choose
$p\in I$   and let
$U_{m}$ be
the $\ip$-leaf
$m$-cylinders containing $p$.
Let
$A_{m}$ be the rectangle $[I,U_m]$ and
$B_{m}$ be the rectangle $[L,U_m]$.
Then
$f^{m-1}_\ip A_m$ and
$f^{n-1}_\ip B_m$ are $\ip$-spanning
rectangles of some Markov rectangle.
Let $a_{m}$ and $b_{m}$ be the
projections of these into $\TT^\i$.
Then by the invariance of $\mu$,
$\mu (A_m)/\mu (B_m) = \mu_\i
(a_m)/\mu_\i (b_m)$ and therefore
\begin{eqnarray*}
\rho_{\cS,\delta} (I:L)^{1/\d} & = &
\lim_{m\to\infty}\frac{\mu
(A_m)^{1/\d}}{\mu
(B_m)^{1/\d}}  \\
& = &
\lim_{m\to\infty}\frac{\mu_\i
(a_m)^{1/\d}}{\mu_\i
(b_m)^{1/\d}} \\
& = &
\lim_{m\to\infty}
\frac{\kappa(a_{m})^{-1}|a_{m}|_{i_m}e^{-(n+m)P/\d}}
{\kappa(b_{m})^{-1}|b_{m}|_{i_m}e^{-(n+m-1)P/\d}}
\\
& = &   J_{\cS,\delta} (\xi_{I})^{-1} r_{\cS}(I:L) e^{-P/\d}
\end{eqnarray*}
where $|a_m|_{i_m}$ and $|b_m|_{i_m}$ are
measured in a chart $i_m$ of the bounded atlas on
$\TT^\i$, and
$\xi_{I}$ is the leaf primary cylinder segment $f_\i^{n-1}(I)$.
Thus, equation (\ref{fdfdbng}) is satisfied.

We note that the ratio of
the size of $K$ to the size $\l$ of the
totality of gaps $G$ in $L$ is
given by $(\sum_{G \subset L} \gamma_{\cS}
(G:K))^{-1}$ where $\gamma_{\cS}$ is the
gap ratio function and the sum is
over all $n$-gaps in $L$.  But
since the complement of the gaps in
$L$ is the union of $n$-cylinders
we have that the ratio of $\l$ to
the size of $L$ is $1- \sum_{D\subset L} r_{\cS}(D:L)$ where the sum
is over all $n$-cylinders $D$ in
$L$.  Thus we deduce that for
$r_{\cS}(I:K)$ we should take
\begin{equation}
r_{\cS}(I:K) =
r_{\cS}(I:L) \frac{\sum_{G \subset L} \gamma_{\cS} (G:K)}
{1-\sum_{D \subset L} r_{\cS}        (D:L)}
\end{equation}
which proves \eqref{fdfdbng3}.
\end{proof}

\subsection{Cocycle-gap pairs}

In this section, we are going to construct a  \emph{gap-cocyle} map
$b$ which reflects the cylinder-gap
condition  of an $\i$-solenoid
function (see \S ~\ref{dfsfbbvv}), i.e the ratios are well-defined along
the $\i$-boundaries of the Markov rectangles.
Hence, $r$ is an $\i$-ratio function.

Let $\TT^\i$ be a gap train-track and
$\TT^\ip$ a no-gap train-track (as in the case of codimension one attractors
or repellors).
Let $\cQ$ be the set of all periodic orbits $O$ which are
contained in the $\i$-boundaries of the Markov rectangles.
For every periodic orbit $O \in \cQ$,
let us choose a point $x = x(O)$ belonging to the orbit $O$.
Let us denote by $p(x)$ the smallest period of $x$.
Let us denote by $M(1,x)$ and $M(2,x)$
the Markov rectangles containing the point $x$.
Let us denote by $l_i(x)$ the $\i$-leaf $i$-cylinder segment of
Markov rectangle $M(1,x)$
containing the point $x$.
Let $A(f_\i^i(x))$ be  the smallest $\i$-leaf segment containing
all the $\i$-boundary leaf segments
of Markov rectangles
intersecting the global leaf segment passing through the point $f_\i^i(x)$.
Let $q(x)$ be the smallest integer which is a multiple of $p(x)$, such that
$$A(f_\i^i(x)) \subset f_\i^{q(x)+i}(l_i(x))$$
for every $0 \le i < p(x)$.
Let us denote the $\i$-leaf segments $f_\i^{q(x)+i}(l_i(x))$  by
$L_i(x)$.
We note that when using the notation $L_i(x)$,
we will always consider $i$ to be $i$ mod $p(x)$.
For every  $j\in \{1,2\}$, let $J(j,x)$ be the primary $\ip$-leaf segment contained in  $M(j,x)$
with $x$ as an endpoint
such that $R(j,i,x)=
[f_\i^{q(x)+i}(l_i(x)),f_\i^{q(x)+i}(J(j,x))]$ is a rectangle
for every $0 \le i < p(x)$.
Let $\Co(j,i,x) \subset \TT^\ip_o$
be the set of all $\i$-primary leaves $\xi$
of Markov rectangles $M$ such that   $f_\ip \xi \subset L_i(x)$
and  $f_\ip M \cap R(j,i,x)$ has non-empty interior.
Let $\Gap (j,i,x) \subset \cG^\ip$ be the set of all
sister pairs $(\xi_1,\xi_2)$ such that
$m f_\ip \xi_1(=m f_\ip \xi_2)$ is an $\i$-primary
leaf of a Markov rectangle $M$ with the property that
$M \cap R(j,i,x)$ has non-empty interior.
Let $\Co_j=\cup_{O \in \cQ} \cup_{i=0}^{p(x(O))-1}  \Co(j,i,x(O))$
and $\Gap_j=\cup_{O \in \cQ} \cup_{i=0}^{p(x(O))-1}  \Gap (j,i,x(O))$.
Let $\rho$ be an $\i$-measure ratio function with corresponding Gibbs measure
$\nu$.
Let  $\cD_j(\rho, \delta, P)$ be the set of all pairs $(\gamma_j,J_j)$ with the following
properties:
\begin{rlist}
\item
$\gamma_j:\Gap_j \to \reals^+$ is a map;
\item
$J_j:\Co_j \to \reals^+$ is a map satisfying property \eqref{aazzqq3},
with respect to $(\rho, \delta, P)$,
for every $\xi \in \Co_j$ such that
$\xi \subset \cup_{O \in \cQ} \cup_{i=0}^{p(x(O))-1} L_i(x(O))$.
\item
For every $x(O) \in \cQ$, letting $x=x(O)$,
$\prod_{l=0}^{p(x)-1} J_j(f_\ip^l I^j(i,x)=1$,
where $I^j(i,x) \subset \Co_j$ is a $\ip$-primary leaf segment
containing the periodic point $f_\ip^i (x)$.
\end{rlist}
For every $x(O) \in \cQ$, let $x=x(O)$ and let $A(1,x)$ and  $A(2,x)$ be  the $p(x)$-cylinders
of $M(1,x)$ and $M(2,x)$, respectively,  containing the point $x$.
The  points
$$\pi_{\TT^\ip_o} A(j,x), \pi_{\TT^\ip_o} f_\i A(j,x), \ldots, \pi_{\TT^\ip_o} f^{p(x)-1}_\i A(j,x)$$
in $\TT^\ip_o$ form a periodic orbit,
under $f_\i$, with period $p(x)$,
where $\pi_{\TT^\ip_o}: \L \to \TT^\ip_o$
is the natural projection.
The  primary cylinders contained in the sets $\Co(j,i,x)$
are    pre-orbits of the  points $\pi_{\TT^\ip_o} f^{i}_\i A(j,x)$ in $\TT^\ip_o$, under $f_\i$.
Hence, we note that, if   $\prod_{i=0}^{p(x)-1} J(\pi_{\TT^\ip_o} f^{i}_\i A(j,x))=1$, then, by Liv\v sic's
theorem (e.g.\ see
\cite{KatokHassBook}),
there is a map $k$ such that, for every $\xi \in \Co(j,i,x)$,
$J(\xi)=k(\xi)/(k \circ f_\ip)(\xi)$.

We say that $C$ is an \emph{out-gap segment of a rectangle} $R$
if $C$ is a gap segment of $R$ and is not a leaf $n$-gap segment
of any Markov rectangle $M$ such that $M \cap R$ is a rectangle with
non-empty interior.

We say that $C$ is  a \emph{leaf $n$-cylinder segment of a rectangle} $R$,
if $C$ is a leaf $n$-cylinder segment of  a Markov rectangle $M$ such that $M \cap R$ is a rectangle with
non-empty interior.
We say that $C$ is  a \emph{leaf $n$-gap segment of a rectangle} $R$,
if $C$ is a leaf $n$-gap segment of  a Markov rectangle $M$ such that $M \cap R$ is a rectangle with
non-empty interior.
We say that $C$ is  an \emph{$n$-leaf segment of a rectangle} $R$,
if $C$ is a leaf $n$-cylinder segment of $R$ or if $C$ is a leaf $n$-gap segment of $R$.

\begin{lemma}
\label{ggdgdhhbna0}
Let $(\gamma_1,J_1) \in \cD_1(\rho, \delta, P)$.
Let $x=x(O)$, where $O \in \cQ$.
For every $i \in \{0,1, \ldots, p(x)-1\}$
and
for all    $2$-leaf segments $C \subset L_i(x)$ of $R(1,i,x)$,
the  ratios
$r(C:mC)$
are uniquely determined such that
they are invariant under $f$, satisfy the matching condition,
and satisfy equalities \eqref{fdfdbng}  and \eqref{fdfdbng3}.
\end{lemma}

\begin{proof}
If $C \subset L_i(x)$  is a leaf $2$-cylinder   segment of $R(1,i,x)$,
then
we define the ratio $r(C:mC)$, using  \eqref{fdfdbng}, by
\begin{equation}
\label{pofdopgp1}
r(C:mC) = J (\xi_C)\,
\rho(C:mC)^{1/\d}\, e^{P/\d}
\end{equation}
where  $\xi_C=f_\i C    \in \Co_1$.
For every sister pair $(\xi_1:\xi_2) \in \Gap_1$
we define the ratio $r(f_\ip \xi_1:f_\ip \xi_2)$ equal to
$\gamma (\xi_1:\xi_2)$.
If $C \subset L_i(x)$  is a leaf $2$-gap   segment of $R(1,i,x)$,  then
we define the ratio $r(C:mC)$ by
\begin{equation}
\label{pofdopgp3}
r(C:mC)  = \frac
{1-\sum_{D \subset mC} r (D:mC)}
{\sum_{G \subset mC} r (G:C)} \ ,
\end{equation}
where the   sum, in the numerator,  is over all $2$-cylinders $D \subset mC$
of $R(1,i,x)$,
and the sum, in the denominator, is over all $2$-gaps $G \subset mC$ of $R(1,i,x)$.
Hence, $$\sum_{C \subset mC}r(C:mC)=1$$
where the sum is over all $2$-leaf segments
$C \subset mC$ of $R(1,i,x)$.
\end{proof}

\begin{lemma}
\label{ggdgdhhbna1}
Let $(\gamma_1,J_1) \in \cD_1(\rho, \delta, P)$.
Let $x=x(O)$, where $O \in \cQ$, and let $i \in \{0,1, \ldots, p(x)-1\}$.
For all $n \ge 0$, and for all out-gaps and all
$2$-leaf segments $C \subset f^{n+i}_\i \ell_i(x)$ of $f^{n+i}_\i M(1,x)$,
the  ratios $r(C:f^{n+i}_\i \ell_i(x))$
are uniquely determined such that
they are invariant under $f$, satisfy the matching condition,
and satisfy equalities \eqref{fdfdbng}  and \eqref{fdfdbng3}.
\end{lemma}

\begin{proof}
Let us denote $f_\i^n M(1,x)$ by $M_n$ and
$f^n_\i \ell_i(x)$ by $L^n_i$.
The proof follows by induction on  $n \ge 0$.
For  the case $n=0$, the ratios $r(C:L^{n+i}_i)$
are uniquely determined by Lemma
\ref{ggdgdhhbna0}.
Let us prove that the
ratios $r(C:L^{n+1+i}_i)$  are uniquely determined
using  the   induction hypotheses with respect to $n$.
For every  out-gap  and every primary cylinder segment
$C \subset L^{n+1+i}_i$ of $f^{n+i}_\i M(1,x)$,
$f_\ip C$ is a  out-gap  or a $2$-leaf segment.
Hence, by   the   induction hypotheses,
the ratio $r(f_\ip C:L^{n+i}_i)$
is well-defined. Therefore, using the invariance of $f$,
we define
\begin{equation}
\label{pofdopgp5}
r(C:L^{n+1+i}_i) =r(f_\ip C:L^{n+i}_i) \ .
\end{equation}
For every  $2$-leaf segment
$C \subset L^{n+1+i}_i$ of $f^{n+i}_\i M(1,x)$,
the ratio  $r  (C:m C)$
is well-defined by Lemma \ref{ggdgdhhbna0}.
Hence, by  \eqref{pofdopgp5}, we define
$$
\label{pofdopgp7}
r(C:L^{n+1+i}_i)=r ( C:mC) r(mC:L^{n+1+i}_i)
$$
which ends the proof of the induction.
\end{proof}

\begin{lemma}
\label{ggdgdhhbna2}
Let $(\gamma_1,J_1) \in \cD_1(\rho, \delta, P)$.
Let $x=x(O)$, where $O \in \cQ$, and let $i \in \{0,1, \ldots, p(x)-1\}$.
Let $n \ge 0$ and $j \in \{0, \ldots, n\}$.
For all out-gaps and all
$j+2$-leaf segments $C \subset f^{n}_\ip L_i(x)$ of $f^{n}_\ip R(1,i,x)$,
the  ratios $r(C:f^{n}_\ip L_i(x))$
are uniquely determined such that
they are invariant under $f$, satisfy the matching condition,
and satisfy equalities \eqref{fdfdbng}  and \eqref{fdfdbng3}.
\end{lemma}

\begin{proof}
The proof follows by induction in $n \ge 0$.
For  the case $n=0$, noting that $L_i(x) = f^{q(x)+i}_\i \ell_i(x)$,
the ratios $r(C:L_i(x))$
are well-defined  by   Lemma
\ref{ggdgdhhbna1}.
Hence,
using the matching condition, the ratio
$r(f^{n+1}_\ip L_{i+1}(x):f^{n}_\ip L_i(x))$
is well-defined.
Let us prove that for all out-gaps and $j+2$-leaf segments
$C \subset f^{n+1}_\ip L_i(x)$ of $f^{n+1}_\ip R(1,i,x)$, with $1 \le j \le n+1$,
the  ratios $r(C:f^{n+1}_\ip L_i(x))$
are uniquely determined
using  the   induction hypotheses with respect to $n$.
By the induction hypotheses  and by the matching condition,
the ratio $r (f_\i C:f^{n}_\ip L_i(x))$
is well-defined.
By invariance of $f$, we define
$r(C:f^{n+1}_\ip L_i(x)) = r (f_\i C:f^{n}_\ip L_i(x))$.
which ends the proof of the induction.
\end{proof}

Let us attribute the ratios for the cylinders and gaps of $R(2,i,x)$
such that they agree with the ratios previously defined  in $R(1,i,x)$.

\begin{lemma}
\label{ggdgdhhbna3}
Let $(\gamma_1,J_1) \in \cD_1(\rho, \delta, P)$.
Let $x=x(O)$, where $O \in \cQ$, and let $i \in \{0,1, \ldots, p(x)-1\}$.
Let $n \ge 0$ and $j \in \{1, \ldots, n\}$.
For all out-gaps and all
$j+2$-leaf segments $C \subset f^{n}_\ip L_i(x)  \setminus f^{n+1}_\ip L_{i+1}(x)$ of $f^{n}_\ip R(2,i,x)$,
the  ratios $r(C:f^{n}_\ip L_i(x))$
are uniquely determined such that
they are invariant under $f$, satisfy the matching condition,
satisfy equalities \eqref{fdfdbng}  and \eqref{fdfdbng3},
and are well-defined along
the $\i$-boundaries of the Markov rectangles.
Hence, $r$ is an $\i$-ratio function.
\end{lemma}

\begin{proof}
The proof follows by induction in $n \ge 0$.
Let us prove the case $n=0$.
By construction, $L_i(x) \supset A(f_\i^i(x))$, i.e  $L_i(x)$ contains  all the $\i$-boundary leaf segments
of Markov rectangles
intersecting the global leaf segment passing through the point $f_\i^i(x)$.
Hence, if  $G_2 \subset  L_i(x)  \setminus f_\ip L_{i+1}(x)$
is  an out-gap   of $R(2,i,x)$,
then
there is an out-gap or a leaf $2$-gap segment $G_1$ of
$R(1,i,x)$  such that $G_1=G_2$.
Therefore, we define $r(G_2: L_i(x))=r(G_1: L_i(x))$.
Since $L_i(x) \supset A(f_\i^i(x))$, if   $G_2 \subset  L_i(x)  \setminus f_\ip L_{i+1}(x)$
is   a  leaf $2$-gap segment of $R(2,i,x)$
then
there is an out-gap or a leaf $2$-gap segment $G_1$ of
$R(1,i,x)$  such that $G_1=G_2$.
Hence, we define $r(G_2: L_i(x))=r(G_1: L_i(x))$.
If $C_2 \subset  L_i(x)  \setminus f_\ip L_{i+1}(x)$
is   a leaf $2$-cylinder   segment of $R(2,i,x)$,
then
there is  a primary leaf segment or a leaf $2$-cylinder   segment $C_1$ of $R(1,i,x)$
such that $C_2=C_1$.
Therefore, we define $r(C_2: L_i(x))=r(C_1: L_i(x))$.
Let us prove that for all out-gaps and $j+2$-leaf segments
$C \subset f^{n+1}_\ip L_i(x) \setminus f^{n+2}_\ip L_{i+1}(x)$
of $f^{n+1}_\ip R(2,i,x)$, with $1 \le j \le n+1$,
the  ratios $r(C:f^{n+1}_\ip L_i(x))$
are uniquely determined
using  the   induction hypotheses with respect to $n$.
By the induction hypotheses  and by the matching condition,
the ratio $r (f_\i C:f^{n}_\ip L_i(x))$
is well-defined.
By invariance of $f$, we define
$r(C:f^{n+1}_\ip L_i(x)) = r (f_\i C:f^{n}_\ip L_i(x))$.
which ends the proof of the induction.
\end{proof}

\begin{lemma}
\label{ggdgdhhbna4}
Let $(\gamma_1,J_1) \in \cD_1(\rho, \delta, P)$.
Let $x=x(O)$, where $O \in \cQ$, and let $i \in \{0,1, \ldots, p(x)-1\}$.
For all out-gaps and all
$2$-leaf segments $C \subset L_i(x)$ of $R(2,i,x)$,
the  ratios $r(C:L_i(x))$
are uniquely determined such that
they are invariant under $f$, satisfy the matching condition,
satisfy equalities \eqref{fdfdbng}  and \eqref{fdfdbng3},
and are well-defined along
the $\i$-boundaries of the Markov rectangles.
Hence, $r$ is an $\i$-ratio function.
\end{lemma}

\begin{proof}
By construction of $L_i(x)  \setminus f_\ip L_{i+1}(x)$, there is $k = k(n,i,x)$ such that
$L_i(x)  \setminus f_\ip L_{i+1}(x) = \cup_{l=1}^{k} D_l$,
where $D_l$ are out-gaps,
primary leaf segments and $2$-leaf segments of $R(1,i,x)$.
Therefore, $f^{n}_\ip L_i(x)  \setminus f^{n+1}_\ip L_{i+1}(x) = \cup_{l=1}^{k} f^{n}_\ip D_l$
where $f^{n}_\ip D_l$ are out-gaps and
$j+2$-leaf segments of $R(1,i,f^{n}_\ip (x))$ with $0 \le j \le n$.
Hence, by Lemma \ref{ggdgdhhbna3} and
using the matching condition, the ratio
$r(f^{n}_\ip L_i(x)  \setminus f^{n+1}_\ip L_{i+1}(x):f^{n}_\ip L_i(x))$
is well-defined.
Hence, using the matching condition, we define
$$r(f^{n+1}_\ip L_{i+1}(x):f^{n}_\ip L_i(x)) = 1-
r(f^{n}_\ip L_i(x)  \setminus f^{n+1}_\ip L_{i+1}(x):f^{n}_\ip L_i(x)) \ .
$$
Therefore, using again the matching condition, we define
\begin{equation}
\label{wwwqqqp}
r(f^{n+1}_\ip L_{i+n+1}(x):L_i(x))=
\prod_{j=0}^{n} r(f^{j+1}_\ip L_{i+j+1}(x):f^{j}_\ip L_{i+j}(x)) \ .
\end{equation}
Let $M(i,x)$ be the $2$-cylinder of $R(2,i,x)$ containing the point $x$.
Take $N>0$, large enough, such that $f^{N+1}_\ip L_{i+N+1}(x) \subset M(i,x)$.
Hence,
there is $m = m(N,i,x)$ such that
$M(i,x) =  \left( \cup_{l=0}^{m} D_l \right) \cup f^{N+1}_\ip L_{i+N+1}(x)$
where $D_l$ are out-gaps or
$j+2$-leaf segments of $R(1,i,x)$ for some $0 \le j \le N$.
Hence, by Lemma \ref{ggdgdhhbna3}, \eqref{wwwqqqp}  and
using the matching condition, the ratio is well-defined by
$$
r(M(i,x): L_i(x))=
\sum_{l=0}^{m} r( D_l: L_i(x)) + r(f^{N+1}_\ip L_{i+N+1}(x): L_i(x)) \ .
$$
If $C \subset  L_i(x)\setminus M(i,x)$ is a out-gap or a
$2$-leaf segment of $R(2,i,x)$, then, by Lemma \ref{ggdgdhhbna3},
the ratio is well-defined by $r(C: L_i(x))$.
By construction of the ratios, in Lemmas \ref{ggdgdhhbna0}-\ref{ggdgdhhbna4},
they are compatible with the cylinder-gap condition.
\end{proof}

\begin{definition}
\label{aswqss}
Let $(\gamma_1,J_1) \in \cD_1(\rho, \delta, P)$.
Let $x=x(O)$, where $O \in \cQ$, and let $i \in \{0,1, \ldots, p(x)-1\}$.
Let  the  ratios $r(C:L_i(x))$ for all out-gaps and all
$2$-leaf segments $C \subset L_i(x)$ of $R(2,i,x)$ be as
given in Lemma \ref{ggdgdhhbna4}.
For all $\xi \in \Co(2,i,x)$, letting $I=f_\ip \xi \subset L_i(x)$, we define
$$
J_2 (\xi)  =r(I:L_i(x))    r(L_i(x):m I)  \rho(I:K_i)^{-1/\d}  e^{-P/\d} \ .
$$
For all $(C,D) \in \Gap (2,i,x)$, we define
$$ \gamma(C:D)=r(f_\ip C:L_i(x)) r(L_i(x):f_\ip D) \ .$$
\end{definition}

\begin{lemma}
\label{ggdgdhhbna5}
Let $\cD_1(\rho, \delta, P) \ne \emptyset$.
The \emph{gap cocycle map}
$b = b_{\rho, \delta, P}:\cD_1(\rho, \delta, P) \to \cD_2(\rho, \delta, P)$
is well-defined by $b(\gamma_1,J_1) = (\gamma_2,J_2)$
where $\gamma_2$ and $J_2$ are as given in Definition \ref{aswqss}.
Furthermore,
the cocycle gap map $b$ is a bijection.
\end{lemma}

\begin{proof}
Let us check that $(\gamma_2,J_2)$ satisfies properties (i)-(iii)
of $\cD_2(\rho, \delta, P)$.
By construction of the
ratios $r$, in Lemmas \ref{ggdgdhhbna0}-\ref{ggdgdhhbna4},
$(\gamma_2,J_2)$ satisfies properties (i) and (ii) in the definition of
$\cD_2(\rho, \delta, P)$.
Let us check property (iii).
Let us denote by $A$ and  $B$   the $p(x)$-cylinders
of $M(1,x)$ and $M(2,x)$, respectively,  containing the point $x$.
By invariance of $r$, we have that
$r(A:B)=r(f_\i^{p(x)} A:f_\i^{p(x)} B)$, and so
$r(A:f_\i^{p(x)} A)=r(B:f_\i^{p(x)} B)$.
By invariance of the $\i$-measure ratio function $\rho$, we have that
$\rho(A:B)=\rho(f_\i^{p(x)} A:f_\i^{p(x)} B)$, and so
$\rho(A:f_\i^{p(x)} A)=\rho(B:f_\i^{p(x)} B)$.
Since, by hypotheses $\prod_{l=0}^{p(x)-1} J(m^i f^{p(x)-i}_\i A)=1$, we get, from \eqref{fdfdbng},
that
$r(A:f_\i^{p(x)} A)= \rho(A:f_\i^{p(x)} A) e^{p(x)P/\delta}$.
Therefore,
\begin{eqnarray*}
r(B:f_\i^{p(x)} B) & = & r(A:f_\i^{p(x)} A)\\
& = & \rho(A:f_\i^{p(x)} A) e^{p(x)P/\delta}\\
& = & \rho(B:f_\i^{p(x)} B) e^{p(x)P/\delta}
\end{eqnarray*}
and so, using \eqref{fdfdbng}, we obtain that  $\prod_{i=0}^{p(x)-1} J(\pi_{\TT^\ip_o} f^{i}_\i B)=1$.
\end{proof}

\begin{definition}
\label{ghhhhgdd}
Let $\TT^\i$ be a gap train-track.
Let $\delta > 0$ and $P \in\reals$.
Let $\rho$ be an $\i$-measure ratio function
and $\nu=\nu_\rho$
the corresponding Gibbs measure on $\Theta$.
We say that a pair $(\gamma,J)$
is a $(\nu,\delta,P)$
\emph{$\i$ cocycle-gap pair}
if $(\gamma,J)$ has the following properties:
\begin{rlist}
\item
$\gamma$ is an $\i$-gap ratio function.
\item
$J$ is an  $\i$ measure-length ratio cocyle.
\item
If $\TT^\ip$ is a no-gap train-track
then $(\gamma,J)$ satisfies the following
\emph{gap-cocyle property}:
$b(\gamma|\Gap_1,J|\Co_1)=
(\gamma|\Gap_2,J|\Co_2)$, where
$b=b_{\nu,\delta,P} $ is the gap-cocyle map.
\end{rlist}
Let ${\mathcal JG}^\i(\nu,\delta,P)$
be the set of all
$(\nu,\delta,P)$ $\i$ cocycle-gap pairs.
\end{definition}

\begin{lemma}
\label{ggdgdhhbna}
Let $\TT^\i$ be a gap train-track. Let $\delta >0$ and $P \in \reals$.
Let $\rho$ be an $\i$-measure ratio function with corresponding Gibbs measure
$\nu$.
\begin{rlist}
\item
If there is a  $(\rho,\delta,P)$  $\i$-measure-length ratio cocycle, then
the set ${\mathcal JG}^\i(\nu,\delta,P)$ is
an infinite dimensional space.
\item
If $\cS$ is a $C^{1+}$ self-renormalisable structure
with natural geometric measure $\mu_{\cS,\delta}=i_*\nu$
and pressure $P$, then   $(\gamma_\cS,J_{\cS,\delta}) \in {\mathcal JG}^\i(\nu,\delta,P)$.
\item
If the set ${\mathcal JG}^\i(\nu,\delta,P) \neq \emptyset$,
then there is a well-defined injective map $(\gamma,J) \to r(\gamma,J)$
which associates to each cocycle-gap pair $(\gamma,J) \in {\mathcal JG}^\i(\nu,\delta,P)$
an $\i$-ratio function
$r(\gamma,J)$ satisfying  \eqref{fdfdbng} and \eqref{fdfdbng3}.
\end{rlist}
\end{lemma}

\begin{remark}
\label{dsfdgdhghv} Let $0<\delta < 1$ and $P=0$.
Let $\rho$ be an $\i$-measure ratio function with corresponding Gibbs measure
$\nu$. Since $J = 1$ is a $(\rho,\delta,P)$  $\i$-measure-length ratio cocycle,
then, by Lemma \ref{ggdgdhhbna}, the set
${\mathcal JG}^\i(\nu,\delta,P)$ is an infinite dimensional space.
\end{remark}

\noindent
\emph{Proof of Lemma  \ref{ggdgdhhbna}}.
{\it Proof of (i).}
Choose a map  $\gamma_1:\Gap_1 \to \reals^+$.
Let $J_0$ be a $(\rho,\delta,P)$  $\i$-measure-length ratio cocycle,
and let $J_1 = J_0|\Co_1$.
Since $(\gamma_1,J_1) \in \cD_1(\rho, \delta, P)$, by Lemma \ref{ggdgdhhbna5},
the pair $(\gamma_2,J_2) = b_{\rho, \delta, P}(\gamma_1,J_1) \in \cD_2(\rho, \delta, P)$ is well-defined.
Let $k_0$ and $k_2$ be maps such that
$J_0 = k_0/(k_0 \circ f_\ip)$ and   $J_2 = k_2/(k_2 \circ f_\ip)$.
For every $x(O) \in \cQ$, let $x=x(O)$, and let $B$  be the $p(x)$-cylinder
of $M(2,x)$  containing the point $x$.
Recall that the  points
$\pi_{\TT^\ip_o} B$, $\pi_{\TT^\ip_o} f_\i B$, $\ldots$, $\pi_{\TT^\ip_o} f^{p(x)-1}_\i B$
in $\TT^\ip_o$ form a periodic orbit
under $f_\i$, with period $p(x)$, and that the  primary cylinders contained in the set $\Co(2,i,x)$
are    pre-orbits of the  points $\pi_{\TT^\ip_o} f^{i}_\i B$ in $\TT^\ip_o$, under $f_\i$.
Therefore, there is a small neighbourhood $V$ of $\Co_2$ in $\TT^\ip_o$,
there is  $\epsilon >0$, small enough, and
there is an H\"older continuous map $k: \TT^\ip_o \to \reals^+$
with the following properties:
\begin{rlist}
\item
$k|\Co_2 = k_2$,
$k|(\TT^\ip_o \setminus V) = k_0$ and $\Co_1 \subset \TT^\ip_o \setminus V$.
\item
Let
$a=\min_{\xi \in \Co_2} \{J_0 (\xi), J_2 (\xi)\}$
and $b=\max_{\xi \in \Co_2} \{J_0 (\xi), J_2 (\xi)\}$,
and let $J = k/(k \circ f_\ip)$.
For every $\xi \in V$,
we have that $a- \epsilon \le J (\xi) \le b+\epsilon$,
and, so, $J$ satisfies the cocycle-gap property.
\end{rlist}

Choosing an H\"older continuous map $\gamma:\cG^\ip \to \reals^+$ such that
$\gamma|\Gap_1 = \gamma_1$ and $\gamma|\Gap_2 = \gamma_2$ and
by property (i) above, the pair $(\gamma, J)$ satisfies \eqref{aazzqq3}.
Therefore, the pair $(\gamma, J)$ is contained in ${\mathcal JG}^\i(\nu,\delta,P)$.
Using that \eqref{aazzqq3} is an open condition,
the above construction allows us to construct
an infinite set of $\i$-measure-length ratio cocycles
and  an infinite set of gap ratio functions
such that  the corresponding pairs are contained in ${\mathcal JG}^\i(\nu,\delta,P)$.

\noindent
{\it Proof of (ii).}
Let  $\cS$ be a $C^{1+}$ self-renormalisable structure
with natural geometric measure $\mu_{\cS,\delta}=i_*\nu$
and pressure $P$.
By Lemma \ref{gfthttt}, $J_{\cS,\delta}$ is a $(\rho,\delta,P)$ $\i$-measure-length ratio cocycle
and, by \eqref{ffbgnnnn}, $\gamma_\cS$ is an $\i$-gap ratio function.
If $\TT^\ip$ is a no-gap train-track,
using \eqref{fdfdbng} and \eqref{fdfdbng3},
the pair
$(\gamma_\cS,J_{\cS,\delta})$
satisfies the cocycle-gap condition
because
the ratio function $r_\cS$
associated to $\cS$ is well-defined along
the $\i$-boundaries of the Markov rectangles.

\noindent
{\it Proof of (iii).}
The equations  \eqref{fdfdbng} and \eqref{fdfdbng3}
give us an inductive construction, on the level $n$
of the $n$-cylinders and $n$-gaps, of a ratio function $r$ in terms
of $(\rho,J,\gamma,\delta,P)$ with the property that the ratio between
a leaf $n$-cylinder segment $C$ and a leaf $n$-cylinder or $n$-gap segment $D$
with a common endpoint with $C$ is bounded away from zero and infinity independent of $n$ and of
the cylinders and gaps considered.
The construction gives that $r$ is invariant under $f$.
The H\"older continuity of $\gamma$, $J$ and $\rho$ implies that $r$ satisfies
\eqref{1eq:Holderness_of_ratios}.
If $\TT^\ip$ is a no-gap train-track, by the construction of the cocycle-gap condition,
the ratio function $r$  is well-defined along
the $\i$-boundaries of the Markov rectangles.
Hence, $r$ is an $\i$-ratio function.
\qed

\section{Realisations of Gibbs measures}
\label{fgrrffghtrd}

In this section, we are going to give
an explicit geometric construction of all $C^{1+}$
hyperbolic diffeomorphisms
which have a natural geometric measure,
and we will prove the theorems stated in the introduction of the paper.

\subsection{One-dimensional realisations}

Let $\cS$ be a $C^{1+}$ self-renormalisable structure
on a  train-track $\TT^\i$. In Lemma \ref{dsgdgnnww}
we have shown
that the map
\begin{equation}
\label{afro}
(\cS, \delta) \to \rat_{\cS,\delta}
\end{equation}
is well-defined where $ \rat_{\cS,\delta}$ is the
$\i$-measure ratio function associated to a Gibbs measure
$\nu_{\cS,\delta}=\nu$
such that $\mu_{\cS,\delta} = (i_\i)_*\nu_\i$ is a natural
geometric measure of $\cS$.

\begin{lemma}
\label{trdfg} (Rigidity)
Let $\TT^\i$ be a no  gap train-track (and $\delta=1$).
The map $\cS  \to \rat_{\cS,\delta}$
is an one-to-one correspondence between $C^{1+}$ self-renormalisable structures
on $\TT^\i$ and $\i$-measure ratio functions.
Furthermore, $\rat_{\cS,\delta}=r_\cS$ where
$r_\cS$ is the  ratio function determined by $\cS$.
\end{lemma}

However, if $\TT^\i$ is a   gap train-track
then the set of pre-images of the map
$(\cS, \delta) \to \rat_{\cS,\delta}$
is an infinite dimensional space (see Lemma \ref{gdvaswwa} below).

\begin{proof}
By Theorem \ref{thm:SRB}, the
$C^{1+}$ self-renormalisable structure
$\cS$ realises a Gibbs measure $\nu = \nu_{\cS,\delta}$.
By Lemma \ref{dsgdgnnww}, we get that $\rat_{\cS,\delta}=r_\cS$.
Since, by Lemma \ref{fghfhttf}, the ratio function $r_\cS$ determines
uniquely the $C^{1+}$ self-renormalisable structure $\cS$,
the map $\cS \to \rat_{\cS,\delta}$
is an one-to-one correspondence.
\end{proof}

\begin{definition}
\label{rdfrgfff}
Let $\TT_\i$ and $\TT_\ip$ be  (gap or  no-gap) train-tracks.
Let $\rho$ be an $\i$-measure
ratio function and $\nu=\nu_\rho$ on $\Theta$
the corresponding Gibbs measure (see Remark \ref{prator}).
Let us denote by $\cD^\i(\nu,\delta,P)$
the set of all
$C^{1+}$ self-renormalisable structures $\cS$
with geometric natural measure $\mu_{\cS,\delta}=(i_\i)*\nu_\i$
and pressure $P$.
\end{definition}

By Lemma \ref{trdfg}, if  $\TT^\i$ is a no-gap train-track,
and   $\delta =1$ and  $P=0$, the set
$\cD^\i(\nu,\delta,P)$  is a singleton.

Let $\TT^\i$ be a gap train-track and
$\cS$  a $C^{1+}$ self-renormalisable structure in $\cD^\i(\nu,\delta,P)$.
In Theorem \ref{gapratioeee},
we associate to the $C^{1+}$ self-renormalisable structure
$\cS$ a measure-length ratio cocycle $J_\cS$,
and,
in \S ~\ref{gfthttt}, we associate to the
$C^{1+}$ self-renormalisable structure
$\cS$
a gap ratio function $\gamma_\cS$.
By Lemma \ref{ggdgdhhbna},
if $\TT^\ip$ is a no-gap train-track
then    the cylinder-gap condition of $r_\cS$ implies that
the pair $(\gamma_\cS,J_{\cS,\delta})$
satisfies the cocycle-gap condition.
Therefore,
the map
\begin{equation}
\label{zzzsss}
\cS \to (\gamma_\cS,J_{\cS,\delta})
\end{equation}
between  $C^{1+}$ self-renormalisable structures contained in $\cD^\i(\nu,\delta,P)$
and pairs contained in ${\mathcal JG}^\i(\nu,\delta,P)$
is well-defined.

\begin{definition}
\label{dfgdrfgggx}
The
\emph{$(\delta_\i,P_\i)$-bounded solenoid equivalence class}  of
a Gibbs measure $\nu$
is the set of all solenoid functions $\s_\i$
with  the following properties:
There is $C=C(\s_\i)>0$ such that for every pair
$(\xi,D) \in \msc_\i$
$$
\left|\delta_\i \log s_\i (  D_\L \cap   \xi_\L:   \xi_\L)
-\log  \rho_{\xi} (D) -nP_\i  \right|
< C   \  ,
$$
where  (i) $\rho$ is the $\i$-extended measure scaling function of $\nu$,
(ii) $s_\i$ is the scaling function determined by $\s_\i$,
(iii)  $\xi_\L=i(\pi_\ip^{-1} \xi)$ is an $\ip$-leaf primary cylinder segment
and (iv) $D_\L=i(\pi_\i^{-1} D)$ and so
$D_\L \cap   \xi_\L$ is an $\i$-leaf   $n$-cylinder segment.
\end{definition}

\begin{remark}
Let  $\s_{1,\i}$ and
$\s_{2,\i}$  be two solenoid functions
in the same $(\delta_\i,P_\i)$-bounded solenoid equivalence  class
of a Gibbs measure $\nu$.
Using the fact that  $\s_{1,\i}$ and
$\s_{2,\i}$ are bounded away from zero, we obtain that
the corresponding scaling functions also satisfy inequality \eqref{dsdfgreww}
for all pairs $(J,m^iJ)$ where $J$ is an
$\i$-leaf $(i+1)$-cylinder.
Hence, the solenoid functions  $\s_{1,\i}$ and
$\s_{2,\i}$ are in the same bounded   equivalence class (see Definiton \ref{fgrdegerdg}).
\end{remark}

By Lemma \ref{gdvaswwa}, below, the set  ${\mathcal JG}^\i(\nu,\delta,P)$
gives a parametrization of all $C^{1+}$ self-renormalisable
structures $\cS$ which are pre-images of the $\i$-measure ratio function
$\rat_{\nu,\i}$, under the map
$\cS \to \rat_{\cS,\delta}$ given in \eqref{afro},
with a natural geometric  measure $\mu_\i = (i_\i)_*\nu_\i$
and pressure $P(\cS,\delta)=P$.
Hence,  ${\mathcal JG}^\i(\nu,\delta,P)$
forms a moduli space for the set of
all $C^{1+}$ self-renormalisable structures
in $\cD^\i(\nu,\delta,P)$.

\begin{lemma} (Flexibility)
\label{gdvaswwa}
Let $\TT_\i$ be a gap
train-track.
Let $\rho$ be an $\i$-measure ratio function
and $\nu=\nu_\rho$
the corresponding Gibbs measure on $\Theta$.
\begin{rlist}
\item
Let $\delta>0$ and $P \in \reals$
be such that ${\mathcal JG}^\i(\nu,\delta,P) \ne \emptyset$.
The map $\cS \to (\gamma_\cS,J_{\cS,\delta})$
determines a one-to-one correspondence
between $C^{1+}$ self-renormalisable structures
in $\cD^\i(\nu,\delta,P)$ and cocycle-gap pairs
in ${\mathcal JG}^\i(\nu,\delta,P)$.
\item A $C^{1+}$ self-renormalisable structure
$\cS$ is contained in $\cD^\i(\nu,\delta,P)$ if,
and only if, the $\i$-solenoid function $\s_{\cS}$ is contained in the
$(\delta,P)$-bounded solenoid equivalence class of $\nu$ (see Definition \ref{dfgdrfgggx}).
\end{rlist}
\end{lemma}

\noindent
\emph{Proof of (i)}.
Let us prove that $(J,\gamma) \in {\mathcal JG}^\i(\nu,\delta,P)$
determines an unique $C^{1+}$ self-renormalisable structure $\cS$
with a natural geometric measure $\mu_{\cS,\delta}=(i_\i)_* \nu_\i$.
By Lemma \ref{ggdgdhhbna},
the pair $(J,\gamma)$ determines an unique $\i$-ratio function
$r=r_\i(J,\gamma)$.
By Lemma \ref{fghfhttf}, the $\i$-ratio function $r$
determines an unique $C^{1+}$ self-renormalisable structure $\cS$
with an atlas $\cB(r)$.
Let us prove that $\mu_\i =(i_\i)_*\nu_\i$ is a natural geometric measure of $\cS$
with the given $\delta$ and $P$. Let $\rho$ be the $\i$-measure ratio function
associated to the Gibbs measure $\nu$.
By Lemma \ref{frgvnhgdfd}, for every leaf $n$-cylinder or $n$-gap segment
$I$ we obtain that
\begin{equation}
\label{sdkbjbjib!1}
\mu_\i (I) = \cO( \rho (I \cap \xi: \xi))
\end{equation}
for every $\xi \in \pi_\ip(I)$.
On the other hand,
by construction of the ratio function $r_\i$ and using \eqref{fdfdbng}, we get
$$
\rho (I \cap \xi: \xi) =
e^{-n P} r(I \cap \xi: \xi) ^{\delta}
\prod_{j=0}^{{n-1}} \left(
J \left( \tau_\ip^j (\xi)\right) \right)^{-\delta}\ .
$$
Since $J$ is a H\"older cocycle,
it follows that
$\prod_{j=0}^{{n-1}}
J \left( \tau_\ip^j (\xi)\right)=k(\xi)/k( \tau_\ip^n (\xi)) $
is uniformly bounded away from zero and infinity,
where $k$ is an H\"older continuous positive function.
By \eqref{eqn:12133}, we get that
\begin{equation}
\label{sdkbjbjib!23}
r(I \cap \xi: \xi) = \cO \left( |I|_j \right)
\end{equation}
where $j \in \cB(r)$
and $I$ is contained in the domain of $j$.
Hence,
\begin{equation}
\label{sdkbjbjib!2}
\rho (I \cap \xi : \xi) = \cO \left( |I|_j^\delta e^{-n P} \right) \ .
\end{equation}
Putting together equations \eqref{sdkbjbjib!1} and
\eqref{sdkbjbjib!2}, we deduce that
$ \mu_\i (I ) = \cO \left( |I|_j^\delta e^{-n P} \right) $, and so
$\mu_\i =(i_\i)^*\nu_\i$ is a natural geometric measure of $\cS$
with the given $\delta$ and $P$.

\noindent
\emph{Proof of (ii)}.
Let $\cS$ be a $C^{1+}$ self-renormalisable structure in $\cD^\i(\nu,\delta,P)$.
Then, putting together \eqref{sdkbjbjib!23} and \eqref{sdkbjbjib!2},
there is $\kappa>0$ such that
$$
|\delta \log r_\i(I \cap \xi : \xi) - \log \rho (I \cap \xi : \xi) -np |
< \kappa
$$
for every leaf $n$-cylinder $I$ and $\xi \in \pi_\ip(I)$.
Thus the solenoid function $r|\sol^\i$ is in the
$(\delta,P)$-bounded solenoid equivalence class of $\nu$.

Conversely, let $\cS$ be a $C^{1+}$ self-renormalisable structure
in the
$(\delta,P)$-bounded solenoid equivalence class of $\nu$
and $\mu_\i =(i_\i)_*\nu_\i$, i.e.
there is $\kappa>0$ such that
\begin{equation}
\label{prati}
|\delta \log r_\i(I \cap \xi : \xi) - \log \rho (I \cap \xi : \xi) -np |
< \kappa
\end{equation}
for every leaf $n$-cylinder $I$ and $\xi \in \pi_\ip(I)$.
Hence, using \eqref{sdkbjbjib!1}
and \eqref{sdkbjbjib!23} in \eqref{prati},
we get
$\mu_\i (I) = \cO \left( |I|_j^\delta e^{-n P} \right)$.
Since $\mu_\i =(i_\i)_*\nu_\i$
we get that $\cS$ is contained in $\cD^\i(\nu,\delta,P)$.
\qed

\subsection{Two-dimensional realisations}

We start by giving the definition of
a natural geometric measure for a
$C^{1+}$ hyperbolic diffeomorphism.

\begin{definition}
\label{fgtyjtyjyu}
For $\i \in \{s,u\}$,
if $\TT^\i$ is a gap train-track assume $0<\delta_\i <1$, and
if $\TT^\i$ is a no-gap train-track take $\delta_\i =1$.
\begin{rlist}
\item
Let $g$ be a
$C^{1+}$ hyperbolic diffeomorphism in $\cT(\fL )$.
We say that $g$ has a \emph{natural geometric measure}
$\mu=\mu_{g,\delta_s,\delta_u}$
with pressures $P_s=P_s(g,\delta_s,\delta_u)$ and $P_u=P_u(g,\delta_s,\delta_u)$ if,
there is $\kappa >1$ such that
for all leaf $n_s$-cylinder $I_s$, for all leaf $n_u$-cylinder $I_u$,
\begin{equation}
\label{fre}
\kappa ^{-1} <
\frac{\mu (R)}
{|I_u|^{\delta_{u}} |I_s|^{\delta_s}
e^{-n_{s}P_s-n_{u}P_u}}
< \kappa \ ,
\end{equation}
where $R$ is the $(n_s,n_u)$-rectangle $[I_s,I_u]$
and
where the lengths $|\cdot |$ are measured in the
stable and unstable $C^{1+}$foliated
lamination atlasses $\cA_s(g,\rho)$ and $\cA_u(g,\rho)$
of $g$ with respect to some Riemannian metric $\rho$.
\item
We say that a $C^{1+}$ hyperbolic diffeomorphism
with a natural geometric measure
$\mu=\mu_{g,\delta_s,\delta_u}$
with pressures $P_s=P_s(g,\delta_s,\delta_u)$ and $P_u=P_u(g,\delta_s,\delta_u)$
is a $C^{1+}$
\emph{realisation of a Gibbs measure} $\nu=\nu_{g,\delta_s,\delta_u}$
if $\mu=i_* \nu$.
We denote by $\cT(\nu,\delta_s,P_s,\delta_u,P_u)$ the set
of all these $C^{1+}$ hyperbolic diffeomorphisms $g \in \cT(\fL )$.
\end{rlist}
\end{definition}

A $C^{1+}$ hyperbolic diffeomorphism $g \in \cT(\fL )$
determines an unique pair $(\cS(g,s),\cS(g,u))$
of $C^{1+}$ stable and unstable
self-renormalisable structures (see Lemma \ref{ghjjjyuu}).
By Theorem \ref{thm:SRB}, for $\delta_s>0$ and $\delta_u>0$,
the pair $(\cS(g,s),\cS(g,u))$ of self-renormalisable structures
determines an unique pair of natural geometric measures
$(\mu_{\cS(g,s),\delta_s},\mu_{\cS(g,u),\delta_u})$
corresponding to a unique  pair of Gibbs measures
$(\nu_{\cS(g,s),\delta_s},\nu_{\cS(g,u),\delta_u})$.
Furthermore, by Lemma \ref{dsgdgnnww}, the self-renormalisable structures $(\cS(g,s),\cS(g,u))$
determine a pair of
measure ratio functions $(\rho_{\cS(g,s),\delta_s},\rho_{\cS(g,u),\delta_u})$
of $(\nu_{\cS(g,s),\delta_s},\nu_{\cS(g,u),\delta_u})$.

\begin{lemma}
\label{aasskk}
Let $g$ be a
$C^{1+}$ hyperbolic diffeomorphism
contained in $\cT(\fL )$.
The following statements are equivalent:
\begin{rlist}
\item
$g$ has a natural geometric measure
$\mu_{g,\delta_s,\delta_u}$;
\item
$g$ is a $C^{1+}$ realisation of a Gibbs measure
$\nu_{g,\delta_s,\delta_u}$;
\item
$\nu_{\cS(g,s),\delta_s} = \nu_{\cS(g,u),\delta_u}$;
\item
The $s$-measure ratio function $\rho_{\cS(g,s),\delta_s}$
is dual to the $u$-measure ratio function $\rho_{\cS(g,u),\delta_u}$.
\end{rlist}
Furthermore, if $g$ has a natural geometric measure
$\mu_{g,\delta_s,\delta_u}$, then
$(\pi_s)_*\mu_{g,\delta_s,\delta_u} = \mu_{\cS(g,s),\delta_s}$
and $(\pi_u)_*\mu_{g,\delta_s,\delta_u} = \mu_{\cS(g,u),\delta_u}$.
\end{lemma}

\begin{proof}
By Lemma \ref{ploky}, (iii) is equivalent to (iv).
By definition
if $g$ is a $C^{1+}$ realisation of a Gibbs measure
$\nu_{g,\delta_s,\delta_u}$ then
$\mu_{g,\delta_s,\delta_u}=i_*\nu_{g,\delta_s,\delta_u}$
is a natural geometric measure of $g$, and so
(ii) implies (i).
Let us prove first that  (i) implies (ii) and (iii), and secondly that
(iii) implies (i). Then the last paragraph of this lemma follows from \eqref{anho}
below which ends the proof.

\noindent
\emph{(i) implies (ii) and (iii)}.
Let $\mu_{g,\delta_s,\delta_u}$ be the natural geometric measure
of $g$.
Since the stable and unstable lamination atlasses $\cA_s(g,\rho)$ and $\cA_u(g,\rho)$ of $g$ are
$C^{1+}$ foliated
(see \S ~\ref{sfgvgbg})
and by construction of the  $C^{1+}$ train-track atlasses $\cB_s(g,\rho)$ and $\cB_u(g,\rho)$, in
\S ~\ref{fddhhhuy}, we obtain
that there is $\kappa_1 \ge 1$ with the property that,  (for $\i = s$ and $u$) and
for every $\i$-leaf $n_\i$-cylinder $I$,
\begin{equation}
\label{freee}
\kappa_1^{-1} |I|_\rho \le |I'|_j \le \kappa_1 |I|_\rho
\end{equation}
where $I' = \pi_{\TT^\i} (I)$, where $|I'|_j$ is measured in any chart
$j \in \cB_\i(g,\rho)$ and where $|I|_\rho$ is the length in the Riemannian
metric $\rho$ of the minimal full
$\i$-leaf containing $I$.
Let $I'_\Lambda$ be the $(1,n_\i)$-rectangle in $\Lambda$
such that $\pi_{\TT^\i} (I'_\Lambda)=I'$.
Noting that
$(\pi_{\TT^\i})_* \mu_{g,\delta_s,\delta_u}(I') = \mu_{g,\delta_s,\delta_u}(I'_\Lambda)$,
by \eqref{fre} and \eqref{freee},
there is $\kappa_2 \ge 1$ such that
\begin{equation}
\label{fredee}
\kappa_2^{-1} \le
\frac{(\pi_{\TT^\i})_* \mu_{g,\delta_s,\delta_u}(I')}
{|I'|_j^{\delta_{\i}}
e^{-n_{\i}P_\i}}
\le \kappa_2 \ ,
\end{equation}
for every $n_\i$-cylinder $I'$ on the train-track.
By Theorem \ref{thm:SRB},  the natural geometric measure
determined by
the  $C^{1+}$
self-renormalisable structure $\cS(g,\i)$ and by $\delta_\i>0$
is uniquely determined by \eqref{fredee}. Hence,
\begin{equation}
\label{anho}
(\pi_{\TT^s})_* \mu_{g,\delta_s,\delta_u} = \mu_{\cS(g,s),\delta_s}
~~~{\rm and}~~~
(\pi_{\TT^u})_* \mu_{g,\delta_s,\delta_u} = \mu_{\cS(g,u),\delta_u} \ .
\end{equation}
Therefore, the Gibbs measures
$\nu_{\cS(g,s),\delta_s}$ and $ \nu_{\cS(g,u),\delta_u}$ on $\Theta$
are equal which proves   (iii), and
$\mu_{g,\delta_s,\delta_u} = i_* \nu_{\cS(g,s),\delta_s} =i_* \nu_{\cS(g,u),\delta_u}$
which proves (ii).

\noindent
\emph{(iii) implies (i)}.
Let us denote $\nu_{\cS(g,s),\delta_s} = \nu_{\cS(g,u),\delta_u}$ by $\nu$.
Let $\mu = i_* \nu$.
For $\i\in\{ s,u\}$,
by definition of a $C^{1+}$ realisation of a Gibbs measure
as a self-renormalisable structure $\cS(g,\delta_\i)$,
for every $\i$-leaf $n_\i$-cylinder $I_\i$,
there is $\kappa_3 \ge 1$ such that
$$
\kappa_3^{-1} \le
\frac{\mu_\i (I'_\i)}
{|I'_\i|_j^{\delta_{\i}} e^{-n_{\i}P_\i}}
\le \kappa_3 \ ,
$$
where $I'_\i = \pi_{\TT^\i} (I)$ and $|I'_\i|_j$ is measured in any chart
$j \in \cB_\i(g,\rho)$.
Hence, by \eqref{freee}, for $\i = s$ and $u$, we obtain that
\begin{equation}
\label{anmkt}
\mu_\i (I'_\i) = \cO \left(|I_\i|_\rho^{\delta_{\i}}
e^{-n_{\i}P_\i} \right)
\ .
\end{equation}
Let $R$ be the rectangle  $[I_s,I_u]$.
By Theorem \ref{thm:ratiodecomp_general},
$$
\mu (R) =
\int_{I_\ip'}
\rat_{\i,\xi}
(R:M) \mu_\ip (d\xi ) \ ,
$$
where $M$ is the Markov rectangle containing  $R$.
By Theorem \ref{thm:SRB} (i) and (ii), we get that
$\rat_{\i,\xi} (R:M) = \cO(\mu_\i(I_\i'))$
for every $\xi \in \pi_{\TT^\ip}(R)$. Hence
\begin{equation}
\label{amjk}
\mu (R) = \cO (\mu_s(I'_s) \mu_u(I'_u)) \ .
\end{equation}
Putting together \eqref{anmkt} and \eqref{amjk},
we get
\begin{equation}
\label{fredd}
\mu (R) = \cO  \left(|I_u|_\rho^{\delta_{u}}
|I_s|_\rho^{\delta_s} e^{-n_{s}P_s-n_{u}P_u} \right)
\end{equation}
and so $\mu$ is a natural geometric measure.
\end{proof}

\begin{lemma}
\label{fghbhhjbsde}
The map $g \to (\cS(g,s),\cS(g,u))$ gives a one-to-one
correspondence
between $C^{1+}$ conjugacy classes of hyperbolic diffeomorphisms
contained in $\cT(\nu,\delta_s,P_s,\delta_u,P_u)$
and pairs of $C^{1+}$ self-renormalisable structures contained in
$\cD^s (\nu,\delta_s,P_s) \times \cD^u(\nu,\delta_u,P_u)$.
\end{lemma}

\begin{proof}
By Lemma \ref{aasskk}, if $g \in \cT(\nu,\delta_s,P_s,\delta_u,P_u)$
then, for $\i\in\{ s,u \}$,
$\cS(g,\i) \in \cD^\i(\nu,\delta_\i,P_\i)$.
Conversely, by Lemma \ref{ghjjjyuu},
a pair $(\cS_s,\cS_u) \in \cD^s (\nu,\delta_s,P_s) \times \cD^u(\nu,\delta_u,P_u)$ determines
a  $C^{1+}$
hyperbolic diffeomorphism $g$ such that
$\cS(g,s)=\cS_s$ and $\cS(g,u)=\cS_u$ and
$\nu_{\cS(g,s),\delta_s} = \nu_{\cS(g,u),\delta_u} = \nu$.
Therefore, by Lemma \ref{aasskk}, we obtain that $g$ is a $C^{1+}$
realisation of the Gibbs measure $\nu$.
\end{proof}

\begin{lemma} (Dual-rigidity)
\label{dsgdgnnww3}
Let $\TT^\ip$ be a no-gap train-track (and so $\delta_\ip=1$ and $P_\ip=0$).
For every $\delta_\i > 0$ and every
$C^{1+}$ $\i$-self-renormalisable structure $\cS_\i$
there is an unique
$C^{1+}$ $\ip$-self-renormalisable structure $\cS_\ip$
such that the $C^{1+}$ hyperbolic diffeomorphism $g$ corresponding
to the pair $(\cS_s,\cS_u)=(\cS(g,s),\cS(g,u))$ has a natural geometric measure
$\mu_{g,\delta_s,\delta_u}$.
Furthermore, $\mu_{\cS_s,\delta_s} = (\pi_{\TT^s})_*\mu_{g,\delta_s,\delta_u}$
and $\mu_{\cS_u,\delta_u} = (\pi_{\TT^u})_*\mu_{g,\delta_s,\delta_u}$.
\end{lemma}

\begin{proof}
By Theorem \ref{thm:SRB}, a $C^{1+}$ self-renormalisable structure
$\cS_\i$ and $\delta_\i >0$
determine an unique Gibbs measure $\nu=\nu_{\cS_\i,\delta_\i}$ and $P_\i \in \reals$ such that
$\cS_\i \in \cD^\i(\nu,\delta_\i,P_\i)$ is a $C^{1+}$
realisation of $\nu$.
By Lemma \ref{dsgdgnnww},
the $C^{1+}$ self-renormalisable structure
$\cS_\i$ determines   an $\i$-measure ratio function  $\rat_{\cS_\i,\delta_\i}$ for the Gibbs
measure $\nu$.
By Lemma \ref{ploky}, the $\i$-measure ratio function  $\rat_{\cS_\i,\delta_\i}$ determines an
unique
$\ip$-measure ratio function
$\rho_\ip$  of $\nu$ on $\Theta$.
By Lemma \ref{trdfg}, there is  an unique
$C^{1+}$ self-renormalisable structure
$\cS_\ip$, with $\ip$-measure ratio function
$\rho_{\cS_\ip,1}=\rho_\ip$, which
is a $C^{1+}$ realisation of
the Gibbs measure  $\nu$.
By Lemma  \ref{ghjjjyuu}, the pair $(\cS_s,\cS_u)$ determines
a  $C^{1+}$
hyperbolic diffeomorphism $g$ such that
$\cS(g,s)=\cS_s$ and $\cS(g,u)=\cS_u$.
Hence, $\nu_{\cS(g,s),\delta_s}=\nu$
and $\nu_{\cS(g,u),\delta_u}=\nu$
which implies that
$\nu_{\cS(g,s),\delta_s} = \nu_{\cS(g,u),\delta_u}$. Therefore,
by Lemma \ref{aasskk},   $g$ is a $C^{1+}$
realisation of the Gibbs measure $\nu$
with natural geometric measure $\mu_{g,\delta_s,\delta_u}=i_*\nu$.
Thus, $\mu_{\cS_s,\delta_s} = (\pi_{\TT^s})_*\mu_{g,\delta_s,\delta_u}$
and $\mu_{\cS_u,\delta_u} = (\pi_{\TT^u})_*\mu_{g,\delta_s,\delta_u}$.
\end{proof}

Recall the definition of the   maps $g \to (\cS(g,s),\cS(g,u))$
and $\cS(g,\i) \to (\gamma_{\cS(g,\i)},J_{\cS(g,\i),\delta_\i})$
for $\i$ equal to $s$ and $u$.

\begin{theorem} (Flexibility)
\label{fghhhytr}
Let $\TT_\i$  be a gap
train-track. Let $\nu$ be a Gibbs measure determining  an
$\i$-measure ratio function.
Let $\delta_\i>0$ and $P_\i \in \reals$
be such that
${\mathcal JG}^\i (\nu,\delta_\i,P_\i) \ne \emptyset.$
\begin{rlist}
\item
(Smale horseshoes)
Let $\delta_\ip>0$ and $P_\ip \in \reals$
be such that
${\mathcal JG}^\ip (\nu,\delta_\ip,P_\ip) \ne \emptyset.$
The  map
$$g \to (\gamma_{\cS(g,s)},J_{\cS(g,s),\delta_s},\gamma_{\cS(g,u)},J_{\cS(g,u),\delta_u})$$
gives an one-to-one correspondence
between $C^{1+}$ conjugacy classes of hyperbolic diffeomorphisms
in $\cT(\nu,\delta_s,P_s,\delta_u,P_u)$
and pairs of stable and unstable cocycle-gap pairs
in ${\mathcal JG}^s(\nu,\delta_s,P_s) \times {\mathcal JG}^u(\nu,\delta_u,P_u)$.
\item
(Codimension one attractors and repellors)
Let $\delta_\ip=1$ and $P_\ip =0$.
The  map
$g \to (\gamma_{\cS(g,\i)},J_{\cS(g,\i),\delta_\i})$
gives an one-to-one correspondence
between $C^{1+}$ conjugacy classes of hyperbolic diffeomorphisms
in $\cT(\nu,\delta_s,P_s,\delta_u,P_u)$
and pairs of stable and unstable cocycle-gap pairs
in ${\mathcal JG}^\i(\nu,\delta_\i,P_\i)$.
\end{rlist}
\end{theorem}

\begin{proof}
Statement (i) follows from putting together the
results of lemmas  \ref{gdvaswwa} and \ref{fghbhhjbsde}.
Statement (ii) follows as statement (i) using the fact that,
by Lemma \ref{trdfg}, the $C^{1+}$ self-renormalisable structure
$\cS(g,\i)$ uniquely determines $\cS(g,\ip)$ in this case.
\end{proof}

\section{Eigenvalues}

In this section, we show that  the set of
stable and unstable eigenvalues of all  periodic points
of the hyperbolic diffeomorphisms
is a complete invariant of the Lipschitz conjugacy classes
extending  the results of De la Llave, Marco and Moriyon.
Furthermore, we extend the  eigenvalue formula of A. N.  Liv\v sic  and Ja. G. Sinai
from Anosov diffeomorphisms
to $C^{1+}$ hyperbolic diffeomorphisms.

\subsection{Lipschitz conjugacy classes}

\begin{lemma}
\label{akju1}
Let $\TT^\i$ be a (gap or a no-gap) train-track.
Let $\delta>0$ and $P \in \reals$.
Let $\cS_1 \in \cD^\i(\nu_1,\delta,P)$
and $\cS_2 \in \cD^\i(\nu_2,\delta,P)$
be $C^{1+}$ self-renormalisable structures
The following statements are equivalent:
\begin{rlist}
\item
The $C^{1+}$ self-renormalisable structures $\cS_1$ and $\cS_2$ are Lipschitz conjugate;
\item
The Gibbs measures $\nu_1$ and $\nu_2$ are equal;
\item
The solenoid functions $s_{\cS_1}$ and $s_{\cS_2}$
are in the same bounded equivalence class (Definition \ref{fgrdegerdg}).
\end{rlist}
\end{lemma}

\noindent
\emph{Proof  that (i) is equivalent to (ii)}.
Using \eqref{dcsfere}, if $\nu_1 = \nu_2$
then the $C^{1+}$ self-renormalisable structure
$\cS_1$ is Lipschitz conjugate to $\cS_2$.
Conversely,
if $\cS_1$ is Lipschitz conjugate to $\cS_2$
then the $C^{1+}$ self-renormalisable structure
$\cS_1$ (and $\cS_2$) satisfies \eqref{dcsfere} with respect
to the measures
$\mu_{\i,1} = (i_\i)_* \nu_{\i,1}$
and $\mu_{\i,2} = (i_\i)_* \nu_{\i,2}$.
By Theorem \ref{thm:SRB}, there is an unique
$\tau$-invariant Gibbs measure satisfying \eqref{dcsfere}
and so $\nu_1 = \nu_2$.

\noindent
\emph{Proof  that (ii) is equivalent to (iii)}.
Using \eqref{dsdfgreww} and \eqref{eqn:12133}, we obtain that
the  $C^{1+}$  self-renormalisable structures   $\Str$ and $\Str'$ on $\TT^\i$ are
in the same Lipschitz equivalence class if, and only if,
the corresponding solenoid functions $r_{\Str,\i}|\sol^\i$
and $r_{\Str',\i}|\sol^\i$ are in the same bounded equivalence class.
Hence, statement  (ii) is equivalent to statement (iii).
\qed

\begin{lemma}
\label{asdro}
Let $g_1$ and $g_2$ be $C^{1+}$ hyperbolic diffeomorphisms.
The following statements are equivalent:
\begin{rlist}
\item
The diffeomorphism $g_1$ is Lipschitz conjugate to $g_2$.
\item
For $\i$ equal to $s$ and $u$, $\cS(g_1,\i)$ is Lipschitz conjugate to $\cS(g_2,\i)$.
\item
For $\i$ equal to $s$ and $u$, the solenoid functions $s_{g_1,\i}$ and $s_{g_2,\i}$
are in the same bounded equivalence class (Definition \ref{fgrdegerdg}).
\end{rlist}
\end{lemma}

\begin{figure}[tbp]
\centerline{\includegraphics[width=12.5cm]{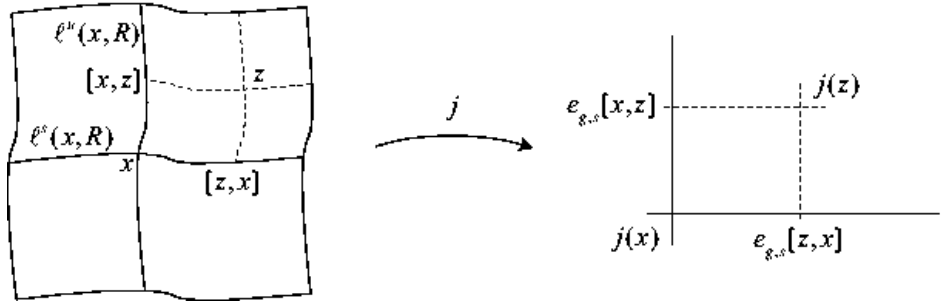} }
%\ \picill 170pt by 117pt (fig2_basicholos.eps) }
\caption{An orthogonal chart.}
\label{orthogonal_chart}
\end{figure}

\noindent
\emph{Proof  that (i) is equivalent to (ii)}.
For all $x\in \Lambda$, let $A$ be a small open set of $M$
containing $x$, and let $R$ be a rectangle
(not necessarily a Markov rectangle)
such that $A \cap \Lambda \subset R$.
We construct an orthogonal chart $j: R \to \reals ^2$
as follows.
Let $e_{g,s}: \ell^s(x,R) \to \reals$
be a chart contained in $\cA ^s (g,\rho)$ and $e_{g,u}:\ell^u(x,R) \to \reals$
be a chart contained in $\cA ^u (g,\rho)$.
The orthogonal chart $j$  on $R$ is now given by
$j (z) = (e_{g,s}[z,x]),e_{g,u}[x,z])) \in \reals ^2$ (see Figure \ref{orthogonal_chart}).
By \cite{HR}, the orthogonal chart $j : R \to \reals ^2$
has an extension ${\hat j} : B \to \reals ^2$ to an open set $B$
of the surface such that
${\hat j}$ is $C^{1+}$ compatible with the charts in the
$C^{1+}$ structure $\cC(g)$ of the surface $M$.
Hence, using the orthogonal charts,
any two  $C^{1+}$ hyperbolic diffeomorphisms $g_1$ and $g_2$
are Lipschitz conjugate if, and only if
the charts in $\cA ^\i (g_1,\rho_1)$ are bi-Lipschitz compatible
with the charts in $\cA ^\i (g_2,\rho_2)$ for $\i$ equal to $s$ and $u$.
By construction
of the train-track atlasses
$\cB ^\i (g_1,\rho_1)$  and $\cB ^\i (g_2,\rho_2)$
from the lamination atlasses $\cA ^\i (g_1,\rho_1)$  and $\cA ^\i (g_2,\rho_2)$,
the charts in $\cA ^\i (g_1,\rho_1)$ are bi-Lipschitz compatible
with the charts in $\cA ^\i (g_2,\rho_2)$ if, and only if,
the charts in $\cB ^\i (g_1,\rho_1)$ are bi-Lipschitz compatible
with the charts in $\cB ^\i (g_2,\rho_2)$.
Hence, the $C^{1+}$ hyperbolic diffeomorphisms $g_1$ and $g_2$
are Lipschitz conjugate if, and only if,
for  $\i$ equal to $s$ and $u$, the corresponding
$C^{1+}$  self-renormalisable structures
$\cS(g_1, \i)$ and $\cS(g_2, \i)$ are Lipschitz conjugate.
Therefore, statement (i) is equivalent to statement (ii).

\noindent
\emph{Proof  that (ii) is equivalent to (iii)}. Follows from Lemma \ref{akju1}.
\qed

\begin{lemma}
\label{fgdgfgdds1}
Let $\delta_s>0$, $\delta_u>0$ and $P_s,P_u \in \reals$.
\begin{rlist}
\item
A $C^{1+}$ hyperbolic diffeomorphism
$g$ is contained in $\cT(\nu,\delta_s,P_s,\delta_u,P_u)$
if, and only if,
for $\i$ equal to $s$ and $u$,
the $\i$-solenoid function $\s_{g,\i}$ is contained in the
$(\delta_\i,P_\i)$-bounded solenoid equivalence class of $\nu$
(see Definition \ref{dfgdrfgggx}).
\item
If $g_1 \in \cT(\nu_1,\delta_s,P_s,\delta_u,P_u)$
and $g_2 \in \cT(\nu_2,\delta_s,P_s,\delta_u,P_u)$
are $C^{1+}$ hyperbolic diffeomorphisms
then $g_1$ is Lipschitz conjugate to $g_2$
if, and only if, $\nu_1 = \nu_2$.
\end{rlist}
\end{lemma}

\noindent
\emph{Proof of (i)}.
By Lemma \ref{ghjjjyuu}, the $C^{1+}$ hyperbolic diffeomorphism
$g$  determines an unique pair $(\cS(g,s),\cS(g,u))$
of $C^{1+}$ self-renormalisable structures such that
$\sigma_{g,s} = \sigma_{\cS (g,s),s}$ and
$\sigma_{g,u} = \sigma_{\cS (g,u),u}$.
By Lemma \ref{fghbhhjbsde},
$g \in \cT(\nu,\delta_s,P_s,\delta_u,P_u)$
if, and only if,
$(\cS(g,s),\cS(g,u)) \in
\cD^s (\nu,\delta_s,P_s) \times \cD^u(\nu,\delta_u,P_u)$.
By Lemma \ref{gdvaswwa} (ii), for $\i$ equal to $s$ and $u$,
$\cS(g,\i) \in \cD^\i (\nu,\delta_\i,P_\i)$
if, and only if,
$\cS(g,\i)$ is contained in the
$(\delta_\i,P_\i)$-bounded solenoid equivalence class of $\nu$
which ends the proof.

\noindent \emph{Proof of (ii)}.
By Lemma \ref{fghbhhjbsde},
$g_1 \in \cT(\nu_1,\delta_s,P_s,\delta_u,P_u)$
and $g_2 \in \cT(\nu_2,\delta_s,P_s,\delta_u,P_u)$
if, and only if, for $\i$ equal to $s$ and $u$,
$
\cS(g_1,\i)  \in
\cD^\i (\nu_1,\delta_\i,P_\i)
$
and
$
\cS(g_2,\i)  \in
\cD^\i (\nu_2,\delta_\i,P_\i)
$.
By Lemma \ref{akju1},
$\cS(g_1,\i)$ and $\cS(g_2,\i)$
are Lipschitz conjugate if, and only if, $\nu_1 = \nu_2$.
Since, by Lemma \ref{asdro},
$g_1$ and $g_2$ are Lipschitz conjugate if, and only if,
for $\i$ equal to $s$ and $u$,
$\cS(g_1,\i)$ and $\cS(g_2,\i)$
are Lipschitz conjugate, we get that $g_1$ and $g_2$ are Lipschitz conjugate if, and only
if, $\nu_1 = \nu_2$.
\qed

\subsection{Extending the result of De la Llave, Marco and Moriyon}
\label{dfegdaddrggggd}

Let $\cP$ be the set of all periodic points
in $\L$ under $f$. Let $p(x)$ be the (smallest) period of
the periodic point $x \in \cP$.
For every $x \in \cP$ and $\i \in \{s,u\}$,
let $j:J \to \reals$ be a chart in $\cA(g,\rho_g)$
such that $x \in J$.
The  \emph{eigenvalue $\lambda_{g,\i}^\i(x)$ of   $x$}
is the derivative of the map
$j^{-1} f^p j$ at $j(x)$.

For $\i \in \{s,u\}$,
by construction of the train-tracks,
$\cP^\i=\pi_{\TT^\i}(\cP)$
is the set of all periodic points
in $\TT^\i$ under the Markov map $f_\i$.
Furthermore, $\pi_{\TT^\i}|\cP$ is an injection
and the periodic points   $x \in \L$
and  $\pi_{\TT^\i}(x) \in \TT^\i$
have the same period $p(x)=p(\pi_{\TT^\i}(x))$.
Let us denote $\pi_{\TT^\i}(x)$ by $x_\i$.
Let $\cS_\i$ be a $C^{1+}$ self-renormalisable structure.
Let  $j:J \to \reals$ be a train-track
chart of $\cS_\i$
such that $x_\i \in J$.
The \emph{eigenvalue $\lambda_{\cS_\i}(x_\i)$ of $x_\i$}
is the derivative of the map
$j \circ \tau_\i^{p(x_\i)} \circ j^{-1}$
at $j(x_\i)$, where $\tau_\i$ is the Markov map
on the train-track $\TT^\i$.

For every $x \in \cP$, every $\i \in \{s,u\}$
and every $n \ge 0$,  let $I^\i_{n}(x)$
be an  $\i$-leaf $(np(x)+1)$-cylinder segment such that
$x \in I^\i_{n}(x)$ and $f^{p(x)}_\i I^\i_{n+1}(x)=I^\i_{n}(x)$.

\begin{lemma}
\label{fgregtgeeew}
For $\i \in \{s,u\}$, let  $\cS_\i \in \cD^\i(\nu,\delta_\i,P_\i)$ be a
$C^{1+}$ $\i$ self-renormalisable structure.
For every $x \in \cP$,
\begin{eqnarray}
\lambda_{\cS_\i}(x_\i)  &= & r_{\cS_\i} (I^\i_{0}(x):I^\i_{1}(x))
\label{grgrgr000} \\
&= &  \rho_{\nu,\i} (I^\i_{0}:I^\i_{1})^{-1/\d_\i} e^{-p(x) P_\i/\d_\i}
\label{grgrgr001} \\
&= &  \rho_{\nu,\ip} (I^\ip_{0}:I^\ip_{1})^{-1/\d_\i} e^{-p(x) P_\i/\d_\i}
\label{grgrgr003} \ ,
\end{eqnarray}
where $r_{\cS_\i}$ is the $\i$-ratio function of $\cS_\i$,
$\rho_{\nu,\i}$ is the $\i$-measure ratio function of
the Gibbs measure $\nu$,
and $\rho_{\nu,\ip}$ is the $\ip$-measure ratio function of
the Gibbs measure $\nu$.
\end{lemma}

\begin{proof}
For every $x \in \cP$,
let us denote by $p$     the period $p(x)$  of  $x$, and
let us denote by $I^\i_{n}$ the interval $I^\i_{n}(x)$.
We note that the $p$-mother $m^p I^\i_{n+1}$ of
$I^\i_{n+1}$ is $I^\i_n$, and so
$f^p_\i I^\i_{n+1}=m^p I^\i_{n+1}$.
By \eqref{eqn:12133},
$$
r_{\cS_\i} (I^\i_{0}:I^\i_{1}) =  \lim_{n \to \infty}
\frac{|I^\i_{n}|}  {|I^\i_{n+1}|} \ .
$$ Hence,
\begin{eqnarray*}
\lambda_{\cS_\i}(x_\i)
& = & \lim_{n \to \infty}
\frac{|f^p_\i I^\i_{n+1}|}  {|I^\i_{n+1}|} \\
& = & \lim_{n \to \infty}
\frac{|I^\i_{n}|}  {|I^\i_{n+1}|} \nonumber  \\
& = & r_{\cS_\i} (I^\i_{0}:I^\i_{1})
\end{eqnarray*}
which proves \eqref{grgrgr000}.
By Theorem \ref{dsgdgnnww}, the
$\i$-measure ratio function  $\rat_{\cS_\i,\delta_\i}$ is
the $\i$-measure ratio function  $\rho_{\nu,\i}$
of the Gibbs measure $\nu$. Hence, by \eqref{fdfdbng}, we get
\begin{eqnarray}
\label{grgrgr35}
r_{\cS_\i} (I^\i_{1}:I^\i_{0}) & = &
\prod_{l=0}^{p-1} r_{\cS_\i} (m^{l} I^\i_{1}:m^{l+1} I^\i_{1}) \nonumber \\
& = &
\prod_{l=0}^{p-1} \left( J_{\cS,\delta_\i} (\xi_{l})
\rho_{\nu,\i}(m^{l} I^\i_{1}:m^{l+1} I^\i_{1})^{1/\d_\i} e^{P_\i/\d_\i} \right)
\end{eqnarray}
where $\xi_{l}=   f_\i^{p-l} m^l I^\i_{1} \in \TT^\i_o$.
We note that    $f_\ip \xi_{l}= \xi_{l+1}$  and $f_\ip \xi_{p-1}= \xi_{0}$ in $\TT_o^\ip$.
Since $J_{\cS_\i,\delta_\i}= \kappa/(\kappa\comp f_\ip)$ for some function $\kappa$,
we get
\begin{equation}
\label{grgrgr3}
\prod_{l=0}^{p-1} J_{\cS_\i,\delta_\i} (\xi_{l}) =
\prod_{l=0}^{p-1} \frac{\kappa (\xi_{l})}{\kappa (\xi_{l+1})}
= 1 \ .
\end{equation}
Furthermore,
\begin{equation}
\label{grgrgr5}
\prod_{l=0}^{p-1} \rho_{\nu,\i} (m^{l} I^\i_{1}:m^{l+1} I^\i_{1})
= \rho_{\nu,\i} (I^\i_{1}:I^\i_{0})  \ .
\end{equation}
Using \eqref{grgrgr3} and \eqref{grgrgr5} in \eqref{grgrgr35}
we obtain that
$$
r_{\cS_\i} (I^\i_{1}:I^\i_{0}) =
\rho_{\nu,\i} (I^\i_{1}:I^\i_{0})^{1/\d_\i} e^{p P_\i/\d_\i} \ .
$$
Therefore, by \eqref{grgrgr000}, we have
\begin{eqnarray*} \label{grgrgr05}
\lambda_{\cS_\i}(x_\i) & = &
r_{\cS_\i} (I^\i_{0}:I^\i_{1}) \\
& = &
\rho_{\nu,\i} (I^\i_{0}:I^\i_{1})^{-1/\d_\i} e^{-p P_\i/\d_\i}
\end{eqnarray*}
which proves \eqref{grgrgr001}.
By Lemma \ref{lemma:ratio}, there is
$0 < \varepsilon < 1$ such that  for every $n \ge 0$
\begin{equation}
\label{grgrgr13331}
\rho_{\nu,s} (I^s_{n+1}:I^s_{n})
\in
(1 \pm \varepsilon^{n})
\frac{\mu ([I^s_{n+1},I^u_{1}]) }
{\mu ([I^s_{n},I^u_{1}]) }
\end{equation}
and
\begin{equation}
\label{grgrgr13333}
\rho_{\nu,u} (I^u_{n+1}:I^u_{n})
\in
(1 \pm \varepsilon^{n})
\frac{\mu ([I^s_{1},I^u_{n+1}]) }
{\mu ([I^s_{1},I^u_{n}]) } \ .
\end{equation}
Since  $f^{np}\left([I^s_1,I^u_{n+1}]\right)=
\left([I^s_{n+1},I^u_{1}]\right)$ and
by invariance of $\mu$, we obtain that
\begin{equation}
\label{grgrgr17}
\frac{\mu ([I^s_{n+1},I^u_{1}]) }
{\mu ([I^s_{n},I^u_{1}])  }
= \frac{\mu ([I^s_1,I^u_{n+1}])}
{\mu ([I^s_{1},I^u_{n}])}
\end{equation}
Putting together \eqref{grgrgr13331}, \eqref{grgrgr13333}
and \eqref{grgrgr17}, we obtain that
$$
\rho_{\nu,s} (I^s_{n+1}:I^s_{n}) \in
(1 \pm \varepsilon^{n}) \rho_{\nu,u} (I^u_{n+1}:I^u_{n}) \ .
$$
Hence, by invariance of $\rho_{\cS_s,s}$ and $\rho_{\cS_s,u}$ under $f$,
we obtain
\begin{eqnarray*}
\label{grgrgr23}
\rho_{\nu,s} (I^s_1:I^s_0)
& = &  \lim_{n \to \infty} \rho_{\nu,s} (I^s_{n+1}:I^s_{n})  \\
& = &  \lim_{n \to \infty} \rho_{\nu,u} (I^u_{n+1}:I^u_{n})  \\
& = &   \rho_{\nu,u} (I^u_1:I^u_0)
\end{eqnarray*}
which proves \eqref{grgrgr003}.
\end{proof}

\begin{lemma}
\label{akju2}
Let $\TT^\i$ be a (gap or a no-gap) train-track.
\begin{rlist}
\item
The $C^{1+}$
self-renormalisable structures $\cS_1 \in \cD^\i(\nu_1,\delta,P)$
and  $\cS_2 \in \cD^\i(\nu_2,\delta,P)$
have  the same   eigenvalues for all
periodic orbits  if, and only if,
$\nu_1$ is equal to $\nu_2$.
\item
The set of eigenvalues of all periodic orbits of a
$C^{1+}$ self-renormalisable structure
is  a complete invariant of  each
Lipschitz conjugacy class.
\end{rlist}
\end{lemma}

Statement (ii) of the above lemma for Markov maps is also in
\cite{sullivan:nested}.

\noindent
\emph{Proof of  (i)}.
By Lemma \ref{akju1},
the $C^{1+}$ self-renormalisable structures $\cS_1 \in \cD^\i(\nu_1,\delta,P)$
and $\cS_2 \in \cD^\i(\nu_2,\delta,P)$ are Lipschitz conjugate
if, and only if,
the Gibbs measures $\nu_1$ and $\nu_2$ are equal.
By Lemma \ref{fgregtgeeew},
if the Gibbs measures $\nu_1$ and $\nu_2$ are equal,
then
$\cS_1$ and $\cS_2$ have the same eigenvalues
for all periodic orbits.
Hence,
to finish the proof of statement (i),
we are going to prove that if
the $C^{1+}$
self-renormalisable structures $\cS_1$
and  $\cS_2$
have  the same   eigenvalues for all
periodic orbits,
then the $C^{1+}$
self-renormalisable structures $\cS_1$
and  $\cS_2$ are Lipschitz conjugate.

Without loss of generality,
let us assume that $\cS_1$
and  $\cS_2$ are unstable $C^{1+}$
self-renormalisable structures.
For $j \in \{1,2\}$, the {\emph (restricted) $u$-scaling function}
$z_{u,j}: \Theta^u \to \reals^+$ of $\cS$ is well-defined by
(see \S ~\ref{foldjgjlke})
$$z_{u,j} (w_0 w_1 \ldots) = \lim_{n \to \infty}
\frac{|\pi_{\TT^s} \circ f^{n+1} \circ \pi^{-1}_{\TT^u} \circ i_u (w_0 w_1 \ldots)|_{k_n}}
{|\pi_{\TT^s} \circ f^{n} \circ \pi^{-1}_{\TT^u} \circ i_u (w_1 w_2 \ldots)|_{k_n}},$$
where $k_n$ is a train-track chart in a $C^{1+}$ self-renormalisable atlas $\cB_j$
determined by $\cS_j$ such that the domain of the chart $k_n$ contains
$\pi_{\TT^s} \circ f^{n} \circ \pi^{-1}_{\TT^u} \circ i_u (w_1 w_2 \ldots)$.
For every stable-leaf $(i+1)$-cylinder $J$,
let $w(J) \in \Theta^u$ be such that
$i_u(w(J)) = \pi_{\TT^u}(J)$.
Hence, for every $l \in \{0, \ldots, i-1\}$ we have that
$$\pi_{\TT^u}^{-1} \circ i_u(f_u^l w(J)) = f^{-i+l}( m^l J) \ ,$$
where $f^{-i+l}( m^l J)$ are stable-leaf primary cylinders.
By construction of the (restricted) $u$-scaling function $z_{u,j}$ and of
the $u$-scaling function $s_{u,j}$ of $\cS_j$,
we have that
\begin{equation}
\label{tryri10}
s_{u,j}(J:m^i J) = \prod_{l=0}^{i-1} z_{u,j} (f_u^l (w(J))) \ .
\end{equation}
Let $\cP_{\Theta^u}$ be the set of all periodic point
under the shift.
For every $w = w_0 w_1 \ldots \in \cP_{\Theta^u}$ let $p(w)$ be the smallest period of $w$.
By construction of the train-tracks,
for every $w$, there is a unique periodic point $x(w) \in \Lambda$ with period $p(w)$
with respect to the map $f$
such that $i_u(w) = \pi_{\TT^u} x(w)$.
Furthermore, there is a unique periodic point $\pi_sx(w) \in \TT^s$ with period $p(w)$ for the Markov map.
By \eqref{tryri10},
for every $w \in \cP_{\Theta^u}$
we have that
\begin{equation}
\label{aalli}
\prod_{i=0}^{p(w)-1} z_{u,j} (f_u^i(w)) = \lambda_{\cS_j} (\pi_{\TT^s} x(w)) \ .
\end{equation}
Since the $C^{1+}$
self-renormalisable structures $\cS_1$
and  $\cS_2$
have  the same   eigenvalues  for all
periodic orbits, by \eqref{aalli}, we have that
\begin{equation}
\label{aallir}
\prod_{i=0}^{p(w)-1} \frac{z_{u,1} (f_u^i(w))}{z_{u,2} (f_u^i(w))} = 1 \ ,
\end{equation}
for every $w \in \cP_{\Theta^u}$.
From Liv\v sic's
theorem (e.g.\ see
\cite{KatokHassBook}) we get that
\begin{equation}
\label{aallir3}
\frac{z_{u,1} (w)}{z_{u,2} (w)} =\frac{\kappa(w)}{\kappa\comp f_u (w)}
\end{equation}
where $\kappa:\Theta^u \to \reals^+$ is a
positive H\"older continuous
function.
By \eqref{tryri10} and \eqref{aallir3},
for every stable-leaf $(i+1)$-cylinder $J$
we obtain that
\begin{eqnarray}
\label{tryri}
\frac{ s_{u,1}(J:m^i J)}{s_{u,2}(J:m^i J)} & = &
\prod_{l=0}^{i-1} \frac{z_{u,1} (f_u^l (w(J)))}{{z_{u,2} (f_u^l (w(J)))}} \nonumber \\
& = & \frac{\kappa(w)}{\kappa \comp f_u^i (w)}\ .
\end{eqnarray}
Since $\kappa$ is bounded away from zero and infinity,
there is $C>1$ such that for all $w \in \Theta^u$ and $i \ge 1$
we have that
\begin{equation}
\label{aallir5}
C^{-1} < \frac{\kappa(w)}{\kappa\comp \tau_u^i (w)} < C \ .
\end{equation}
Putting together \eqref{tryri} and \eqref{aallir5}, we obtain that
$$\frac{ s_{u,1}(J:m^i J)}{s_{u,2}(J:m^i J)} \ . $$
Therefore,
the $\i$-solenoid functions $\s_{u,1}:\sol^\i \to \reals^+$
and $\s_{u,2}:\sol^\i \to \reals^+$ corresponding to the $C^{1+}$
self-renormalisable structures $\cS_1$
and  $\cS_2$
are in the same bounded equivalence class
(see Definition \ref{fgrdegerdg}).
Hence,
by Lemma \ref{akju1},
the $C^{1+}$ self-renormalisable structures $\cS_1$ and $\cS_2$ are Lipschitz conjugate.

\noindent
\emph{Proof of  (ii)}.
Statement (ii) follows from
putting together  Lemma \ref{akju1} and  Statement (i) of this lemma with $P=0$.
\qed

\begin{theorem}
\label{fgdgfgdds3}
\begin{rlist}
\item If $g \in \cT(\fL )$ then
$\lambda_{g,s}(x)= \lambda_{\cS_s}(x_s)^{-1}$ and
$\lambda_{g,u}(x)=\lambda_{\cS_s}(x_u)$
where, for $\i \in \{s,u\}$, $\lambda_{g,\i}(x)$
is the  eigenvalue of
the $C^{1+}$ hyperbolic diffeomorphism  $g$
and $\lambda_{\cS_\i}$
is the  eigenvalue of
the $C^{1+}$ self-renormalisable structure $\cS_\i=\cS(g,\i)$.
\item
The set of stable and unstable eigenvalues of all periodic orbits of a
$C^{1+}$ hyperbolic diffeomorphism
$g \in \cT(\fL )$
is  a complete invariant of  each
Lipschitz conjugacy class.
\end{rlist}
\end{theorem}

\noindent
\emph{Proof of (i)}.
By construction of the train-track atlas
$\cB ^\i (g,\rho_g)$    from the
lamination atlas $\cA ^\i (g,\rho_g)$   in \S ~\ref{fddhhhuy},
if $\lambda_{g,\i}(x)$ is the eigenvalue of $x \in \cP$
then the eigenvalue of $x_\i \in \cP_\i$
is either $\lambda_{g,\i}(x)$ if $\i=u$,
or  $\lambda_{g,\i}^{-1}$ if $\i=s$.

\noindent \emph{Proof of (ii)}.
By  Lemma \ref{akju2},
the set of eigenvalues of all periodic orbits
of a $C^{1+}$
self-renormalisable   structure is
a complete invariant of each
Lipschitz conjugacy class of $C^{1+}$
self-renormalisable   structures.
Hence, using Lemma \ref{asdro},
we get that
the set of stable and unstable
eigenvalues of all periodic orbits of a
$C^{1+}$ hyperbolic diffeomorphism $g$
is  a complete invariant of  each
Lipschitz conjugacy class.
\qed

\subsection{Extending the eigenvalue formula of A. N.  Liv\v sic  and Ja. G. Sinai}

In Lemma \ref{fgregtgeeew}, we give an explicit formula for the stable and unstable eigenvalues of all
periodic points of the $C^{1+}$ hyperbolic diffeomorphisms in terms of measure ratio functions.
We use this formula, in Theorem \ref{fgdgfgdds31} below,
to extend the eigenvalue formula of A. N.  Liv\v sic  and Ja. G. Sinai for Anosov diffeomorphisms
to $C^{1+}$ hyperbolic diffeomorphisms.

\begin{theorem}
\label{fgdgfgdds31}
A $C^{1+}$ hyperbolic diffeomorphism $g \in \cT(\fL )$
has a natural geometric measure $\mu_{g,\delta_s,\delta_u}$ with
pressures $P_s=P_s(g,\delta_s,\delta_u)$ and
$P_u=P_u(g,\delta_s,\delta_u)$ if, and only if,
for all $x \in \L$
\begin{equation}
\label{aqdfdfd5}
\lambda_{g,s}(x_s)^{-\delta_s} e^{p(x)P_s} =
\lambda_{g,u}(x_u)^{\delta_u} e^{p(x)P_u}   \   .
\end{equation}
\end{theorem}

\begin{proof}
By Theorem \ref{thm:SRB},
the $C^{1+}$ self-renormalisable structures $S(g,s)$ and $S(g,u)$
are $C^{1+}$ realisations of  Gibbs measures  $\nu_1=\nu_{S(g,s),\delta_s}$ and $\nu_2=\nu_{S(g,u),\delta_s}$.
By Lemma \ref{fgregtgeeew} and the statement (i) of Theorem \ref{fgdgfgdds3},
for all $x \in \cP$
we have
\begin{eqnarray}
\label{aqdfdfd1}
\lambda_{g,u}(x_u)   &= & \lambda_{S(g,u)}(x_u) \nonumber \\
&= &
\rho_{\nu_2, u} (I^u_{0}:I^u_{1})^{-1/\d_u} e^{-p(x) P_u/\d_u}
\end{eqnarray}
and
\begin{eqnarray}
\label{aqdfdfd3}
\lambda_{g,s}(x_s)   &= & \lambda_{S(g,s)}(x_s)^{-1} \nonumber \\
&= & \rho_{\nu_1, u} (I^u_{0}:I^u_{1})^{1/\d_s} e^{p(x) P_s/\d_s} \ .
\end{eqnarray}
Let us prove that if  the $C^{1+}$ hyperbolic diffeomorphism $g$
has a natural geometric measure then \eqref{aqdfdfd5} holds.
Hence, by Lemma \ref{aasskk}, the Gibbs measures  $\nu_1$ and $\nu_2$ are equal.
By \eqref{aqdfdfd1}, we have
\begin{eqnarray}
\label{aqdfdfd35r}
\rho_{\nu_1, u} (I^u_{0}:I^u_{1}) & = & \rho_{\nu_2, u} (I^u_{0}:I^u_{1}) \nonumber \\
& = &  \lambda_{g,u}(x_\i)^{-\d_u}  e^{-p(x) P_u} \ .
\end{eqnarray}
By \eqref{aqdfdfd3}, we obtain that
\begin{equation}
\label{aqdfdfd35fr}
\rho_{\nu_1, u} (I^u_{0}:I^u_{1}) = \lambda_{g,s}(x_\i)^{\d_s}  e^{-p(x) P_s} \ .
\end{equation}
Putting together \eqref{aqdfdfd35r}
and \eqref{aqdfdfd35fr}, we obtain that
$$
\lambda_{g,s}(x_s)^{-\delta_s} e^{p(x)P_s} =
\lambda_{g,u}(x_u)^{\delta_u} e^{p(x)P_u}   \ ,
$$
and so \eqref{aqdfdfd5} holds.
Conversely, let us prove that if  \eqref{aqdfdfd5} holds then
  the $C^{1+}$ hyperbolic diffeomorphism $g$
has a natural geometric measure.
Putting together  \eqref{aqdfdfd5} and \eqref{aqdfdfd3}, we obtain that
$$
\lambda_{g,u}(x_u)
=
\rho_{\nu_1, u} (I^u_{0}:I^u_{1})^{-1/\d_u} e^{-n P_u/\d_u} \ .
$$
Hence, the Gibbs measure $\nu_1$     determines the same set of eivenvalues
for all periodic orbits of self-renormalisable structures in $\TT^u$
as the Gibbs measure $\nu_2$. Therefore,
by Lemma \ref{akju2}, $\nu_1 =\nu_2$ and consequently, by Lemma   \ref{aasskk},
the $C^{1+}$ hyperbolic diffeomorphism $g$
has a natural geometric measure.
\end{proof}

\section{Invariant Hausdorff measures}
\label{fghtghrr}

In this section, we present the proofs of all theorems stated in the Introduction.
Let $\cS_\i$ be a $C^{1+}$
$\i$ self-renormalisable structure.
By Lemma \ref{cor_geommeas}, a  natural geometric measure $\mu_{\cS_\i,\delta_\i}$
with pressure $P(\cS_\i,\delta_\i)=0$ is an invariant measure absolutely continuous
with respect to the Hausdorff measure of $\TT^\i$   and
$\delta_\i$ is the Hausdorff dimension of $\TT^\i$ with respect to the charts of $\cS_\i$.
Let us denote $\cD^\i(\nu,\delta_\i,0)$  and  ${\mathcal JG}^\i(\nu,\delta_\i,0)$ respectively
by $\cD^\i(\nu,\delta_\i)$ and  ${\mathcal JG}^\i(\nu,\delta_\i)$.
By Theorem \ref{thm:SRB},  for every $C^{1+}$
$\i$ self-renormalisable structure $\cS_\i$ there is
an unique  Gibbs measure $\nu_{\cS_\i}$
such that $\cS_\i \in  \cD^\i(\nu,\delta_\i)$.
Using Lemma \ref{fghbhhjbsde},
we obtain that the sets
$[\nu] \subset  \cTf( \delta_s,\delta_u)$
defined in the introduction
are equal  to the sets
$\cT(\nu,\delta_s,0,\delta_u,0)$ (see Definition \ref{fgtyjtyjyu}).

Theorem \ref{dfsdfaaa1} follows from Theorem \ref{fgdgfgdds31},
Theorem \ref{gthrthrtr}  follows from Theorem \ref{fghhhytr},
Theorem \ref{pfiejd} follows from   Lemma \ref{ghjjjyuu}, and
Theorem \ref{rfgrghhbvcde} follows from   Lemma \ref{fghbhhjbsde}.

\bigskip
\noindent
\emph{Proof of Theorem \ref{dfsdfaaa3}}.
\emph{Proof of  statement (i)}.
By Lemma \ref{fgdgfgdds1} (ii), the sets
$[\nu]\subset \cTf(\delta_s,\delta_u)$ are
Lipschitz conjugacy classes in
$ \cTf( \delta_s, \delta_u)$, and the map
$\nu \rightarrow \cT(\nu, \delta_s,\delta_u)$ is injective.
If $g \in  \cTf( \delta_s,\delta_u)$
then $g$ has a natural geometric measure $\mu_{g,\delta_s,\delta_u}$
with pressures $P_s(g,\delta_s,\delta_u)$ and
$P_u(g,\delta_s,\delta_u)$ equal to zero.
By Lemma \ref{aasskk}, there is a Gibbs measure
$\nu=\nu_{g,\delta_s,\delta_u}$ on $\Theta$
such that $i_*\nu=\mu_{g,\delta_s,\delta_u}$ and so
$g \in [\nu]\subset \cTf(\delta_s,\delta_u)$.
Hence, the map $\nu \rightarrow \cT(\nu, \delta_s,\delta_u)$
is surjective into the Lipschitz  conjugacy classes
in $ \cTf( \delta_s,\delta_u)$.

\noindent
\emph{Proof of  statement (ii)}.
By  Theorem \ref{fgdgfgdds3} (ii),
the set of stable and unstable eigenvalues of all periodic orbits of a
$C^{1+}$ hyperbolic diffeomorphisms
$g \in  \cTf( \delta_s,\delta_u)$
is  a complete invariant of  each
Lipschitz conjugacy class, and by statement (i)
of this lemma the sets $\cT(\nu, \delta_s,\delta_u)$
are the Lipschitz conjugacy classes in $ \cTf( \delta_s,\delta_u)$.
\qed

\bigskip
\noindent
\emph{Proof of Theorem  \ref{self1}}.
We will separate the proof in three parts.
In part (i), we prove that if $\cS_\i \in \cD^\i(\nu,\delta_\i)$
then $\s_{\nu,\i}$ satisfies  the properties indicated in   Theorem  \ref{self1}.
In part (ii), we prove the converse of part (i). In part (iii),
we prove that
$\cD^\i(\nu,\delta_\i) \ne \emptyset$ if, and ony if,
$\cD^\ip(\nu,\delta_\ip) \ne \emptyset$.

\noindent
\emph{Part (i).}
Let $\cS_\i \in \cD^\i(\nu,\delta_\i)$.
By Lemma \ref{dsgdgnnww},  $\cS_\i $ and $\delta_\i$ determine an unique   $\i$-measure ratio
function $\rho_{\nu,\i}$ of the Gibbs measure $\nu$.
Hence,  the function $\rho_{\nu,\i}|\msol^\ip$
is the $\i$-measure solenoid function  $\s_{\nu,\i}$
of $\nu$ and, by Lemma \ref{altum},
$\s_{\nu,\i}$ satisfies  the properties indicated in   Theorem  \ref{self1}.

\noindent
\emph{Part (ii).}
Conversely, if $\nu$ has an $\i$-solenoid function $\s_{\nu,\i}$
satisfying the properties indicated in   Theorem  \ref{self1},
by lemmas \ref{fghbhbfw3} and  \ref{altum}, $\s_{\nu,\i}$ determines an unique
$\i$-measure ratio function $\rho_{\nu,\i}$ of $\nu$.
If $\TT^\i$ is a   no-gap train-track, by Lemma \ref{trdfg},
there is a $C^{1+}$ self-renormalisable structure $\cS_\i \in \cD^\i(\nu,\delta_\i)$
with $\delta_\i=1$.
If $\TT^\i$ is a   gap train-track then, by Remark \ref{dsfdgdhghv},
the set ${\mathcal JG}^\i(\nu,\delta_\i)$ is non-empty (in fact
it is an   infinite dimensional space).
Hence, by Lemma \ref{gdvaswwa}, the set $\cD(\nu,\delta_\i)$ is also non-empty
which ends the proof.

\noindent
\emph{Part (iii).}
To prove that  $\cD^\i(\nu,\delta_\i) \ne \emptyset$ if, and ony if,
$\cD^\ip(\nu,\delta_\ip) \ne \emptyset$, it is enough to prove one of the implications.
Let us prove that if
$\cD^\i(\nu,\delta_\i) \ne \emptyset$  then $\cD^\ip(\nu,\delta_\ip) \ne \emptyset$,
Let $\cS_\i \in \cD^\i(\nu,\delta_\i)$.
By Lemma \ref{dsgdgnnww},  $\cS_\i $ and $\delta_\i$ determine an unique   $\i$-measure ratio
function $\rho_{\nu,\i}$ of the Gibbs measure $\nu$.
By Lemma \ref{ploky},
the $\i$-measure
ratio function $\rho_{\nu,\i}$
determines an unique dual $\ip$-measure
ratio function $\rho_{\nu,\ip}$ of $\nu$.
Hence,  the function $\rho_{\nu,\ip}|\msol^\ip$
is the $\i$-measure solenoid function  $\s_{\nu,\ip}$
of $\nu$ and, by Lemma \ref{altum},
$\s_{\nu,\ip}$ satisfies  the properties indicated in   Theorem  \ref{self1}.
Now the proof follows as in part (ii), with $\i$ changed by $\ip$,
which shows that $\sigma_{\nu,\ip}$ determines a non-empty set
$\cD^\ip(\nu,\delta_\ip)$.
\qed

\bigskip
\noindent
\emph{Proof of Theorem \ref{anosov1} and Theorem  \ref{uyevvgd}}.
\emph{Proof that  statement (i) implies statements (ii) and (iii)}.
If $g \in [\nu]\subset \cTf(\delta_s,\delta_u)$,
by Lemma \ref{fghbhhjbsde}, the sets
$\cD^s (\nu,\delta_s)$ and $\cD^u (\nu,\delta_u)$
are both non-empty.
Hence, by  Theorem   \ref{self1},
the stable measure solenoid function of
the Gibbs measure $\nu$ satisfies (ii) and
the unstable measure solenoid function of
the Gibbs measure $\nu$ satisfies (iii).

\noindent
\emph{Proof that statement (ii) implies statement (i),
and that  statement (iii) implies  statement (i)}.
By  Theorem   \ref{self1},
the properties of the the $\i$-solenoid function $\s_{\nu,\i}$
indicated in this theorem imply that
$\cD^\i (\nu,\delta_\i) \ne \emptyset$.
Again, by  Theorem   \ref{self1} and $\cD^\ip (\nu,\delta_\ip)\ne \emptyset$.
Hence, by Lemma \ref{fghbhhjbsde}, the set $[\nu]\subset \cTf(\delta_s,\delta_u)$
is non-empty. Therefore, every $g \in \cT(\nu, \delta_s,\delta_u)$
is a $C^{1+}$-Hausdorff realisation of $\nu$ which ends the proof.
\qed

\bigskip
\noindent
\emph{Proof of Theorem \ref{htrytgvg}}.
Let $\nu$ be a  Gibbs measure.
By  Theorem  \ref{self1}, the set $\cD^s (\nu,\delta_s)$
and   $\cD^u (\nu,\delta_u)$ are both non-empty.
Hence, by Lemma \ref{fghbhhjbsde},  the set
$\cT(\nu, \delta_s,\delta_u)$ is also non-empty.
Therefore, every $g \in \cT(\nu, \delta_s,\delta_u)$
is a $C^{1+}$-Hausdorff realisation of $\nu$ which ends the proof.
\qed

\subsection{Moduli space ${\mathcal SOL}^\i$}

Recall the definiton of the set ${\mathcal SOL}^\i$ given in \S ~\ref{fsefsefff}.
By Theorem \ref{83dfsf}, below,
the set of all $\i$-measure solenoid functions $\s_\nu$
with the properties indicated in  Theorem  \ref{self1}
determine an  infinite dimensional metric  space ${\mathcal SOL}^\i$ which gives
a nice parametrization of all  Lipschitz conjugacy classes
$\cD^\i(\nu,\delta)$
of    $C^{1+}$ self-renormalisable structures $\cS_{\i}$ with
a given Hausdorff dimension $\delta$.

Theorem \ref{333dfsf}  follows from Theorem \ref{83dfsf}.

\begin{theorem}
\label{83dfsf}
If $\TT^\i$ is a gap train-track assume $0< \delta_\i <1$
and
if $\TT^\i$ is a no-gap train-track assume  $\delta_\i =1$.
\begin{rlist}
\item
The map $\cS \to \rat_{\cS,\delta_\i}$
induces an one-to-one correspondence between the sets
$\cD^\i(\nu,\delta_\i)$ and the elements of ${\mathcal SOL}^\i$.
\item
The map $g \to \rat_{\cS(g,\i),\delta_{\i}}$
induces an one-to-one correspondence between the sets
$[\nu]$ contained in  $ \cTf(\delta_s,\delta_u)$
and the elements of  ${\mathcal SOL}^\i$.
\end{rlist}
\end{theorem}

\bigskip
\noindent
\emph{Proof of Theorem \ref{83dfsf}}.
\emph{Proof of (i).}
If $\cS \in \cD^\i(\nu,\delta_\i)$ then the Hausdorff dimension
of $\cS$ is $\delta_\i$,
and $\cS$ determines an $\i$-measure
ratio function $\rho_{\cS,\delta_\i}=\rho_{\nu,\i}$
which does not depend upon $\cS \in \cD^\i(\nu,\delta_\i)$.
By Lemma \ref{altum}, $\rho_{\nu,\i}|\Msol^\i$
is an element of ${\mathcal SOL}^\i$.
Hence,
the map $\cS \to \rat_{\cS,\delta_\i}$
associates to each set
$\cD^\i(\nu,\delta)$ an unique  element of ${\mathcal SOL}^\i$.
Conversely, let $\hat{\sigma} \in \Msol^\i$.
By Lemma \ref{altum}, $\hat{\sigma}$ determines an unique
$\i$-measure ratio function $\rho_{\hat{\sigma}}$ such that
$\rho_{\hat{\sigma}}|\Msol^\i=\hat{\sigma}$.
By Lemma   \ref{fgdfgbbcc3}, the
$\i$-measure ratio function $\rho_{\hat{\sigma}}$
determines a Gibbs measure  $\nu_{\hat{\sigma}}$.
If $\TT^\i$ is a no-gap train-track then, by Lemma \ref{trdfg},
$\rho_{\hat{\sigma}}$ determines a non-empty set
$\cD^\i(\nu_{\hat{\sigma}},\delta_\i)$.
If $\TT^\i$ is a gap train-track then,
by Remark \ref{dsfdgdhghv},
the set  ${\mathcal JG}^\i(\nu,\delta_\i)$ is non-empty
and so, by Lemma \ref{gdvaswwa}, the set $\cD^\i(\nu_{\hat{\sigma}},\delta_\i)$ is
also non-empty. Therefore, each element $\hat{\sigma} \in \Msol^\i$
determines an unique non-empty set $\cD^\i(\nu_{\hat{\sigma}},\delta_\i)$
of $C^{1+}$ self-renormalisable structures $\cS$
with $\rat_{\cS,\delta_\i}|\Msol^\i=\hat{\sigma}$.

\noindent
\emph{Proof of (ii).}
By Lemma \ref{aasskk}, if $g \in [\nu]$
then $\cS(g,\i) \in \cD^\i(\nu,\delta_\i)$ and so,
by statement (i) of this lemma,
$\rat_{\cS(g,\i),\delta_{\i}}|\Msol^\i$ is an element of  ${\mathcal SOL}^\i$
which does not depend upon $g \in [\nu]$.
Conversely, let $\hat{\sigma} \in \Msol^\i$.
By statement (i) of this lemma,  $\hat{\sigma}$
determines an $\i$-measure ratio function $\rho_{\hat{\sigma},\i}$,
and
a non-empty set $\cD^\i(\nu_{\hat{\sigma}},\delta_\i)$.
By Lemma \ref{aasskk}, $\rho_{\hat{\sigma},\i}$ determines an unique dual
$\ip$-ratio function $\rho_{\hat{\sigma},\ip}$ associated to
the Gibbs measure $\nu_{\hat{\sigma}}$.
Again, by statement (i) of this lemma,
$\rho_{\hat{\sigma},\ip}|\Msol^\ip$ determines
a non-empty set $\cD^\ip(\nu_{\hat{\sigma}},\delta_\ip)$.
By Lemma \ref{fghbhhjbsde}, the set
$\cD^s(\nu_{\hat{\sigma}},\delta_s) \times \cD^u(\nu_{\hat{\sigma}},\delta_u)$ determines
an unique  non-empty set
$ [\nu_{\hat{\sigma}}]\subset \cTf(\delta_s,\delta_u)$
of hyperbolic diffeomorphisms $g \in [\nu_{\hat{\sigma}}]$
such that
$\rat_{\cS(g,\i),\delta_{\i}}|\Msol^\i=\hat{\sigma}$.
\qed

\subsection{Moduli space of cocycle-gap pairs}

By Lemma \ref{akju1}, each set $\cD^\i(\nu,\delta)$ is
a Lipschitz conjugacy class. Hence, by Theorem \ref{dfgrggrrd} proved below, if
$\TT^\i$ is a no-gap train-track then the
Lipschitz conjugacy class consists of  a
single   $C^{1+}$ self-renormalisable structure.
Furthermore, by Lemma \ref{akju2},  the set of   eigenvalues of all
periodic orbits of $\cS_\i$
is a complete invariant of
each set    $\cD^\i(\nu,\delta)$.

\bigskip
\noindent
\emph{Proof of Theorem \ref{dfgrggrrd}}.
Statement (i) follows from Lemma \ref{trdfg}. Now, let us prove statement (ii).
By Remark \ref{dsfdgdhghv},
the set ${\mathcal JG}^\i(\nu,\delta)$ is an   infinite dimensional space,
and  by  Lemma \ref{gdvaswwa}, the set $\cD^\i(\nu,\delta)$ is
parametrized by the cocycle-gap pairs in ${\mathcal JG}^\i(\nu,\delta)$
which ends the proof.
\qed

\bigskip
\noindent
\emph{Proof of Theorem  \ref{fdgdrgrrqs}}.
By Lemma \ref{fghbhhjbsde}, if $g \in \cT(\nu, \delta_s,\delta_u)$
then $\cS_\ip(g) \in\cD^\ip (\nu,\delta_\ip)$.
Conversely, let $\cS_\ip$ be a $C^{1+}$ self-renormalisable structure contained in
$\cD^\ip (\nu,\delta_\ip)$.
By Lemma \ref{fghbhhjbsde},
a pair $(\cS_\i,\cS_\ip)$   determines
a $C^{1+}$  hyperbolic diffeomorphism $g \in \cT(\nu, \delta_s,\delta_u)$.
if, and only if, $\cS_\ip \in \cD^\ip (\nu,\delta_\ip)$.
By  Theorem  \ref{self1}, the
set $\cD (\nu,\delta_\i)$ is non-empty.
Noting that
$\delta_\i=1$,  it follows from  Theorem \ref{dfgrggrrd} (ii) that
the set $\cD^\i (\nu,\delta_\i)$ contains only
one
$C^{1+}$ self-renormalisable structure $\cS_\i$ which finishes the proof.
\qed

\subsection{$\delta_\i$-bounded solenoid equivalence class of Gibbs measures}
When we speak of a $\delta_\i$-bounded solenoid equivalence class of  $\nu$
we mean a $(\delta_\i,0)$-bounded solenoid equivalence class  of a Gibbs measure $\nu$ (see Definition  \ref{dfgdrfgggx}).
In \S ~\ref{dfgggvvd}, we use the cocycle-gap pairs
to construct explicitly the solenoid functions in the
$\delta_\i$-bounded solenoid equivalence classes of the Gibbs measures $\nu$.
By Theorem \ref{ddfsffe1} (ii) proved below,
given an $\i$-solenoid function $\s_\i$
there is an unique  Gibbs measure $\nu$
such that $\sigma_\i$ belongs to the $\delta_\i$-bounded
solenoid equivalence class of  $\nu$.

\bigskip
\noindent
\emph{Proof of Theorem  \ref{ddfsffe1}}.
Statement (i) follows from Lemma \ref{fgdgfgdds1} (i).
Statement (ii) follows  from    Lemma \ref{trdfg} if $\TT^\i$ is a no-gap train-track,
and  from  Lemma \ref{gdvaswwa} (ii) if $\TT^\i$ is a  gap train-track.
Statement (iii) follows from statement (ii) and Theorem \ref{dfgrggrrd}.
\qed

\bigskip
\noindent
\emph{Proof of Theorem  \ref{ddfsffe3}}.
By Theorem \ref{ddfsffe1} (ii), the $\i$-solenoid function $\s_\i$ determines
an unique $C^{1+}$    self-renormalisable  structure $\cS_\i \in \cD^\i(\nu,\delta_\i)$.
By  Theorem  \ref{self1}, the set $\cD^\ip(\nu,\delta_\ip)$
is nonempty.
Let $\cS_\ip \in \cD^\ip(\nu,1)$.
By Theorem  \ref{ddfsffe1} (ii),
the $C^{1+}$    self-renormalisable  structure
$\cS_\ip$ determines an unique $\ip$-solenoid function
$\s_\ip$ such that, by Theorem \ref{ddfsffe1} (i),
the pair $(\s_\i,\s_\ip)$ determines an unique
$C^{1+}$ conjugacy
class $\cT(\nu, \delta_s,\delta_u)$
of hyperbolic diffeomorphisms $g \in  \cT(\nu, \delta_s,\delta_u)$
with an invariant measure
$\mu=i_* \nu$ absolutely continuous with respect to
the Hausdorff measure.
\qed

%>>>>>>>>>>>>>> fandr_refs.tex <<<<<<<<<<<<<<<<
\par\bigskip \noindent {\Large\bf
Acknowledgments} \par\bigskip

We are grateful to Dennis Sullivan, Mikhail Lyubich, Krerley Oliveira, Carlos Matheus, and Fl\'avio Ferreira for a number
of very fruitful and useful discussions on this work and also
for their friendship and encouragement.
We thank the Calouste Gulbenkian Foundation, PRODYN-ESF, POCTI by FCT
and Minist\'erio da Ci\^encia e da Tecnologia, and Centro de Matem\'atica da
Universidade do Porto for their
financial support of A. A. Pinto
and the UK Engineering and Physical Sciences
Research Council and the EU for their
financial support of D. A. Rand.
Part of this research was
started during a visit by the
authors to the IHES and CUNY and
continued at the Isaac Newton
Institute, IMPA and SUNY.
We thank them for their  hospitality.

\end{document}